\documentclass[11pt]{amsart}
\usepackage[a4paper, margin=1in]{geometry}
\usepackage{macros}
\setlist[itemize]{leftmargin=1.5em,itemsep=3pt,topsep=4pt}
\allowdisplaybreaks
\usepackage{dirtytalk}
\usepackage{wasysym}

\title[Induced packing treewidth III]{Induced packing treewidth III.\\ Branch decompositions}

\author{Amir Nikabadi$^{\natural}$}\thanks{$^{\natural}$IT University of Copenhagen, Denmark (\texttt{amir@itu.dk})}
\author{Paweł Rzążewski$^{\parallel}$}\thanks{$^{\parallel}$Warsaw University of Technology, Poland (\texttt{pawel.rzazewski@pw.edu.pl}). Supported by the National Science Centre grant 2024/54/E/ST6/00094.}

\date{}

\begin{document}
\begin{titlepage}
\begin{abstract}

The notion of induced packing treewidth aims to unify classes defined by forbidden induced subgraphs or induced minors with classes defined by the existence of certain structured tree decompositions.

In this paper, we introduce induced packing branchwidth, a branch-decomposition analogue
of induced packing treewidth. For a graph $H$, the induced $H$-packing
branchwidth of a graph $G$, denoted by $\branchpi_H(G)$, is the minimum,
over all branch decompositions of $G$, of the largest number of pairwise
anticomplete induced copies of $H$ that cross one displayed cut.
For $H=P_2$, it equals sim-width.
We prove several structural properties of graphs of bounded $\branchpi_H$, for various choices of $H$. In particular, we prove the following results:

\begin{itemize}
    \item For every fixed graph $H$ that is a subdivided claw or a broom, and every integer $k\geq 0$, the class $\{G\mid\branchpi_H(G)\le k\}$ is $\chi$-bounded. These extend and answer a question of Oum [{BIRS} workshop on \say{New Perspectives in Coloring and Structure}, 2024] on $\chi$-boundedness of graphs of bounded sim-width.

\item For every fixed forest $H$ with at least two vertices and integer $k\geq 0$, graphs of induced $H$-packing branchwidth at most $k$ that exclude $K_{s,s}$ as a subgraph, where $s\ge2$, have degeneracy bounded by a polynomial in $s$. Consequently, for every fixed $s$, the corresponding class excluding an induced copy of $K_{s,s}$ is polynomially $\chi$-bounded. This extends a result of Scott, Seymour, and Spirkl [Journal of Graph Theory, 2023].

\item For every fixed $t\ge2$, bounded induced $P_t$-packing branchwidth
implies the existence of dominated weighted balanced separators, with
explicit bounds on the number of dominating vertices. As an application, we prove that the treewidth of any graph is linearly equivalent to the induced $P_t$-packing branchwidth of its line graph.
The case $t=2$ answers a question of Brettell, Munaro,
Paulusma, and Yang on the sim-width of line graphs [European
Journal of Combinatorics, 2025].

\item For every fixed $t\ge2$ and $k\ge 1$, graphs with
$\branchpi_{P_t}\le k$ are (uniformly) quasi-isometric to graphs of
treewidth at most $3k/2$. This extends
Hickingbotham's result for bounded sim-width.
\end{itemize}
\end{abstract}
\maketitle
\thispagestyle{empty}
\vspace{-0.5cm}
\hypersetup{linkcolor=red!50!black}
\tableofcontents

\paragraph{AI disclosure.}
The original ideas and results presented in this paper are the product of human authorship.
Claude Opus 5.5 was used to check the paper for mistakes, as well as to polish the presentation. The authors take full responsibility for the content of the paper.
\thispagestyle{empty}
\end{titlepage}

\newpage\setcounter{page}{1}

\section{Introduction}\label{sec:introduction}

\subsection{Motivation}
Graph decompositions provide a common setting for studying graph
structure and width parameters. In the first paper of this
series~\cite{nikabadi2026induced}, we introduced \emph{induced packing treewidth}.
For a family
$\cH$ of graphs, an induced $\cH$-packing is a collection of pairwise
anticomplete induced subgraphs, each isomorphic to a member of $\cH$.
The induced $\cH$-packing treewidth, denoted by $\treepi_{\cH}$, measures
the maximum number of members of such a packing that can intersect one
bag, minimized over tree decompositions. For $\cH=\{P_1\}$ and
$\cH=\{P_2\}$, this recovers tree-independence
number~\cite{dallard2024treewidthII,yolov2018minor} and induced matching
treewidth~\cite{yolov2018minor,DBLP:journals/jcss/LimaMMORS26}, respectively. The second paper~\cite{nikabadiII} studies structural
consequences of bounded induced packing treewidth and, in particular, its relationship
with sim-width.

We now consider the corresponding notion for branch decompositions.
A \emph{branch decomposition} of a graph $G$ is a tree in which every
internal node has degree three and whose leaves correspond bijectively
to the vertices of $G$. Deleting an edge $e$ of the tree splits its
leaves, and hence $V(G)$, into two sides; this partition is the
\emph{cut displayed} by $e$. Several width parameters, including
sim-width, minimize over all branch decompositions a measure of the
displayed cuts.
Instead of counting the members of an induced packing that meet a bag,
we count those that meet both sides of a displayed cut. The
\emph{induced $\cH$-packing branchwidth} of $G$, denoted
$\branchpi_{\cH}(G)$, is the minimum, over all branch decompositions of
$G$, of the largest size of an induced $\cH$-packing whose members all
cross one displayed cut. For a singleton family $\{H\}$, we write
$\branchpi_H(G)$ for $\branchpi_{\{H\}}(G)$.
For example, $\branchpi_{P_2}(G)$ is exactly
\emph{sim-width}~\cite{KangKwonStrommeTelle2017}, which we denote by
$\simw(G)$. If every member of $\cH$ is connected, then
$\branchpi_{\cH}(G)\le\min\{\simw(G),\treepi_{\cH}(G)\}$.
Thus this parameter provides a common generalization of sim-width and
induced packing treewidth. Its precise definition and basic properties
are given in~\cref{sec:prelim}.

\subsection{Our results}

\begin{figure}[b]
    \centering

    \begin{subfigure}[c]{2.4cm}
        \centering
        \begin{tikzpicture}[
            font=\fontsize{7.2}{8.6}\selectfont,
            line width=.60pt,
            line cap=round,
            line join=round,
            every node/.style={inner sep=2.2pt},
            edge/.style={draw=black}]

            \coordinate (mc-x) at (0,0);

            \foreach \i in {1,...,5} {
                \coordinate (mc-u\i) at (45:{0.30*\i});
            }
            \draw[edge]
                (mc-x)--(mc-u1)--(mc-u2)--(mc-u3)
                --(mc-u4)--(mc-u5);

            \foreach \i in {1,...,5} {
                \coordinate (mc-v\i) at (0:{0.30*\i});
            }
            \draw[edge]
                (mc-x)--(mc-v1)--(mc-v2)--(mc-v3)
                --(mc-v4)--(mc-v5);

            \foreach \i in {1,...,5} {
                \coordinate (mc-w\i) at (-45:{0.30*\i});
            }
            \draw[edge]
                (mc-x)--(mc-w1)--(mc-w2)--(mc-w3)
                --(mc-w4)--(mc-w5);

            \fill (mc-x) circle (1.08pt);

            \foreach \i in {1,...,5} {
                \fill (mc-u\i) circle (1.08pt);
                \fill (mc-v\i) circle (1.08pt);
                \fill (mc-w\i) circle (1.08pt);
            }

            \node[left] at (mc-x) {};
            \node[above right] at (mc-u5) {};
            \node[right] at (mc-v5) {};
            \node[below right] at (mc-w5) {};

        \end{tikzpicture}
    \end{subfigure}%
    \hspace{0.35cm}%
    \begin{subfigure}[c]{2.9cm}
        \centering
        \begin{tikzpicture}[
            font=\fontsize{7.2}{8.6}\selectfont,
            line width=.60pt,
            line cap=round,
            line join=round,
            every node/.style={inner sep=2.2pt},
            edge/.style={draw=black}]

            \foreach \i in {1,...,6} {
                \coordinate (br-v\i) at ({0.30*(\i-1)},0);
            }

            \draw[edge]
                (br-v1)--(br-v2)--(br-v3)
                --(br-v4)--(br-v5)--(br-v6);

            \path (br-v6) ++(60:.68) coordinate (br-u1);
            \path (br-v6) ++(20:.68) coordinate (br-u2);
            \path (br-v6) ++(-20:.68) coordinate (br-u3);
            \path (br-v6) ++(-60:.68) coordinate (br-u4);

            \foreach \i in {1,...,4} {
                \draw[edge] (br-v6)--(br-u\i);
            }

            \foreach \i in {1,...,6} {
                \fill (br-v\i) circle (1.08pt);
            }

            \foreach \i in {1,...,4} {
                \fill (br-u\i) circle (1.08pt);
            }

            \node[left] at (br-v1) {};
            \node[below left] at (br-v6) {};
            \node[above right] at (br-u1) {};
            \node[right] at (br-u2) {};
            \node[right] at (br-u3) {};
            \node[below right] at (br-u4) {};

        \end{tikzpicture}
    \end{subfigure}

    \caption{Left: The subdivided claw $S_5$. Right: The broom $B_{6,4}$.}
    \label{fig:subdivided-claw-broom}
\end{figure}
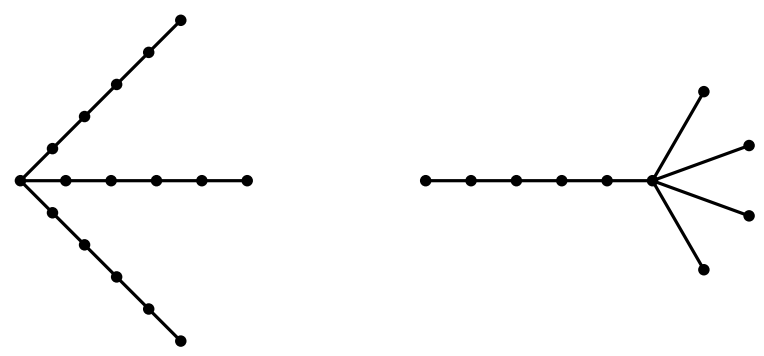

A class $\cC$ of graphs is \emph{$\chi$-bounded} if there is a function
$f$ such that every $G\in\cC$ satisfies $\chi(G)\le f(\omega(G))$,
where $\chi(G)$ and $\omega(G)$ denote the chromatic number and the
clique number of $G$; it is \emph{polynomially $\chi$-bounded} if $f$
can be chosen to be a polynomial. Our first main result shows that
bounded induced $H$-packing branchwidth implies $\chi$-boundedness
when $H$ is a subdivided claw or a broom.

For $t\ge1$, the \emph{subdivided claw} $S_t$ is the tree consisting of
three paths of length $t$ with a common endpoint and
otherwise disjoint vertex sets; see~\cref{fig:subdivided-claw-broom}.
For $t,d\ge1$, a \emph{broom} $B_{t,d}$ is obtained from a path on $t$ vertices, with endpoints $v_1$ and $v_t$, by adding $d$ new vertices, each adjacent to $v_t$.
The \emph{handle} of the broom is the displayed path, its \emph{center} is
$v_t$, and its \emph{root} is $v_1$; see~\cref{fig:subdivided-claw-broom}.
Brooms include all paths with at least two vertices and all stars.

\begin{restatable}{theorem}{subdividedClawsAndBrooms}\label{thm:mainchiboundedness}
Let $H$ be a subdivided claw or a broom. Then for every $k\in\Nzero$,
the class $\{G\mid\branchpi_H(G)\le k\}$ is $\chi$-bounded.
\end{restatable}

Since $P_2=B_{1,1}$ and, for every $G$, $\branchpi_{P_2}(G)=\simw(G)$,
\cref{thm:mainchiboundedness} implies in particular that graphs of
bounded sim-width are $\chi$-bounded. This answers a question that Oum
posed at the BIRS workshop ``New Perspectives in Coloring and
Structure'' in 2024~\cite{Oum}. The same question was answered
independently by Bria\'nski, Giocanti, Legrand-Duchesne, Micek, Oum,
and Walczak~\cite{simwchibdd}. \Cref{thm:mainchiboundedness} is more
general: it applies to every subdivided claw and every broom $H$, and
since $\branchpi_H(G)\le\simw(G)$ for every graph $G$ and every connected $H$,
each of these
classes contains all graphs of bounded sim-width. We prove
\cref{thm:mainchiboundedness} in~\cref{sec:chi-boundedness}.

The binding functions obtained in \cref{thm:mainchiboundedness} are very
large, since our proof, as well as the results it relies
on~\cite{scott1997induced,CSSLH,CSSChandeliers,abrishami2026burling},
applies Ramsey's theorem repeatedly. Our second result gives polynomial bounds for every forest
$H$, provided that a fixed \emph{biclique}, that is, a complete
bipartite graph $K_{s,s}$, is excluded as an induced subgraph.
The assumption that $H$ is a forest cannot be dropped, as we show at
the end of~\cref{sec:degreeboundedness}.

\begin{restatable}{theorem}{polydegree}\label{thm:forest-induced-biclique-chi}
For every $h\ge2$ and $k\in\Nzero$, there are a constant $c=c(h,k)$
and a positive integer $\ell=\ell(h,k)$ with the following property.
Let $H$ be any forest on $h$ vertices, let $s\ge2$ be an integer, and let $G$ be
a graph with $\branchpi_H(G)\le k$ and no induced copy of $K_{s,s}$. Then
\[
 \chi(G)\le\dgn(G)+1
 \le c\bigl(\omega(G)+1\bigr)^{(s-1)\ell}.
\]
\end{restatable}

Since every graph $G$ satisfies $\chi(G)\le\dgn(G)+1$, the polynomial
bound on the degeneracy in \cref{thm:forest-induced-biclique-chi}
immediately gives the following.

\begin{corollary}\label{cor:forest-induced-biclique-polychi}
For every fixed forest $H$ with at least two vertices
and fixed integers $k\ge0$ and $s\ge2$, the class of graphs with
$\branchpi_H\le k$ and no induced copy of $K_{s,s}$ is polynomially
$\chi$-bounded.
\end{corollary}

We prove \cref{thm:forest-induced-biclique-chi} in
\cref{sec:degreeboundedness}. Its main step is a bound, polynomial in
$s$, on the degeneracy of graphs with $\branchpi_H\le k$ and no
$K_{s,s}$ subgraph (\cref{prop:biclique-polynomial-degree}). For $k=0$, such
a bound was proved by Scott, Seymour, and Spirkl~\cite{SSS23}.

Our remaining results concern paths, that is, $H=P_t$ for some fixed $t$. We first show
that bounded $\branchpi_{P_t}$ yields balanced separators dominated by
few vertices: for every vertex-weighting $\wei$, the closed
neighborhood of a bounded number of vertices is a $1/2$-balanced
separator of $(G,\wei)$. Weighted graphs and balanced separators are
defined in~\cref{sec:prelim}.

\begin{theorem}\label{thm:dominated}
Let $t\ge2$, let $k\ge1$, and let $G$ satisfy $\branchpi_{P_t}(G)\le k$.
For every vertex-weighting $\wei$ of $G$, there is
$D\subseteq V(G)$ with $|D|\le t\lfloor3k/2\rfloor+t-2$ such that
$N[D]$ is a $1/2$-balanced separator of $(G,\wei)$.
\end{theorem}

The proof in~\cref{sec:separators} uses a weighted centroid of a branch
decomposition and a clean-partition version of the Gy\'arf\'as path
argument~\cite{gyarfas1987problems}. There we also prove a $2/3$-balanced version in which, in
addition, $G[D]$ has at most $k+1$ components. This is useful when an
induced star is excluded: the closed neighborhood of each component
has bounded independence number. Combined with the separator theorem
of Chudnovsky, Hajebi, Lokshtanov, and Spirkl~\cite[Lemma~7.1]{chudnovsky2026tree},
this gives the linear estimate used in our next result.

For a graph $G$, let $L(G)$ denote its line graph.
Brettell, Munaro, Paulusma, and Yang~\cite[Open Problem~4]{BMPY25} asked
for the asymptotically optimal function $f$ such that
$\simw(L(G))\ge f(\tw(G))$ for every graph $G$. We show that, for every
fixed $t\ge2$, each of $\branchpi_{P_t}(L(G))$ and $\tw(G)$ is bounded
by a linear function of the other. For $t=2$, this shows that the
optimal $f$ is linear, answering their question.

\begin{restatable}{theorem}{linegraphtheorem}\label{thm:linegraph}
Let $t\ge2$. Let $G$ be a graph and let $k=\branchpi_{P_t}(L(G))$. Then $k\le\tw(G)+1\le21((t+1)k+t-1)$. Moreover, if $t=2$ and $k\ge1$, then $\tw(G)+1\le63k$.
\end{restatable}

The proof of~\cref{thm:linegraph} appears in~\cref{sec:linegraphs}, and it bypasses the
Grid-Minor Theorem~\cite{gridminorthm}. We note that the dependence on both the width and $t$ has optimal order
(see~\cref{lem:sharpness}).

Finally, we study the large-scale geometry of graphs of bounded
$\branchpi_{P_t}$. A class of graphs has \emph{bounded coarse treewidth}
if there are $a$ and $R$ such that every graph in the class has a tree
decomposition in which every bag is covered by at most $a$ balls of
radius $R$. Roughly speaking, a \emph{quasi-isometry} between two graphs
is a map that distorts distances by at most a bounded multiplicative
factor and a bounded additive error, and whose image is within bounded
distance of every vertex of the target graph. A class is
\emph{uniformly quasi-isometric to graphs of bounded treewidth} if all
its members admit quasi-isometries with the same constants to graphs of
treewidth bounded by a common constant. Nguyen, Scott, and
Seymour~\cite{NSSCoarse} and Hickingbotham~\cite{Hick} showed that these
two properties are equivalent. We add a third equivalent condition, in
terms of a distance version of induced packing branchwidth: for
$r\ge1$, the parameter $\branchpi_{P_t,r}$ is defined like
$\branchpi_{P_t}$, except that the counted induced copies of $P_t$ must
have pairwise distance greater than $r$. Note that
$\branchpi_{P_t,1}=\branchpi_{P_t}$. Precise definitions are given
in~\cref{sec:metric}.

\begin{restatable}{theorem}{quasiso}\label{thm:qisometri}
	Let $t\ge2$ be fixed. For a graph class $\mathcal C$, the following
	are equivalent.
	\begin{enumerate}
		\item $\mathcal C$ has bounded coarse treewidth.
		\item $\mathcal C$ is uniformly quasi-isometric to graphs of bounded
		treewidth.
		\item There exist integers $r\ge1$ and $k\ge0$ such that
		$\branchpi_{P_t,r}(G)\le k$ for every $G\in\mathcal C$.
	\end{enumerate}
\end{restatable}

We prove \cref{thm:qisometri} in~\cref{sec:metric}, where we also give
explicit constants. For instance, if $\branchpi_{P_t}(G)\le k$ for some
$k\ge1$, then $G$ is quasi-isometric, with constants depending only on
$t$ and $k$, to a graph of treewidth at most $\lfloor3k/2\rfloor$
(\cref{thm:metric-quasi-isometry}). For $t=2$, since
$\branchpi_{P_2}=\simw$, this improves the bound $12k-1$ on the
treewidth in a result of Hickingbotham~\cite[Theorem~14]{Hick} for
graphs of sim-width at most $k$ to $\lfloor3k/2\rfloor$. For every
$t\ge3$, the hypothesis is strictly weaker than bounded sim-width, as
shown at the end of~\cref{sec:metric}, so this also extends his result.

\section{Preliminaries}\label{sec:prelim}

For a positive integer $n$, by $[n]$ we denote the set $\{1,\ldots,n\}$.
For a set $X$, by $2^X$ we denote the family of all subsets of $X$.

\paragraph{Graphs}
Throughout the paper, graphs have finite vertex sets, no loops, and no parallel edges.
Let $G$ be a~graph with vertex set $V(G)$ and edge set $E(G)$.
For $X \subseteq V(G)$, we denote the subgraph of $G$ \emph{induced by $X$} as $G[X]$,
that is, $G[X] = (X, \{uv \mid u, v \in X \mbox{ and } uv \in E(G) \})$,
and we write $G-X$ for $G[V(G) \setminus X]$.
For a graph $H$, an \emph{induced copy of $H$} in $G$ is an induced subgraph of $G$ isomorphic to $H$, and $G$ is \emph{$H$-free} if it contains no induced copy of $H$; for graphs $H_1,\ldots,H_m$, $G$ is \emph{$(H_1,\ldots,H_m)$-free} if it is $H_i$-free for every $i\in[m]$. For graphs $H_1$ and $H_2$, we write $H_1+H_2$ for the disjoint union of $H_1$ and $H_2$. The \emph{line graph} $L(G)$ of $G$ has vertex set $E(G)$, and two of its vertices are adjacent if the corresponding edges of $G$ share an endpoint.

By $N_G(v)$ we denote the set of neighbors of a vertex $v$ in $G$, and we use $N_G[v]$ for $N_G(v) \cup \{v\}$.
For a set $X \subseteq V(G)$, we define $N_G[X] = \bigcup_{v \in X} N_G[v]$ and $N_G(X) = N_G[X] \setminus X$.
A set $D\subseteq V(G)$ \emph{covers} a set $C\subseteq V(G)$ if $C\subseteq N_G(D)$.
Two disjoint sets $X,Y \subseteq V(G)$ are \emph{complete} (resp., \emph{anticomplete}) if all (resp., no) edges between them exist. Analogously, we say that a vertex is (anti)complete to a set.
We write $X\perp Y$ if $X$ and $Y$ are disjoint and anticomplete.
We also say, informally, that $X$ \emph{sees} $Y$ if some vertex of $X$ has a neighbor in $Y$, and that $X$ \emph{misses} $Y$ if $X\perp Y$. For a vertex $v$, we write $v$ instead of $\{v\}$.

The \emph{degree} of a vertex $v$ is $\deg_G(v)=|N_G(v)|$.
For a nonempty graph $G$, we denote by $\delta(G)$ and $\overline d(G)=2|E(G)|/|V(G)|$ its minimum and average degree, respectively, and its \emph{degeneracy} is
\[
 \dgn(G)=\max\{\delta(G[X])\mid\varnothing\ne X\subseteq V(G)\};
\]
we also let $\dgn(\varnothing)=-1$.

By $\alpha(G)$, $\omega(G)$, and $\chi(G)$ we denote, respectively, the sizes of a largest independent set and of a largest clique in $G$, and the chromatic number of $G$.
A class $\cC$ of graphs is \emph{$\chi$-bounded} if some function $f$, called a \emph{$\chi$-binding function} for $\cC$, satisfies $\chi(G)\le f(\omega(G))$ for every $G\in\cC$; it is \emph{polynomially $\chi$-bounded} if $f$ can be chosen to be a polynomial.
A nonempty graph is \emph{vertex-critical} if deleting any vertex decreases its chromatic number.
Every nonempty graph has a vertex-critical induced subgraph of the same chromatic number, and such a subgraph is connected and has minimum degree at least $\chi-1$.
Deleting successively a vertex of minimum degree and then coloring in reverse order gives $\chi(G)\le\dgn(G)+1$; counting each edge when its first end is deleted gives $|E(G)|\le\dgn(G)\,|V(G)|$, so $\overline d(G)\le2\dgn(G)$ if $G$ is nonempty.

When it does not lead to confusion, we will sometimes identify induced subgraphs with their vertex sets. For example, for $X \subseteq V(G)$, we write $\alpha(X)$ for $\alpha(G[X])$.
For a family $\mathcal Q$ of subgraphs of $G$, we write $V(\mathcal Q)=\bigcup_{Q\in\mathcal Q}V(Q)$.

The length of a path is the number of its edges, i.e., $P_t$ has length $t-1$.
The \emph{distance} $\dist_G(u,v)$ between vertices $u$ and $v$ is the length of a shortest $u$--$v$ path in $G$, or $\infty$ if there is no such path.
For $X,Y\subseteq V(G)$ and $\rho\ge0$, let
\[
 \dist_G(X,Y)=\min\{\dist_G(x,y)\mid x\in X,\ y\in Y\},\qquad
 B_G(X,\rho)=\{v\in V(G)\mid \dist_G(v,X)\le \rho\},
\]
where $\min\varnothing=\infty$. For a vertex $v$, we write $\dist_G(v,X)=\dist_G(\{v\},X)$ and $B_G(v,\rho)=B_G(\{v\},\rho)$.
The \emph{weak diameter} $\diam_G(X)$ of a nonempty set $X\subseteq V(G)$ is the maximum of $\dist_G(x,y)$ over all $x,y\in X$; note that distances are measured in $G$, not in $G[X]$.
The \emph{diameter} of $G$ is $\diam(G)=\diam_G(V(G))$.
Throughout the paper, we drop the subscript $G$ from all the notation above whenever this does not lead to confusion.
We always keep it in $\diam_G(X)$, which in general differs from $\diam(G[X])$.

A \emph{subdivision} of a graph $F$ is obtained from $F$ by replacing
each edge with a path of at least one edge, with mutually disjoint
interiors; for $q\ge1$, the subdivision in which all these paths have
$q$ edges is denoted by $F^{(q)}$. A subdivision of the star $K_{1,s}$,
where $s\ge1$, consists of $s$ paths, its \emph{arms}, with a common
end. In particular, $K_{1,s}^{(r)}$ has $s$ arms of length $r$, and
$K_{1,s}^{(1)}=K_{1,s}$.

\paragraph{Weighted graphs and balanced separators.}
A \emph{vertex-weighted graph} is a pair $(G,\wei)$ where $G$ is a graph and $\wei:V(G)\to\mathbb Q_{\ge0}$ is a weight function; we also call such a $\wei$ a \defn{vertex-weighting} of $G$. For a set $X\subseteq V(G)$, we write $\wei(X)=\sum_{v\in X}\wei(v)$, and we write $\wei(G):=\wei(V(G))$.

For a weighted graph $(G,\wei)$ and $\delta\in[1/2,1)$, a set $X\subseteq V(G)$ is a \emph{$\delta$-balanced separator} of $(G,\wei)$ if every component of $G-X$ has weight at most $\delta\wei(G)$.

A weighted tree has a \emph{centroid}, a node whose
deletion leaves components of weight at most half the total. One finds
such a node by repeatedly moving towards a component heavier than
half the total.

\paragraph{Tree decompositions}
For a graph $G$, a \emph{tree decomposition} of $G$ is a pair $(T,\beta)$, where $T$ is a tree and $\beta \colon V(T) \rightarrow 2^{V(G)}$ is a map with the following properties:
\begin{itemize}
    \item For every $v \in V(G)$, there exists $x \in V(T)$ such that $v \in \beta(x)$.
    \item For every $uv \in E(G)$, there exists $x \in V(T)$ such that $u, v \in \beta(x)$.
    \item For every $v \in V(G)$, the subgraph of $T$ induced by $\{x \in V(T) \mid v \in \beta(x)\}$ is connected.
\end{itemize}
For every $x \in V(T)$, we refer to $x$ as a \emph{node} of $T$, and to $\beta(x)$ as a \emph{bag} of $(T, \beta)$.
The \emph{width} of a tree decomposition $(T, \beta)$ is $\max_{x\in V(T)} |\beta(x)|-1$.
The \emph{treewidth} of $G$, denoted by $\tw(G)$, is the minimum width of a tree decomposition of $G$.
Every nonempty graph $G$ satisfies $\dgn(G)\le\tw(G)$. Indeed, for $\varnothing\ne X\subseteq V(G)$, take a tree decomposition of $G[X]$ of width at most $\tw(G)$ (e.g., restrict the bags of one for $G$) with the fewest nodes. If it has a single bag, every vertex has all its neighbors in that bag. Otherwise a leaf bag is not contained in its neighboring bag, so it contains a vertex lying in no other bag, and all neighbors of this vertex lie in the leaf bag. Hence $\delta(G[X])\le\tw(G)$.
The \emph{independence number} of a tree decomposition $(T, \beta)$ of $G$ is $\max_{x\in V(T)} \alpha(G[\beta(x)])$. The \emph{tree-independence number} of $G$, denoted $\treealpha(G)$, is the minimum independence number of a tree decomposition of $G$.

Every tree decomposition of a weighted graph $(G,\wei)$ has a bag that is a $1/2$-balanced separator of $(G,\wei)$ (see, e.g.,~\cite[Lemma~7.19]{cygan2015parameterized}).
Conversely, balanced separators of bounded independence number, found for
\emph{every} vertex-weighting, certify bounded tree-independence number. The
following is due to~\cite[Lemma~7.1]{chudnovsky2026tree}. Note that this result is stated for normalized vertex-weightings, i.e., those with
$\wei(G)=1$.

\begin{lemma}[Chudnovsky, Hajebi, Lokshtanov, and Spirkl~\cite{chudnovsky2026tree}]\label{lem:sep-vs-ta}
Let $G$ be a graph, let $\delta \in [1/2 , 1)$, and let $m$ be a positive integer. If for every normalized
vertex-weighting $\wei$ of $G$, there is a $\delta$-balanced separator $Z_{\wei}$ of $(G,\wei)$ with $\alpha(Z_{\wei}) \leq m$, then $\treealpha(G) \leq \frac{3-\delta}{1-\delta}m$.
\end{lemma}

\paragraph{Induced packings}
Induced packings and induced packing treewidth were defined in~\cite{nikabadi2026induced}.
Let $G$ be a graph, let $\cH$ be a family of graphs, and let $B \subseteq V(G)$.
An \emph{induced $\cH$-packing} in $G$ is a set of pairwise anticomplete induced subgraphs of $G$, each isomorphic to a member of $\cH$.
An induced subgraph of $G$ is \emph{$B$-rooted} if its vertex set intersects $B$, and an induced $\cH$-packing is \emph{$B$-rooted} if all its members are $B$-rooted.

The size of a largest $B$-rooted induced $\cH$-packing in $G$ is denoted by $\pi_{\cH}(G,B)$.
For a tree decomposition $\cT=(T,\beta)$ of $G$, we define $\pi_{\cH}(G,\cT) = \max_{x \in V(T)} \pi_{\cH}(G, \beta(x))$.
The \emph{induced $\cH$-packing treewidth} of $G$ is $\treepi_{\cH}(G) = \min_{\cT} \pi_{\cH}(G,\cT)$, where the minimum is taken over all tree decompositions $\cT$ of $G$.

\paragraph{Branch decompositions}
Branch decompositions provide a common framework for several width
parameters that measure the complexity of cuts in a graph.
A \emph{cut} of $G$ is a pair $(A,B)$ with $A\subseteq V(G)$ and
$B=V(G)\setminus A$; the sets $A$ and $B$ are its
\emph{sides}, and the cut is \emph{nontrivial} if both sides are nonempty.
A set of vertices, or a subgraph, of $G$ \emph{crosses} a cut if it meets both sides.
Let $G$ be a graph with at least two vertices. A \defn{branch decomposition} of $G$ is a
pair $\mathcal D=(T,\lambda)$, where $T$ is a subcubic tree in which
every internal node has degree three, and $\lambda$ is a
bijection from $V(G)$ to the set of leaves of $T$. Every edge $e\in E(T)$
\emph{displays} the cut $(A_e,B_e)$ of $G$ such that the two components
of $T-e$ contain respectively the leaves in $\lambda(A_e)$ and
$\lambda(B_e)$. The names of the two sides may be
interchanged.
Suppressing a node of degree two in a tree whose leaves are labelled
bijectively by $V(G)$ does not change the cuts of $G$ determined by its
edges in this way. Hence allowing such nodes, as some authors do,
changes none of the width parameters defined below.
For $X\subseteq V(G)$ with $|X|\ge2$, the \emph{restriction} of
$\mathcal D$ to $X$ is the branch decomposition $(T_X,\lambda_X)$ of
$G[X]$, where $T_X$ is obtained from the minimal subtree of $T$
containing $\lambda(X)$ by suppressing its nodes of degree two, and
$\lambda_X$ is the restriction of $\lambda$ to $X$.
Every edge of $T_X$ replaces a path of $T$, and it displays the cut
$(A_e\cap X,B_e\cap X)$ of $G[X]$ for every edge $e$ of that path.
Moreover, for every nonempty $U\subseteq X$, a node of $T_X$ lies in
the minimal subtree of $T_X$ containing $\lambda(U)$ if and only if it
lies in the minimal subtree of $T$ containing $\lambda(U)$.

Just as some bag of a tree decomposition is a $1/2$-balanced
separator, every branch decomposition displays a cut that splits the
weight in a balanced way.

\begin{lemma}\label{lem:balanced-core}
Let $(G,\wei)$ be a weighted graph with at least two vertices such that
$\wei(v)\le2\wei(G)/3$ for every $v\in V(G)$, and let $(T,\lambda)$ be a
branch decomposition of $G$. Then some displayed cut $(A,B)$ satisfies
$\wei(A),\wei(B)\ge\wei(G)/3$.
In particular, for every $W\subseteq V(G)$ with $|W|\ge2$, some
displayed cut $(A,B)$ satisfies $|A\cap W|,|B\cap W|\ge\lceil|W|/3\rceil$.
\end{lemma}
\begin{proof}
Give each leaf $\lambda(v)$ the weight $\wei(v)$ and each internal node
weight zero, and let $x$ be a centroid of the resulting weighted tree.
If $x$ is internal, then the at most three components of $T-x$ have
weights summing to $\wei(G)$, each at most $\wei(G)/2$. One of them has
weight in $[\wei(G)/3,\wei(G)/2]$, and the edge joining it to $x$
displays the required cut. If $x=\lambda(v)$ is a leaf, then the only
component of $T-x$ has weight $\wei(G)-\wei(v)\le\wei(G)/2$, so
$\wei(G)/2\le\wei(v)\le2\wei(G)/3$, and the edge of $T$ at $x$ displays
the required cut $(\{v\},V(G)\setminus\{v\})$.
The second statement follows by applying the first one to the weight
that is $1$ on $W$ and $0$ elsewhere, which satisfies the hypothesis
since $|W|\ge2$.
\end{proof}

Conversely, a tree decomposition yields a branch decomposition in
which every displayed cut is separated by a bag. We use this
in~\cref{lem:basic}, \cref{sec:linegraphs}, and~\cref{sec:metric}.

\begin{lemma}\label{lem:tree-to-branch}
Let $G$ be a graph with at least two vertices, let $T$ be a tree, and
let $\eta\colon V(G)\to V(T)$. Then $G$ has a branch decomposition in
which every displayed cut $(A,B)$ has a node $x\in V(T)$ that lies on
the path of $T$ between $\eta(u)$ and $\eta(w)$ for all $u\in A$ and
$w\in B$.
In particular, if $(T,\beta)$ is a tree decomposition of $G$ and
$v\in\beta(\eta(v))$ for every $v\in V(G)$, then every connected
subgraph of $G$ crossing $(A,B)$ meets $\beta(x)$.
\end{lemma}
\begin{proof}
For each $v\in V(G)$, attach a new leaf $\ell_v$ to $\eta(v)$. Then
replace every node $y\in V(T)$ that now has degree $d>3$ by a path
$y_1\cdots y_d$ of new nodes, and join the $d$ former neighbors of $y$
to distinct nodes of this path. We call $y_1,\ldots,y_d$ the
\emph{copies} of $y$; a node of $T$ that is not replaced is its own
copy. The resulting tree $T'$ has maximum degree at most three, and its
leaves include all $\ell_v$. Let $(T'',\lambda)$ be obtained from the
minimal subtree of $T'$ containing all $\ell_v$ by suppressing its
nodes of degree two, with $\lambda(v)=\ell_v$. Since $|V(G)|\ge2$, this
is a branch decomposition of $G$.

Every edge of $T''$ replaces a path of $T'$ and displays the same cut
$(A,B)$ as every edge of that path. Let $f$ be an edge of this path.
No two leaves $\ell_v$ are adjacent, so $f$ has an end $y'$ that is a
copy of some node $x\in V(T)$. Let $u\in A$ and $w\in B$. The path $R$
of $T'$ between $\ell_u$ and $\ell_w$ contains $f$, and hence $y'$.
The copies of each node of $T$ induce a subpath of $T'$, so $R$ meets
them in a subpath. Hence contracting the copies maps $R$ onto the path
between $\ell_u$ and $\ell_w$ in the tree obtained from $T$ by
attaching the leaves. This path consists of $\ell_u$, the path of $T$
between $\eta(u)$ and $\eta(w)$, and $\ell_w$. Since $x\in V(T)$, it
lies on the path of $T$ between $\eta(u)$ and $\eta(w)$.

For the second statement, a connected subgraph of $G$ crossing $(A,B)$
contains an edge $uw$ with $u\in A$ and $w\in B$. The nodes whose bags
contain $u$ induce a subtree of $T$ containing $\eta(u)$, and those
whose bags contain $w$ induce a subtree containing $\eta(w)$. These
subtrees share a node whose bag contains both $u$ and $w$. Their union
is therefore a subtree containing the path between $\eta(u)$ and
$\eta(w)$, and hence $x$. Thus $u$ or $w$ belongs to $\beta(x)$.
\end{proof}

\paragraph{Maximum matching width and sim-width}
A matching $M$ of $G$ is \emph{induced} if every edge of $G$ with both
ends in $V(M)$ belongs to $M$.
The \emph{maximum matching width} $\mmw(G)$ and the
\emph{sim-width}~\cite{KangKwonStrommeTelle2017} $\simw(G)$ of a graph
$G$ with at least two vertices are the minima, over all branch
decompositions $(T,\lambda)$ of $G$, of the largest size of a matching
and of an induced matching, respectively, of $G$ whose edges all cross
one cut $(A_e,B_e)$ with $e\in E(T)$. For graphs with at most one vertex,
we set $\mmw(G)=\simw(G)=0$. We use the following known
inequality~\cite[Section~1 and Theorem~3.5]{JST}.

\begin{lemma}\label{lem:mmw-tw}
Every graph $G$ with at least one edge satisfies
$\tw(G)+1\le3\mmw(G)$.
\end{lemma}

\paragraph{Clean cuts}
For a graph $H$, a cut of $G$ is \emph{$H$-clean} if no induced
copy of $H$ in $G$ crosses it.
The following observation about $P_t$-clean cuts is used
in~\cref{sec:metric}.

\begin{lemma}\label{lem:clean-cut}
Let $F$ be connected, and let $(A,B)$ be a nontrivial cut of $F$.
If $(A,B)$ is $P_t$-clean, where $t\ge3$, then $\diam(F)\le2t-4$.
Moreover, if $v$ has a neighbor on the opposite side, then
$\dist(x,v)\le t-2$ for every $x\in V(F)$.
\end{lemma}

\begin{proof}
A shortest path between opposite sides is induced. If it had at least
$t$ vertices, a $t$-vertex subpath containing a crossing edge would
contradict cleanliness. Thus opposite-side vertices have distance at
most $t-2$. Joining same-side vertices through a vertex on the other
side gives the diameter bound.

For the stronger assertion, let $vw$ be a crossing edge and suppose
that $x$ is on the same side as $v$ with $\dist(x,v)\ge t-1$.
The preceding bound and the edge $wv$ imply $\dist(x,w)=t-2$.
Appending $v$ to a shortest $x$--$w$ path gives an induced copy of $P_t$:
an edge from $v$ to an earlier vertex would contradict
$\dist(x,v)\ge t-1$. This path crosses the cut, a contradiction.
\end{proof}

\paragraph{Induced packing branchwidth}
In branch decompositions, displayed cuts play the role of bags, and
subgraphs crossing a cut play the role of rooted subgraphs.
Let $\cH$ be a family of graphs.
The size of a largest induced $\cH$-packing in $G$ whose members all cross a cut $(A,B)$ is denoted by $\pi^{\times}_{\cH}(G;A,B)$.
The \emph{induced $\cH$-packing branchwidth} of a branch decomposition
$\cD=(T,\lambda)$ is
\[
  \branchpi_{\cH}(G,\cD)
  =\max_{e\in E(T)}\pi^{\times}_{\cH}(G;A_e,B_e),
\]
and the \emph{induced $\cH$-packing branchwidth} of $G$ is
\[
  \branchpi_{\cH}(G)
  =\min_{\mathcal{D}}\branchpi_{\cH}(G, \mathcal{D}),
\]
where the minimum is taken over all branch decompositions $\cD$ of $G$.
For graphs with at most one vertex, we set $\branchpi_{\cH}(G)=0$.

For a single graph $H$, we write $H$ instead of $\{H\}$ in all of the
notation above; for example, an \emph{induced $H$-packing} is a set of
pairwise anticomplete induced copies of $H$, and
$\branchpi_H=\branchpi_{\{H\}}$. In this paper, we mostly focus on
this setting.
An induced $P_2$-packing is an induced matching, and an induced copy of
$P_2$ crosses a cut exactly when its edge does. Hence
$\branchpi_{P_2}=\simw$.

We also collect the following basic properties of induced packing
branchwidth.

\begin{lemma}\label{lem:basic}
Let $H$ be a graph, and let $\cH$ be a family of graphs. The following
statements hold.
\begin{enumerate}[label=(\arabic*)]
\item\label{item:zero}
If $|V(H)|\ge2$, then $\branchpi_H(G)=0$ if and only if $G$ is $H$-free.
Also, $\branchpi_{P_{t+1}}(G)\le\branchpi_{P_t}(G)$ for $t\ge2$.
\item\label{item:heredity}
If $G'$ is an induced subgraph of $G$, then
$\branchpi_{\cH}(G')\le\branchpi_{\cH}(G)$.
\item\label{item:comparison}
If every member of $\cH$ is connected, then
$\branchpi_{\cH}(G)\le\min\{\simw(G),\treepi_{\cH}(G)\}$.
\end{enumerate}
\end{lemma}

The proofs are straightforward. For the bound
$\branchpi_{\cH}(G)\le\treepi_{\cH}(G)$ in~\ref{item:comparison}, apply
\cref{lem:tree-to-branch} to a tree decomposition $(T,\beta)$ of $G$
with $\pi_{\cH}(G,(T,\beta))=\treepi_{\cH}(G)$, choosing $\eta(v)$ with
$v\in\beta(\eta(v))$. The members of an induced $\cH$-packing crossing a
displayed cut are connected, so they all meet the corresponding bag
$\beta(x)$, and the packing is $\beta(x)$-rooted.

\section{\texorpdfstring{$\chi$}{χ}-boundedness}
\label{sec:chi-boundedness}

In this section, we prove \cref{thm:mainchiboundedness}. We reduce it
to our main technical result, \cref{thm:common-two-control}, which
concerns local control and is stated in \cref{subsec:local-control}.
After recalling Burling graphs in \cref{subsec:burling}, we derive
\cref{thm:mainchiboundedness} from \cref{thm:common-two-control} in
\cref{subsec:control-to-coloring}. The rest of the section proves
\cref{thm:common-two-control}.

All radius, width, and coloring bounds are integers.

\subsection{Local control and node width}
\label{subsec:local-control}

For $\rho\ge1$, let
\[
 \chi^\rho(G)=\max_{v\in V(G)}\chi(G[B(v,\rho)]),
\]
with value zero for the empty graph. A hereditary class $\cC$ is
\emph{$\rho$-controlled} if a nondecreasing function $\varphi$ satisfies
$\chi(G)\le\varphi(\chi^\rho(G))$ for every $G\in\cC$.
All control and $\chi$-binding functions are integer-valued and nondecreasing,
with value zero at zero and value at least $\tau$ at every $\tau\ge1$.
This is no restriction: a witness remains a witness when its values at
positive arguments are increased, and its value at zero matters only
for the empty graph, the only graph with $\chi^\rho(G)=0$ or
$\omega(G)=0$. Since $\chi^\rho\le\chi^{\rho'}$ for
$\rho\le\rho'$, every $\rho$-controlled class is $\rho'$-controlled.

For a branch decomposition $\mathcal D=(T,\lambda)$ of $G$ and nonempty
$U\subseteq V(G)$, the \emph{span} $T[U]$ of $U$ is the minimal subtree
of $T$ containing $\lambda(U)$; for a subgraph $Q$ of $G$, write
$T[Q]=T[V(Q)]$. Let $\sigma_H(G,\mathcal D)$ be the maximum number of
pairwise anticomplete induced copies of $H$ in $G$ whose spans contain a
common node of $T$. The \emph{node width} of $G$ is
\[
 \sigma_H(G)=\min_{\mathcal D}\sigma_H(G,\mathcal D),
\]
where the minimum is taken over all branch decompositions $\mathcal D$
of $G$; for graphs with at most one vertex, we set $\sigma_H(G)=0$.
A branch decomposition $\mathcal D$ of $G$ with
$\sigma_H(G,\mathcal D)\le k$ is a \emph{node-width-$k$ decomposition};
here $H$ is always the graph fixed in the relevant statement. Let
\[
 \cC_k(H)=\{G\mid\sigma_H(G)\le k\}.
\]

Node width approximates induced packing branchwidth, as the following
lemma shows, but it is easier to work with.

\begin{lemma}\label{lem:node-span-width}
Let $H$ have at least two vertices. Then
\[
 \branchpi_H(G)\le\sigma_H(G)
 \le\left\lfloor\frac{3\branchpi_H(G)}2\right\rfloor.
\]
Moreover, for every $k\ge0$, every cut $(A,B)$ displayed by a
node-width-$k$ decomposition of $G$ satisfies $\pi^\times_H(G;A,B)\le k$, and
$\sigma_H$ does not increase under taking induced subgraphs.
\end{lemma}
\begin{proof}
If $|V(G)|\le1$, then all quantities in the two inequalities are zero.
Otherwise, let $\mathcal D=(T,\lambda)$ be a branch decomposition of
$G$, and let $e=xy\in E(T)$. An induced copy $Q$ of $H$ crosses the
cut of $e$ precisely when $e\in E(T[Q])$, and then $x\in V(T[Q])$.
Hence $\pi^\times_H(G;A_e,B_e)\le\sigma_H(G,\mathcal D)$, which gives
the bound for displayed cuts and the first inequality.

For the second inequality, let $k=\branchpi_H(G)$, and let
$\mathcal D=(T,\lambda)$ satisfy $\branchpi_H(G,\mathcal D)=k$. If $k=0$,
then $G$ is $H$-free by \cref{lem:basic}\ref{item:zero}, so
$\sigma_H(G,\mathcal D)=0$. Otherwise,
for a packing $\mathcal Q$ counted at an internal node $x$,
\[
 2|\mathcal Q|
 \le\sum_{e\in E_T(x)}|\{Q\in\mathcal Q\mid e\in E(T[Q])\}|
 \le\deg_T(x)k\le3k,
\]
since each span through $x$ uses at least two edges incident with $x$.
At a leaf, every counted copy contains its labelled vertex, so the
packing has size at most one. Thus $\sigma_H(G,\mathcal D)\le\lfloor3k/2\rfloor$,
which proves the second inequality.

Finally, let $X\subseteq V(G)$ have at least two vertices, and let
$\mathcal D'$ be a branch decomposition of $G$ with
$\sigma_H(G,\mathcal D')=\sigma_H(G)$. In the restriction of
$\mathcal D'$ to $X$, a node lies in the span of an induced copy of
$H$ in $G[X]$ if and only if it lies in the span of that copy in
$\mathcal D'$ (see~\cref{sec:prelim}). Hence
$\sigma_H(G[X])\le\sigma_H(G)$; sets $X$ with at most one vertex are
immediate.
\end{proof}

We often use the bound $\pi^\times_H(G;A,B)\le k$ for cuts displayed
by a node-width-$k$ decomposition without further reference.

The main technical result of this section is the following theorem.

\begin{theorem}[Common control theorem]\label{thm:common-two-control}
For every $k\ge0$, the class $\cC_k(F)$ is $2$-controlled whenever
$F$ belongs to one of the following families:
\begin{enumerate}
 \renewcommand{\theenumi}{\roman{enumi}}
 \renewcommand{\labelenumi}{(\theenumi)}
 \item\label{item:common-paths} $F=P_t$, where $t\ge2$;
 \item\label{item:common-stars} $F=K_{1,s}$, where $s\ge3$;
 \item\label{item:common-long-arms} $F=K_{1,s}^{(r)}$, where $s\ge3$ and $r\ge3$;
 \item\label{item:common-short-claw} $F=K_{1,3}^{(2)}$;
 \item\label{item:common-short-broom} $F=B_{3,d}$, where $d\ge2$.
\end{enumerate}
\end{theorem}

\subsection{Burling graphs}
\label{subsec:burling}

We now recall Burling graphs and show that, for every tree $H$ with at
least two vertices, there are Burling graphs of arbitrarily large induced $H$-packing branchwidth
(\cref{lem:tree-burling-obstruction}). In
\cref{subsec:control-to-coloring}, we use them to turn $2$-control into
$\chi$-boundedness.

Let $(\mathsf{Bur}_j)_{j\ge1}$ be Burling's sequence~\cite{burling1965coloring} of triangle-free graphs with
$\chi(\mathsf{Bur}_j)=j$; see~\cite[Section~2.1]{abrishami2026burling} for a
definition. A graph is a \emph{Burling graph} if it is an induced
subgraph of some $\mathsf{Bur}_j$.

We use three results of Pournajafi and
Trotignon~\cite{pournajafi2023burling,pournajafi2024burling} on
Burling graphs and one of Abrishami et
al.~\cite{abrishami2026burling}, restated in our notation.
A \emph{Burling tree} is a rooted tree with a designated last-born child
at each nonleaf node and a choice function $\gamma$. For a nonroot vertex
$v$ that is not last-born, $\gamma(v)$ is empty or is the vertex set of
a downward path starting at the last-born child of the parent of $v$. The root and
all last-born vertices have empty choices. For any retained set $X$ of
tree vertices, the \emph{derived digraph} on $X$ has an arc $v\to w$
precisely when $w\in\gamma(v)\cap X$.

\begin{lemma}[Pournajafi and Trotignon~{\cite[Theorem~4.9]{pournajafi2023burling}}]
\label{lem:pt-derived}
A graph is Burling if and only if it is the underlying graph of a
derived digraph.
\end{lemma}

\begin{lemma}[Pournajafi and Trotignon~{\cite[Lemma~3.5]{pournajafi2023burling}}]
\label{lem:pt-normalization}
Every derived digraph has a Burling-tree representation in which the
root and all last-born vertices are omitted.
\end{lemma}

Fix a Burling-tree representation of a derived digraph $D$. An arc
$u\to v$ is a \emph{top-arc}, respectively a \emph{bottom-arc}, if $v$
is the out-neighbor of $u$ nearest to, respectively furthest from,
the root in that tree. These notions refer to the fixed representation.

\begin{lemma}[Pournajafi and Trotignon~{\cite[Lemma~3.8]{pournajafi2024burling}}]
\label{lem:pt-permitted-subdivision}
Let $D$ be a digraph derived from a Burling tree. Choose disjoint sets
$E_b,E_t$ of its bottom-arcs and top-arcs, respectively, such that
the tail of each arc in $E_t$ is a source of $D$.
The following simultaneous replacements produce a derived digraph:
\begin{enumerate}
\item replace each $u\to v$ in $E_b$ by a directed $u$--$v$ path of
any prescribed positive length;
\item replace each $u\to v$ in $E_t$ by an arc $w\to v$ and a
directed $w$--$u$ path of any prescribed positive length.
\end{enumerate}
All added vertices are private to their replacement, and all other
arcs are unchanged.
\end{lemma}

\begin{lemma}[Abrishami et al.~{\cite[Theorem~2.9]{abrishami2026burling}}]
\label{lem:2-controlled-Kmm}
For all $m,j\ge1$, the $\mathsf{Bur}_j$-free graphs in a
$2$-controlled class excluding $K_{m,m}^{(2)}$ form a $\chi$-bounded
class.
\end{lemma}

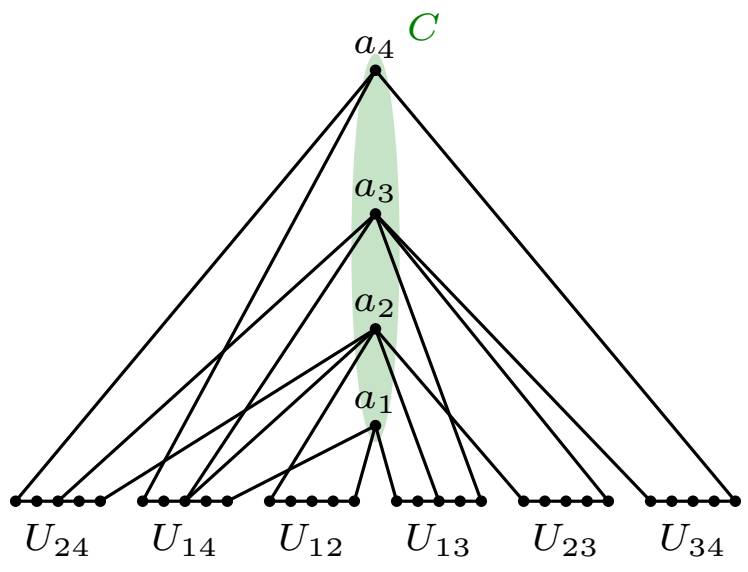
\begin{figure}[t]
\centering
\providecolor{graphgreen}{RGB}{0,128,0}
\begin{tikzpicture}[font=\fontsize{7.2}{8.6}\selectfont,
  line width=.60pt, every node/.style={inner sep=1.2pt}]
  \fill[graphgreen,fill opacity=.22,draw=none]
    (2.38,1.675) ellipse [x radius=.16cm,y radius=1.285cm];
  \coordinate (a1) at (2.38,0.5);
  \coordinate (a2) at (2.38,1.14);
  \coordinate (a3) at (2.38,1.9);
  \coordinate (a4) at (2.38,2.85);
  \coordinate (p241) at (0.56,0);
  \coordinate (p242) at (0.42,0);
  \coordinate (c24) at (0.28,0);
  \coordinate (q242) at (0.14,0);
  \coordinate (q241) at (0,0);
  \coordinate (p141) at (1.4,0);
  \coordinate (p142) at (1.26,0);
  \coordinate (c14) at (1.12,0);
  \coordinate (q142) at (0.98,0);
  \coordinate (q141) at (0.84,0);
  \coordinate (p121) at (2.24,0);
  \coordinate (p122) at (2.1,0);
  \coordinate (c12) at (1.96,0);
  \coordinate (q122) at (1.82,0);
  \coordinate (q121) at (1.68,0);
  \coordinate (p131) at (2.52,0);
  \coordinate (p132) at (2.66,0);
  \coordinate (c13) at (2.8,0);
  \coordinate (q132) at (2.94,0);
  \coordinate (q131) at (3.08,0);
  \coordinate (p231) at (3.36,0);
  \coordinate (p232) at (3.5,0);
  \coordinate (c23) at (3.64,0);
  \coordinate (q232) at (3.78,0);
  \coordinate (q231) at (3.92,0);
  \coordinate (p341) at (4.2,0);
  \coordinate (p342) at (4.34,0);
  \coordinate (c34) at (4.48,0);
  \coordinate (q342) at (4.62,0);
  \coordinate (q341) at (4.76,0);
  \draw (a1) -- (p121);
  \draw (a1) -- (p131);
  \draw (a1) -- (p141);
  \draw (a2) -- (c13);
  \draw (a2) -- (c14);
  \draw (a2) -- (p231);
  \draw (a2) -- (p241);
  \draw (a2) -- (q121);
  \draw (a3) -- (c14);
  \draw (a3) -- (c24);
  \draw (a3) -- (p341);
  \draw (a3) -- (q131);
  \draw (a3) -- (q231);
  \draw (a4) -- (q141);
  \draw (a4) -- (q241);
  \draw (a4) -- (q341);
  \draw (c12) -- (p122);
  \draw (c12) -- (q122);
  \draw (c13) -- (p132);
  \draw (c13) -- (q132);
  \draw (c14) -- (p142);
  \draw (c14) -- (q142);
  \draw (c23) -- (p232);
  \draw (c23) -- (q232);
  \draw (c24) -- (p242);
  \draw (c24) -- (q242);
  \draw (c34) -- (p342);
  \draw (c34) -- (q342);
  \draw (p121) -- (p122);
  \draw (p131) -- (p132);
  \draw (p141) -- (p142);
  \draw (p231) -- (p232);
  \draw (p241) -- (p242);
  \draw (p341) -- (p342);
  \draw (q121) -- (q122);
  \draw (q131) -- (q132);
  \draw (q141) -- (q142);
  \draw (q231) -- (q232);
  \draw (q241) -- (q242);
  \draw (q341) -- (q342);
  \fill[black] (a1) circle[radius=1.08pt];
  \fill[black] (a2) circle[radius=1.08pt];
  \fill[black] (a3) circle[radius=1.08pt];
  \fill[black] (a4) circle[radius=1.08pt];
  \fill[black] (p241) circle[radius=1.08pt];
  \fill[black] (p242) circle[radius=1.08pt];
  \fill[black] (c24) circle[radius=1.08pt];
  \fill[black] (q242) circle[radius=1.08pt];
  \fill[black] (q241) circle[radius=1.08pt];
  \fill[black] (p141) circle[radius=1.08pt];
  \fill[black] (p142) circle[radius=1.08pt];
  \fill[black] (c14) circle[radius=1.08pt];
  \fill[black] (q142) circle[radius=1.08pt];
  \fill[black] (q141) circle[radius=1.08pt];
  \fill[black] (p121) circle[radius=1.08pt];
  \fill[black] (p122) circle[radius=1.08pt];
  \fill[black] (c12) circle[radius=1.08pt];
  \fill[black] (q122) circle[radius=1.08pt];
  \fill[black] (q121) circle[radius=1.08pt];
  \fill[black] (p131) circle[radius=1.08pt];
  \fill[black] (p132) circle[radius=1.08pt];
  \fill[black] (c13) circle[radius=1.08pt];
  \fill[black] (q132) circle[radius=1.08pt];
  \fill[black] (q131) circle[radius=1.08pt];
  \fill[black] (p231) circle[radius=1.08pt];
  \fill[black] (p232) circle[radius=1.08pt];
  \fill[black] (c23) circle[radius=1.08pt];
  \fill[black] (q232) circle[radius=1.08pt];
  \fill[black] (q231) circle[radius=1.08pt];
  \fill[black] (p341) circle[radius=1.08pt];
  \fill[black] (p342) circle[radius=1.08pt];
  \fill[black] (c34) circle[radius=1.08pt];
  \fill[black] (q342) circle[radius=1.08pt];
  \fill[black] (q341) circle[radius=1.08pt];
  \node[above,inner sep=2pt] at (a1) {$a_1$};
  \node[above,inner sep=2pt] at (a2) {$a_2$};
  \node[above,inner sep=2pt] at (a3) {$a_3$};
  \node[above,inner sep=2pt] at (a4) {$a_4$};
  \node[text=graphgreen,above] at (2.7,3) {$C$};
  \node[below,inner sep=4pt] at (0.28,0) {$U_{24}$};
  \node[below,inner sep=4pt] at (1.12,0) {$U_{14}$};
  \node[below,inner sep=4pt] at (1.96,0) {$U_{12}$};
  \node[below,inner sep=4pt] at (2.8,0) {$U_{13}$};
  \node[below,inner sep=4pt] at (3.64,0) {$U_{23}$};
  \node[below,inner sep=4pt] at (4.48,0) {$U_{34}$};
\end{tikzpicture}
\caption{The graph $J_4$. The core $C$ is shaded.}
\label{fig:burling-J}
\end{figure}

We now construct the Burling graphs we need.
For an integer $m\ge2$, define $J_m$ as follows.
Start with an independent set $C=\{a_1,\ldots,a_m\}$, called the
\emph{core}. For each $1\le i<j\le m$, add a new central vertex $c_{ij}$,
and join $a_i$ to $c_{ij}$ and $c_{ij}$ to $a_j$ by two internally
disjoint paths of length three. All new vertices are private to the pair
$\{i,j\}$. Add the edges $c_{ij}a_h$ for $i<h<j$, and no other edges.
For distinct core vertices $u,v$, let $Q_{uv}$ be their specified private
path of length six, and let
$U_{uv}=V(Q_{uv})\setminus\{u,v\}.$
Different sets $U_{uv}$ are pairwise anticomplete, and every noncentral
vertex in such a set has degree two.
The example $J_4$ is shown in \cref{fig:burling-J}.

\begin{lemma}\label{lem:burling-construction}
For every $m\ge2$, the graph $J_m$ is a Burling graph.
\end{lemma}

\begin{figure}[t]
\centering
\definecolor{graphgreen}{RGB}{0,128,0}
\begin{tikzpicture}[
  font=\fontsize{7.2}{8.6}\selectfont,
  line width=.60pt,line cap=round,line join=round,
  every node/.style={inner sep=2.2pt},
  edge/.style={draw=black},
  arc/.style={-{Stealth[length=2.6pt,width=2.4pt]},shorten <=1.7pt,shorten >=1.7pt},
  region/.style={fill=graphgreen,fill opacity=.22,draw=none}]
\begin{scope}
  \path[region] (1.8,.975) ellipse[x radius=.12cm,y radius=1.13cm];
  \coordinate (a0) at (1.8,2.6);
  \coordinate (a1) at (1.8,1.95);
  \coordinate (a2) at (1.8,1.3);
  \coordinate (a3) at (1.8,.65);
  \coordinate (a4) at (1.8,0);
  \coordinate (c12) at (0,1.95);
  \coordinate (c13) at (.6,1.95);
  \coordinate (c14) at (1.2,1.95);
  \coordinate (c23) at (.6,1.3);
  \coordinate (c24) at (1.2,1.3);
  \coordinate (c34) at (1.2,.65);
  \draw[edge] (a0)--(a1)--(a2)--(a3)--(a4);
  \draw[edge] (a0)--(c12) (a0)--(c13) (a0)--(c14);
  \draw[edge] (a1)--(c23) (a1)--(c24) (a2)--(c34);
  \foreach \v in {a0,a1,a2,a3,a4,c12,c13,c14,c23,c24,c34}
    \fill (\v) circle[radius=1.08pt];
  \foreach \i in {0,...,4} \node[right] at (a\i) {$a_{\i}$};
  \foreach \i/\j in {1/2,1/3,1/4,2/3,2/4,3/4}
    \node[below] at (c\i\j) {$c_{\i\j}$};
  \node[graphgreen,anchor=west] at (2.05,1.64) {$\gamma(c_{14})$};
  \node at (1.15,-.50) {(a) $T$, $m=4$};
\end{scope}
\begin{scope}[shift={(3.65,0)}]
  \path[region] (.825,0) ellipse[x radius=1.00cm,y radius=.12cm];
  \coordinate (d1) at (0,0);
  \coordinate (d2) at (.55,0);
  \coordinate (d3) at (1.10,0);
  \coordinate (d4) at (1.65,0);
  \coordinate (d) at (.825,1.95);
  \draw[arc,graphgreen] (d)--(d1);
  \draw[arc] (d)--(d2);
  \draw[arc] (d)--(d3);
  \draw[arc,graphgreen] (d)--(d4);
  \foreach \v in {d,d1,d2,d3,d4} \fill (\v) circle[radius=1.08pt];
  \node[above] at (d) {$c_{14}$};
  \foreach \i in {1,...,4} \node[below] at (d\i) {$a_{\i}$};
  \node[graphgreen,left] at (.34,.98) {$E_t$};
  \node[graphgreen,right] at (1.31,.98) {$E_b$};
  \node at (.825,-.50) {(b) Arcs of $D$ at $c_{14}$};
\end{scope}
\node[graphgreen] at (6.30,1.0) {$\longrightarrow$};
\begin{scope}[shift={(7.15,0)}]
  \path[region] (1.10,0) ellipse[x radius=1.28cm,y radius=.12cm];
  \coordinate (e1) at (0,0);
  \coordinate (e2) at (.7333,0);
  \coordinate (e3) at (1.4667,0);
  \coordinate (e4) at (2.2,0);
  \coordinate (w) at (0,.85);
  \coordinate (z) at (.42,1.64);
  \coordinate (e) at (1.10,1.95);
  \coordinate (p) at (1.74,1.30);
  \coordinate (q) at (2.20,.65);
  \draw[arc,graphgreen] (w)--(e1);
  \draw[arc,graphgreen] (w)--(z);
  \draw[arc,graphgreen] (z)--(e);
  \draw[arc,graphgreen] (e)--(p);
  \draw[arc,graphgreen] (p)--(q);
  \draw[arc,graphgreen] (q)--(e4);
  \draw[arc] (e)--(e2);
  \draw[arc] (e)--(e3);
  \foreach \v in {e1,e2,e3,e4,w,z,e,p,q} \fill (\v) circle[radius=1.08pt];
  \foreach \i in {1,...,4} \node[below] at (e\i) {$a_{\i}$};
  \node[above] at (e) {$c_{14}$};
  \node[left] at (w) {$w_{14}$};
  \node at (1.10,-.50) {(c) Replacements};
\end{scope}
\end{tikzpicture}
\caption{The construction in \cref{lem:burling-construction}. In~(a), the shaded path is $\gamma(c_{14})=\{a_1,a_2,a_3,a_4\}$; each $a_i$ is the last-born child of $a_{i-1}$.~(b) and~(c) show the arcs incident with $c_{14}$ and their replacements.}
\label{fig:burling-construction-proof}
\end{figure}
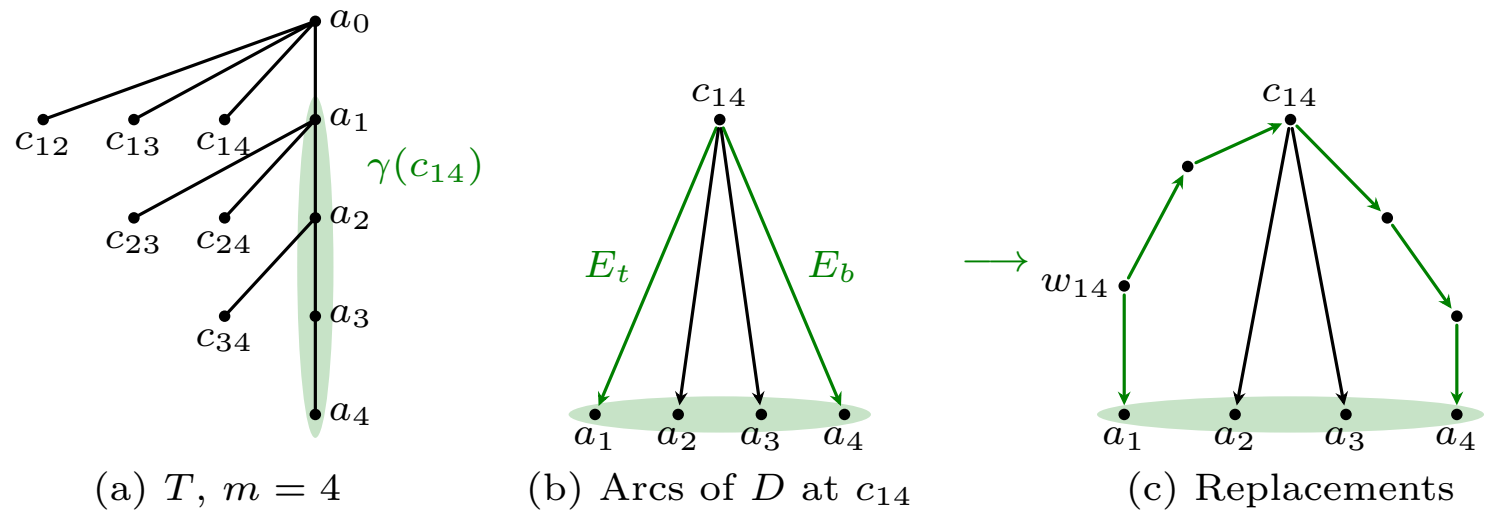

\begin{proof}
Take the rooted tree with spine $a_0a_1\cdots a_m$ and additional
children $c_{ij}$ of $a_{i-1}$ for $1\le i<j\le m$.
Make $a_i$ last-born at $a_{i-1}$, and let
\[
 \gamma(c_{ij})=\{a_i,\ldots,a_j\},\qquad
 \gamma(a_i)=\varnothing.
\]
Retain every vertex except $a_0$. The derived digraph has precisely
the arcs $c_{ij}\to a_h$ with $i\le h\le j$.
Each $c_{ij}$ is a source; its arcs to $a_i,a_j$ are respectively
top and bottom. By \cref{lem:pt-permitted-subdivision}, replace
$c_{ij}\to a_j$ by a directed path of length three, and
$c_{ij}\to a_i$ by $w_{ij}\to a_i$ together with a directed
$w_{ij}$--$c_{ij}$ path of length two.
The underlying graph is exactly $J_m$, and is Burling by
\cref{lem:pt-derived}; see \cref{fig:burling-construction-proof}.
\end{proof}

\begin{lemma}\label{lem:burling-pendant}
Adding a pendant vertex to a Burling graph preserves the Burling
property. Consequently, attaching a private rooted tree at a vertex
of a Burling graph preserves that property.
\end{lemma}

\begin{proof}
By \cref{lem:pt-derived,lem:pt-normalization}, represent $G$ by a
Burling tree with choice function $\gamma$ and retained set $X$
omitting the root and all last-born vertices. For $v\in X$, let $p$
be its parent, and let $p'$ be the last-born child of $p$.
If $\gamma(v)=\varnothing$, replace it by $\{p'\}$; this
does not change the derived digraph on $X$.
Attach a new leaf $w$ to the last vertex $y$ of $\gamma(v)$, and extend
only $\gamma(v)$ by $w$. If $y$ was a leaf of the tree, make $w$ its
last-born child; otherwise retain the old designation. Give $w$ empty
choice and retain $X\cup\{w\}$. All choices remain valid, and the
only new arc is $v\to w$. Thus the underlying graph is the required
pendant extension. Repetition attaches any private rooted tree.
\end{proof}

\begin{lemma}\label{lem:tree-burling-obstruction}
For every tree $H$ with at least two vertices and every $k\in\Nzero$,
there is a Burling graph $F_{H,k}$ with
$\branchpi_H(F_{H,k})>k$.
\end{lemma}

\begin{proof}
Choose a leaf $h$ of $H$ with neighbor $u$, let $J=J_{6k+3}$ with
core $C$, and form $F_{H,k}$ by identifying each $v\in V(J)$ with the
root of a private copy $H_v$ of $(H-h,u)$. By
\cref{lem:burling-construction,lem:burling-pendant}, $F_{H,k}$ is
Burling, and $J$ is an induced subgraph of $F_{H,k}$.

Fix a branch decomposition of $F_{H,k}$. By \cref{lem:balanced-core},
applied to $W=C$, a displayed cut $(A,B)$ has at least $2k+1$ core
vertices on each side. Pair distinct $x_i\in A\cap C$ and
$y_i\in B\cap C$ for $i\in[2k+1]$, and choose an edge $f_i$ of
$Q_{x_iy_i}$ crossing the cut. Since $Q_{x_iy_i}$ has length six,
each $f_i$ contains at most one core vertex. Among the $2k+1$ edges,
at least $k+1$ either all avoid the core or all meet it. In the
former case they lie in distinct, pairwise anticomplete sets
$U_{x_iy_i}$. In the latter case each consists of a core vertex and
its neighbor on $Q_{x_iy_i}$, which is a noncentral vertex of degree
two in $J$; since the core is independent, these edges are again
pairwise anticomplete. We obtain $k+1$ pairwise anticomplete edges
$v_iw_i$ of $J$ crossing the cut. The attachments are private, so the
induced copies $F_{H,k}[V(H_{v_i})\cup\{w_i\}]$ of $H$ are pairwise
anticomplete, and each crosses the cut. Since the decomposition was
arbitrary, $\branchpi_H(F_{H,k})>k$.
\end{proof}

\subsection{From control to coloring}
\label{subsec:control-to-coloring}

To derive \cref{thm:mainchiboundedness} from
\cref{thm:common-two-control}, it remains to turn $2$-control into
$\chi$-boundedness.

\begin{theorem}\label{thm:two-control-colouring}
Let $H$ be a subdivision of a star with at least two vertices, and
let $k\ge0$. Every $2$-controlled class contained in
$\{G\mid\branchpi_H(G)\le k\}$ is $\chi$-bounded.
\end{theorem}

We first derive \cref{thm:mainchiboundedness} from
\cref{thm:common-two-control,thm:two-control-colouring}, and then prove
\cref{thm:two-control-colouring}. Both proofs use the following lemma,
which transfers lower bounds on $\branchpi_F$ to $\branchpi_H$ for
suitable induced subgraphs $H$ of $F$.

\begin{lemma}\label{lem:edge-extension}
Let $F$ and $H$ be graphs with at least two vertices, where $F$ is
connected. If every
edge of $F$ belongs to an induced copy of $H$ in $F$, then
$\branchpi_F(G)\le\branchpi_H(G)$ for every graph $G$.
\end{lemma}
\begin{proof}
In every member $F_i$ of a crossing induced $F$-packing, choose a
crossing edge and an induced copy $H_i\subseteq F_i$ of $H$ containing it.
The $H_i$ remain pairwise anticomplete and cross the same cut. Thus
$\pi_F^\times(G;A,B)\le\pi_H^\times(G;A,B)$ for every cut.
Take maxima over displayed cuts and then minima over decompositions.
\end{proof}

We are now ready to prove \cref{thm:mainchiboundedness}, assuming
\cref{thm:common-two-control,thm:two-control-colouring}.

\subdividedClawsAndBrooms*

\begin{proof}
Fix $H$ and $k\ge0$. We choose a graph $F$ from one of the families in
\cref{thm:common-two-control} such that every edge of $F$ belongs to
an induced copy of $H$ in $F$.

If $H=S_t$, let $F=H=K_{1,3}^{(t)}$. It belongs to
family~(\ref{item:common-stars}) if $t=1$, to
family~(\ref{item:common-short-claw}) if $t=2$, and to
family~(\ref{item:common-long-arms}) if $t\ge3$.

Now let $H=B_{t,d}$. If $H$ is a path, that is, if $d=1$ or
$(t,d)=(1,2)$, let $F=H$, from family~(\ref{item:common-paths}).
Otherwise, let $F=K_{1,t+d-1}$ when $t\le2$, let $F=B_{3,d}$ when
$t=3$, and let $F=K_{1,d+1}^{(t-1)}$ when $t\ge4$; these graphs belong
to families~(\ref{item:common-stars}), (\ref{item:common-short-broom}),
and~(\ref{item:common-long-arms}), respectively. For $t\le3$, $F=H$.
For $t\ge4$, every edge of $F$ belongs to an induced copy of $H$: retain the
entire arm containing it and the first edge of each other arm.

In every case, \cref{lem:edge-extension,lem:node-span-width} give
\[
 \sigma_F(G)\le\left\lfloor\frac{3\branchpi_F(G)}2\right\rfloor
 \le\left\lfloor\frac{3\branchpi_H(G)}2\right\rfloor,
 \quad\text{so}\quad
 \{G\mid\branchpi_H(G)\le k\}\subseteq\cC_{\lfloor3k/2\rfloor}(F).
\]
The class on the right is $2$-controlled by
\cref{thm:common-two-control}. The class on the left is hereditary by
\cref{lem:basic}\ref{item:heredity}, and is therefore $2$-controlled as
well. Since $H$ is a subdivision of a star with at least two vertices,
\cref{thm:two-control-colouring} shows that it is $\chi$-bounded.
\end{proof}

It remains to prove \cref{thm:two-control-colouring}. It combines a
theorem of Abrishami et al.\ (\cref{lem:2-controlled-Kmm}) with two
obstructions: the Burling graphs of \cref{lem:tree-burling-obstruction}
and the subdivided bicliques of the next lemma. The biclique
obstructions also supply the excluded subdivisions used in
\cref{subsec:paths-stars,subsec:long-arms,subsec:arms-two}.

\begin{lemma}\label{lem:short-biclique-obstruction}
Let $k,r,s\ge1$ and $2\le q\le2r$. If $m\ge(sr+2)k$, then
$\branchpi_{K_{1,s}^{(r)}}(K_{m,m}^{(q)})\ge k$.
\end{lemma}
\begin{proof}
Let $X,Y$ be the original parts, and fix an arbitrary branch
decomposition. \Cref{lem:balanced-core}, applied to $W=X$ only, gives
a displayed cut $(A,B)$ with at least $m/3$ vertices of $X$ on each
side. By the pigeonhole principle, one side contains at least $m/2$
vertices of $Y$; after possibly interchanging sides,
\[
 |A\cap X|,|B\cap X|\ge m/3,\qquad |A\cap Y|\ge m/2.
\]
Since $sr+2\ge3$, both $B\cap X$ and $A\cap Y$ have at least $k$
vertices. Choose distinct
$x_i\in B\cap X$, $y_i\in A\cap Y$, $i\in[k]$, and let $Q_i$
be their replacement path. Initially reserve $R_i=\{x_i,y_i\}$.
Whenever construction $i$ uses any part of a replacement path,
reserve both original endpoints in $R_i$, choosing every new endpoint
outside $\bigcup_jR_j$. Then the sets $R_i$ are pairwise disjoint,
and copy $i$ consists of vertices of $R_i$ and interior vertices of
replacement paths with both ends in $R_i$. Original vertices are
pairwise nonadjacent, since $q\ge2$, and an interior vertex of a
replacement path is adjacent only to vertices of that path. Hence the
constructed copies are pairwise anticomplete.

If $q\le r$, grow $s$ arms of length $r$ from $x_i$, starting
one with $Q_i$ and otherwise using fresh original vertices; truncate
the final replacement paths. If $r<q\le2r$, the initial $r$-edge
segments from the two ends of $Q_i$ cover every edge of $Q_i$.
Choose one containing a crossing edge, and complete $s-1$ further
$r$-edge arms from the same endpoint towards fresh original vertices.
Both constructions induce a crossing $K_{1,s}^{(r)}$ and reserve at most
$sr+2$ original vertices. Hence fewer than $(sr+2)k\le m$ original
vertices are reserved whenever a new endpoint is chosen, so every
required vertex can be chosen unreserved in either original part.
The reservation rule gives a crossing induced packing of size $k$ at
$(A,B)$. Since the decomposition was arbitrary, the lemma follows.
\end{proof}

\begin{proof}[Proof of \cref{thm:two-control-colouring}]
Let $\cC$ be the given class. Write $H$ as a subdivision of $K_{1,a}$
with $a\ge1$.
Let $r$ be the maximum length of an arm of $H$, and let
$m=(ar+2)(k+1)$.
Every edge of $K_{1,a}^{(r)}$ belongs to an induced copy of $H$, obtained by
retaining the whole arm containing the edge and truncating the other
arms. Hence
\cref{lem:edge-extension,lem:short-biclique-obstruction} give
$\branchpi_H(K_{m,m}^{(2)})\ge k+1$.
Since $\branchpi_H$ does not increase under taking induced subgraphs
(\cref{lem:basic}\ref{item:heredity}), $\cC$ is $K_{m,m}^{(2)}$-free.
Choose a Burling graph $F_{H,k}$ from
\cref{lem:tree-burling-obstruction}, and $j$ such that
$F_{H,k}$ is induced in $\mathsf{Bur}_j$. For the same reason, every
graph in $\cC$ is $\mathsf{Bur}_j$-free. Hence $\cC$ is $\chi$-bounded
by \cref{lem:2-controlled-Kmm}.
\end{proof}

\subsection{Results from the literature}
\label{subsec:literature}

The rest of this section is devoted to the proof of
\cref{thm:common-two-control}. We first collect the results of other
authors that it uses, restated in our notation.

\paragraph{Chromatic bounds.}
\begin{lemma}[Scott~{\cite[Corollary~2]{scott1997induced}}]
\label{lem:Scott}
If $H$ is a subdivision of a star, then the class of $H$-free graphs
is $\chi$-bounded.
\end{lemma}

\begin{lemma}[Chudnovsky, Scott, and Seymour~{\cite[Theorem~1.4]{CSSLH}}]
\label{lem:long-holes}
For every $\ell\ge4$, the class of graphs with no induced cycle of
length at least $\ell$ is $\chi$-bounded.
\end{lemma}

\paragraph{Coarse partitions.}
For a partition $\mathcal P$ of $V(F)$ into nonempty sets, the
quotient $F/\mathcal P$ has vertex set $\mathcal P$, with two parts
adjacent exactly when an edge of $F$ joins them. The following lemma
is obtained by combining Theorem~2.3 of Nguyen, Scott, and
Seymour~\cite{NSSCoarse} with Theorem~2 of
Distel~\cite{distel2025alternative}.

\begin{lemma}[Nguyen, Scott, and Seymour~{\cite[Theorem~2.3]{NSSCoarse}}; Distel~{\cite[Theorem~2]{distel2025alternative}}]\label{lem:two-Distel}
Let $k\ge1$, and let
$W_{\ref{lem:two-Distel}}(k)=16(k+1)(k+2)^2+2(k+2)$. If $F$ has a tree
decomposition whose bags are unions of at most $k$ cliques, then $F$
has a partition $\mathcal P$ with $\tw(F/\mathcal P)\le k$ and
$\diam_F(P)\le W_{\ref{lem:two-Distel}}(k)$ for every
$P\in\mathcal P$.
\end{lemma}

\paragraph{Multicovers.}
A \emph{$\rho$-cover} of a set $C$ is a sequence
$\mathcal L=(L_0,\ldots,L_{\rho-1})$ of pairwise disjoint sets, its
\emph{levels}, that are disjoint from $C$, such that $L_0$ is a
singleton, $L_a$ covers $L_{a+1}$ for $0\le a<\rho-1$, the last level
covers $C$, nonconsecutive levels are anticomplete, and all earlier
levels are anticomplete to $C$. Write
$V(\mathcal L)=\bigcup_{a=0}^{\rho-1}L_a$. The vertex of $L_0$
is its \emph{apex}, $L_{\rho-1}$ is its \emph{base}, and the
union of the earlier levels is its \emph{core}.

\begin{figure}[t]
\centering
\definecolor{graphpurple}{RGB}{128,0,128}
\begin{tikzpicture}[
  font=\fontsize{7.2}{8.6}\selectfont,
  line width=.60pt,line cap=round,line join=round,
  every node/.style={inner sep=2.2pt},
  edge/.style={draw=black},
  arc/.style={-{Stealth[length=2.6pt,width=2.4pt]},shorten <=1.7pt,shorten >=1.7pt},
  region/.style={fill=graphpurple,fill opacity=.22,draw=none}]
\begin{scope}
  \path[region] (1.2,2.7) ellipse[x radius=0.18cm,y radius=.12cm];
  \node[graphpurple,anchor=east] at (-.15,2.7) {$L_0$};
  \path[region] (1.2,1.9) ellipse[x radius=0.6cm,y radius=.12cm];
  \node[graphpurple,anchor=east] at (-.15,1.9) {$L_1$};
  \path[region] (1.2,1.1) ellipse[x radius=1.0cm,y radius=.12cm];
  \node[graphpurple,anchor=east] at (-.15,1.1) {$L_2$};
  \path[region] (1.2,0.3) ellipse[x radius=1.2cm,y radius=.12cm];
  \node[graphpurple,anchor=east] at (-.15,0.3) {$L_3$};
  \coordinate (lx) at (1.2,2.7);
  \coordinate (lu1) at (0.8,1.9);
  \coordinate (lu2) at (1.6,1.9);
  \coordinate (lv1) at (0.4,1.1);
  \coordinate (lv2) at (1.2,1.1);
  \coordinate (lv3) at (2.0,1.1);
  \coordinate (lb1) at (0.2,0.3);
  \coordinate (lb2) at (0.8667,0.3);
  \coordinate (lb3) at (1.5333,0.3);
  \coordinate (lb4) at (2.2,0.3);
  \draw[edge] (lx)--(lu1);
  \draw[edge] (lx)--(lu2);
  \draw[edge] (lu1)--(lv1);
  \draw[edge] (lu1)--(lv2);
  \draw[edge] (lu2)--(lv2);
  \draw[edge] (lu2)--(lv3);
  \draw[edge] (lv1)--(lb1);
  \draw[edge] (lv1)--(lb2);
  \draw[edge] (lv2)--(lb2);
  \draw[edge] (lv2)--(lb3);
  \draw[edge] (lv3)--(lb3);
  \draw[edge] (lv3)--(lb4);
  \draw[edge] (lb2)--(lb3);
  \foreach \v in {lx,lu1,lu2,lv1,lv2,lv3,lb1,lb2,lb3,lb4} \fill (\v) circle[radius=1.08pt];
  \node[right] at (lx) {$x$};
  \node at (1.05,-.50) {(a) Levels of a $4$-cover};
\end{scope}
\begin{scope}[shift={(4.2,0)}]
  \path[region] (.9,2.05) ellipse[x radius=.56cm,y radius=.88cm];
  \path[region] (3.3,2.05) ellipse[x radius=.56cm,y radius=.88cm];
  \path[region] (.9,.75) ellipse[x radius=.50cm,y radius=.12cm];
  \path[region] (3.3,.75) ellipse[x radius=.50cm,y radius=.12cm];
  \path[region] (2.1,0) ellipse[x radius=1.60cm,y radius=.13cm];
  \node[graphpurple] at (.9,3.08) {$\mathcal L_1$};
  \node[graphpurple] at (3.3,3.08) {$\mathcal L_2$};
  \node[graphpurple,left] at (.31,2.05) {$Z_1$};
  \node[graphpurple,right] at (3.89,2.05) {$Z_2$};
  \node[graphpurple,left] at (.38,.75) {$L_{1,3}$};
  \node[graphpurple,right] at (3.82,.75) {$L_{2,3}$};
  \node[graphpurple,left] at (.45,0) {$C$};
  \coordinate (x1) at (0.9,2.7);
  \coordinate (u1) at (0.9,2.05);
  \coordinate (v11) at (0.6000000000000001,1.4);
  \coordinate (v12) at (1.2,1.4);
  \coordinate (b11) at (0.6000000000000001,0.75);
  \coordinate (b12) at (1.2,0.75);
  \coordinate (x2) at (3.3,2.7);
  \coordinate (u2) at (3.3,2.05);
  \coordinate (v21) at (3.0,1.4);
  \coordinate (v22) at (3.5999999999999996,1.4);
  \coordinate (b21) at (3.0,0.75);
  \coordinate (b22) at (3.5999999999999996,0.75);
  \coordinate (t1) at (0.8,0.0);
  \coordinate (t2) at (2.1,0.0);
  \coordinate (t3) at (3.4,0.0);
  \draw[edge] (x1)--(u1);
  \draw[edge] (u1)--(v11);
  \draw[edge] (u1)--(v12);
  \draw[edge] (v11)--(b11);
  \draw[edge] (v12)--(b12);
  \draw[edge] (x2)--(u2);
  \draw[edge] (u2)--(v21);
  \draw[edge] (u2)--(v22);
  \draw[edge] (v21)--(b21);
  \draw[edge] (v22)--(b22);
  \draw[edge] (b11)--(b12);
  \draw[edge] (b11)--(t1);
  \draw[edge] (b11)--(t2);
  \draw[edge] (b12)--(t3);
  \draw[edge] (b21)--(t1);
  \draw[edge] (b22)--(t2);
  \draw[edge] (b22)--(t3);
  \draw[edge] (t1)--(t2);
  \draw[edge] (t2)--(t3);
  \draw[graphpurple] (b12)--(v21);
  \foreach \v in {x1,u1,v11,v12,b11,b12,x2,u2,v21,v22,b21,b22,t1,t2,t3} \fill (\v) circle[radius=1.08pt];
  \node[right] at (x1) {$x_1$};
  \node[right] at (x2) {$x_2$};
  \node[graphpurple,above,inner sep=1.5pt] at (2.05,1.28) {$L_{1,3}$--$Z_2$};
  \node at (2.1,-.50) {(b) An ordered $4$-multicover};
\end{scope}
\end{tikzpicture}
\caption{The levels of a $4$-cover~(a), and an ordered $4$-multicover of $C$ of length $2$~(b). The shaded regions mark the levels in~(a) and the cores $Z_i$, the bases $L_{i,3}$, and the target $C$ in~(b). The edge between $L_{1,3}$ and $Z_2$ is allowed, since it joins an earlier base to a later core.}
\label{fig:multicover}
\end{figure}
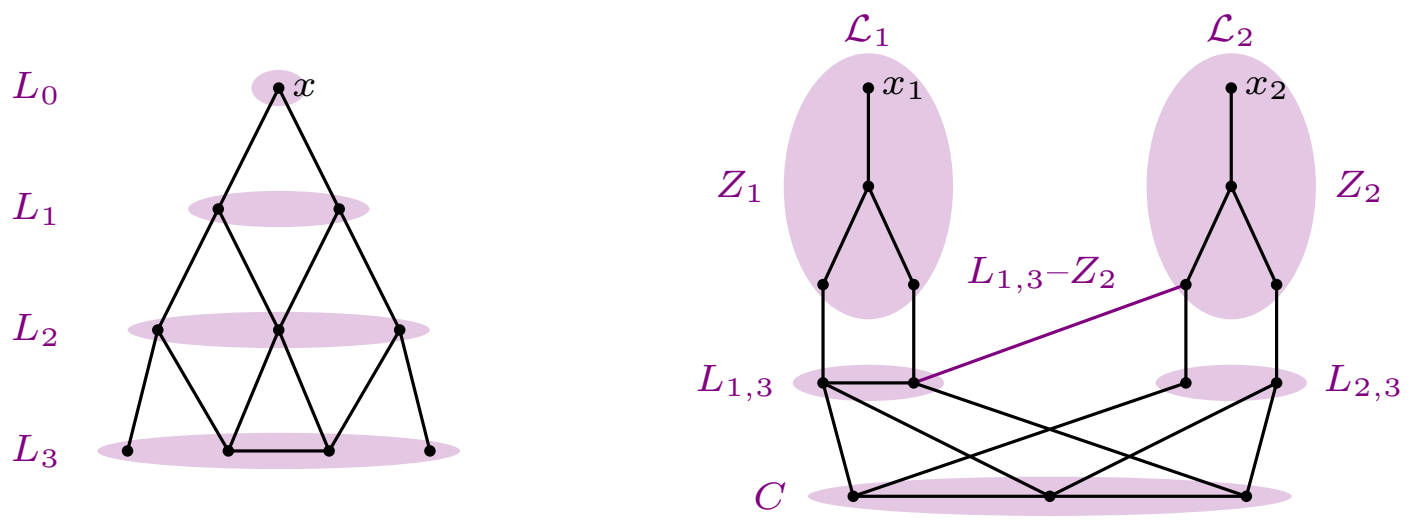

Following~\cite[Section~4]{CSSChandeliers}, an
\emph{ordered $\rho$-multicover} of $C$ of \emph{length} $N$ is a
sequence $\mathcal M=(\mathcal L_i)_{i=1}^N$ of $\rho$-covers
$\mathcal L_i=(L_{i,0},\ldots,L_{i,\rho-1})$ of $C$ with pairwise
disjoint vertex sets such that, for $i<j$, the core of
$\mathcal L_i$ is anticomplete to all of $\mathcal L_j$.
We call $C$ its \emph{target} and write
$V(\mathcal M)=\bigcup_{i=1}^NV(\mathcal L_i)$.
Thus every edge from an earlier cover to a later cover has its
earlier endpoint in the earlier base; in particular, edges from an
earlier base to a later core are allowed. No independence condition
is imposed on the bases.
An ordered $\rho$-multicover $(\mathcal L'_i)_{i=1}^N$ of $C'$ is
\emph{contained} in $(\mathcal L_i)_{i=1}^N$ if
$L'_{i,a}\subseteq L_{i,a}$ for all $i\in[N]$ and $0\le a<\rho$.
Then $\mathcal L'_i$ and $\mathcal L_i$ have the same apex.
See~\cref{fig:multicover} for examples.

\begin{figure}[t]
\centering
\providecolor{graphgreen}{RGB}{0,128,0}
\begin{tikzpicture}[
  font=\fontsize{7.2}{8.6}\selectfont,
  line width=.60pt,line cap=round,line join=round,
  every node/.style={inner sep=2.2pt},
  region/.style={fill=graphgreen,fill opacity=.22,draw=none}]
  \path[region] (6,2.9) ellipse[x radius=.14cm,y radius=.14cm];
  \path[region] (7.426585,1.863525) ellipse[x radius=.14cm,y radius=.14cm];
  \path[region] (6.881678,0.186475) ellipse[x radius=.14cm,y radius=.14cm];
  \path[region] (5.118322,0.186475) ellipse[x radius=.14cm,y radius=.14cm];
  \path[region] (4.573415,1.863525) ellipse[x radius=.14cm,y radius=.14cm];
  \coordinate (a1) at (1.1,2.35);
  \coordinate (a2) at (2.003504,1.693566);
  \coordinate (a3) at (1.658396,0.631434);
  \coordinate (a4) at (0.541604,0.631434);
  \coordinate (a5) at (0.196496,1.693566);
  \coordinate (x1) at (6,2.9);
  \coordinate (x2) at (7.426585,1.863525);
  \coordinate (x3) at (6.881678,0.186475);
  \coordinate (x4) at (5.118322,0.186475);
  \coordinate (x5) at (4.573415,1.863525);
  \coordinate (p121) at (6.356646,2.640881);
  \coordinate (p122) at (7.069939,2.122644);
  \coordinate (p131) at (6.220419,2.221619);
  \coordinate (p132) at (6.661258,0.864856);
  \coordinate (p141) at (5.779581,2.221619);
  \coordinate (p142) at (5.338742,0.864856);
  \coordinate (p151) at (5.643354,2.640881);
  \coordinate (p152) at (4.930061,2.122644);
  \coordinate (p231) at (7.290358,1.444263);
  \coordinate (p232) at (7.017905,0.605737);
  \coordinate (p241) at (6.849519,1.444263);
  \coordinate (p242) at (5.695388,0.605737);
  \coordinate (p251) at (6.713292,1.863525);
  \coordinate (p252) at (5.286708,1.863525);
  \coordinate (p341) at (6.440839,0.186475);
  \coordinate (p342) at (5.559161,0.186475);
  \coordinate (p351) at (6.304612,0.605737);
  \coordinate (p352) at (5.150481,1.444263);
  \coordinate (p451) at (4.982095,0.605737);
  \coordinate (p452) at (4.709642,1.444263);
  \draw (a1) -- (a2);
  \draw (a1) -- (a3);
  \draw (a1) -- (a4);
  \draw (a1) -- (a5);
  \draw (a2) -- (a3);
  \draw (a2) -- (a4);
  \draw (a2) -- (a5);
  \draw (a3) -- (a4);
  \draw (a3) -- (a5);
  \draw (a4) -- (a5);
  \draw (p121) -- (p122);
  \draw[graphgreen] (p121) -- (p131);
  \draw (p121) -- (x1);
  \draw (p122) -- (x2);
  \draw (p131) -- (p132);
  \draw (p131) -- (x1);
  \draw (p132) -- (x3);
  \draw (p141) -- (p142);
  \draw (p141) -- (x1);
  \draw (p142) -- (x4);
  \draw (p151) -- (p152);
  \draw (p151) -- (x1);
  \draw (p152) -- (x5);
  \draw (p231) -- (p232);
  \draw (p231) -- (x2);
  \draw (p232) -- (x3);
  \draw (p241) -- (p242);
  \draw (p241) -- (x2);
  \draw (p242) -- (x4);
  \draw (p251) -- (p252);
  \draw (p251) -- (x2);
  \draw (p252) -- (x5);
  \draw (p341) -- (p342);
  \draw (p341) -- (x3);
  \draw (p342) -- (x4);
  \draw (p351) -- (p352);
  \draw (p351) -- (x3);
  \draw (p352) -- (x5);
  \draw (p451) -- (p452);
  \draw (p451) -- (x4);
  \draw (p452) -- (x5);
  \fill[black] (a1) circle[radius=1.08pt];
  \fill[black] (a2) circle[radius=1.08pt];
  \fill[black] (a3) circle[radius=1.08pt];
  \fill[black] (a4) circle[radius=1.08pt];
  \fill[black] (a5) circle[radius=1.08pt];
  \fill[black] (x1) circle[radius=1.08pt];
  \fill[black] (x2) circle[radius=1.08pt];
  \fill[black] (x3) circle[radius=1.08pt];
  \fill[black] (x4) circle[radius=1.08pt];
  \fill[black] (x5) circle[radius=1.08pt];
  \fill[black] (p121) circle[radius=1.08pt];
  \fill[black] (p122) circle[radius=1.08pt];
  \fill[black] (p131) circle[radius=1.08pt];
  \fill[black] (p132) circle[radius=1.08pt];
  \fill[black] (p141) circle[radius=1.08pt];
  \fill[black] (p142) circle[radius=1.08pt];
  \fill[black] (p151) circle[radius=1.08pt];
  \fill[black] (p152) circle[radius=1.08pt];
  \fill[black] (p231) circle[radius=1.08pt];
  \fill[black] (p232) circle[radius=1.08pt];
  \fill[black] (p241) circle[radius=1.08pt];
  \fill[black] (p242) circle[radius=1.08pt];
  \fill[black] (p251) circle[radius=1.08pt];
  \fill[black] (p252) circle[radius=1.08pt];
  \fill[black] (p341) circle[radius=1.08pt];
  \fill[black] (p342) circle[radius=1.08pt];
  \fill[black] (p351) circle[radius=1.08pt];
  \fill[black] (p352) circle[radius=1.08pt];
  \fill[black] (p451) circle[radius=1.08pt];
  \fill[black] (p452) circle[radius=1.08pt];
  \node[above] at (a1) {$a_1$};
  \node[right] at (a2) {$a_2$};
  \node[below right] at (a3) {$a_3$};
  \node[below left] at (a4) {$a_4$};
  \node[left] at (a5) {$a_5$};
  \node[above] at (x1) {$x_1$};
  \node[right] at (x2) {$x_2$};
  \node[below right] at (x3) {$x_3$};
  \node[below left] at (x4) {$x_4$};
  \node[left] at (x5) {$x_5$};
  \draw[graphgreen,-{Stealth[length=2.6pt,width=2.4pt]}]
    (2.60,1.40) -- (3.82,1.40) node[midway,above] {$\eta$};
  \node at (1.10,-.52) {(a) $F=K_5$};
  \node at (6.00,-.52) {(b) An impression of $F$};
\end{tikzpicture}
\caption{An impression of $F=K_5$ of order three: $\eta(a_i)=x_i$ (shaded), and $\eta(a_ia_j)=P_{ij}$, the three-edge path along the segment $x_ix_j$. The green edge joins interiors of $P_{12}$ and $P_{13}$, which represent incident edges of $F$.}
\label{fig:impression}
\end{figure}
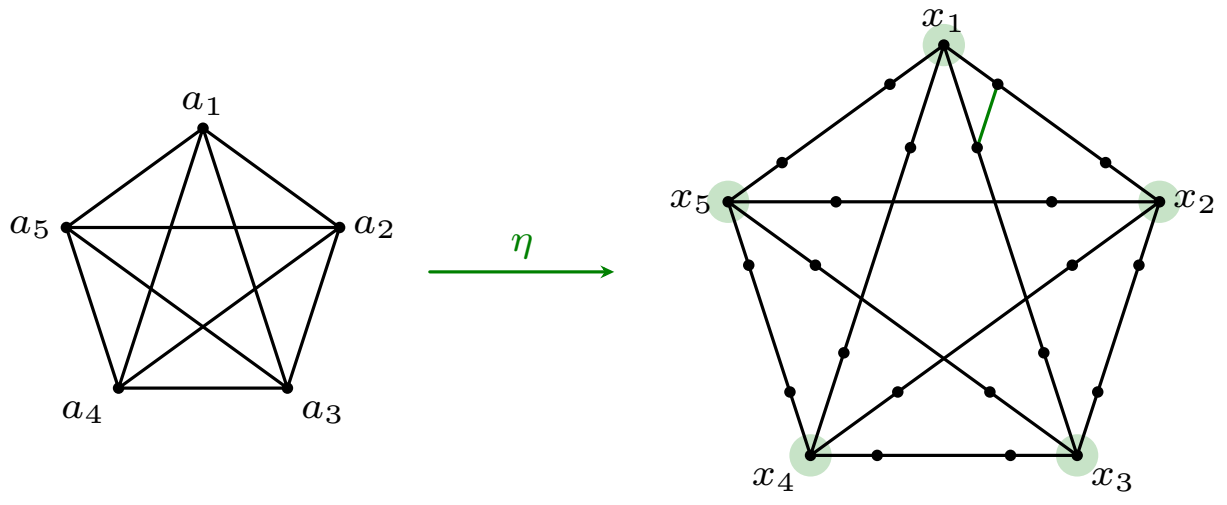

\paragraph{Impressions.}
As in~\cite[Section~3]{CSSChandeliers}, an \emph{impression} $\eta$ of
a graph $F$ in $G$ assigns to the vertices $v$ of $F$ pairwise
distinct and nonadjacent vertices $\eta(v)$ of $G$, and to each edge
$uv$ of $F$ a path $\eta(uv)$ of $G$ between $\eta(u)$ and $\eta(v)$.
Paths for edges with no common endpoint must be disjoint and
anticomplete. Its \emph{order} is the maximum length of its paths;
see \cref{fig:impression}.

\begin{lemma}[Chudnovsky, Scott, and Seymour~{\cite[Theorem~3.2]{CSSChandeliers}}]
\label{lem:css-impressions}
For all $\rho\ge2$ and $m,\nu\ge1$, there is $\mu\ge1$ such that no graph
$G$ with $\omega(G)\le\nu$ and with no induced copy of $K_{m,m}^{(q)}$ for
$2\le q\le\rho$ has an impression of $K_{\mu,\mu}$ of order at most $\rho$.
\end{lemma}

\begin{lemma}[Chudnovsky, Scott, and Seymour~{\cite[Theorem~4.4]{CSSChandeliers}}]
\label{lem:css-two-control}
Let $\rho\ge2$ and $m\ge1$. Every $\rho$-controlled class
that excludes $K_{m,m}^{(q)}$ for every $2\le q\le\rho+3$ is
$2$-controlled.
\end{lemma}

\begin{lemma}[Chudnovsky, Scott, and Seymour~{\cite[Theorem~5.3]{CSSChandeliers}}]
\label{lem:css-exact-paths}
Let $\rho\ge4$, $N,\tau\ge1$, and $c'\ge0$, and let
$c=c'+(N(\rho-1)+1)\tau$. Let $G$ satisfy $\chi^{\rho-1}(G)\le\tau$,
and let $\mathcal M=(\mathcal L_i)_{i=1}^N$ be an ordered
$\rho$-multicover of $C$ with $\chi(G[C])>c$. Then there are
$C'\subseteq C$ with $\chi(G[C'])>c'$, an ordered $\rho$-multicover
$\mathcal M'=(\mathcal L'_i)_{i=1}^N$ of $C'$ contained in
$\mathcal M$, a vertex $z\in C\setminus C'$, and, for every $i\in[N]$,
an induced path $S_i$ with exactly $\rho$ edges from $z$ to the apex of
$\mathcal L_i$, such that $V(S_i)\subseteq V(\mathcal L_i)\cup\{z\}$
and, for $W=C'\cup V(\mathcal M')$,
\[
 \{z\}\perp W,\qquad V(S_i)\perp W\setminus V(\mathcal L'_i).
\]
\end{lemma}

\subsection{Finite-radius control}
\label{subsec:finite-control}

We first show that the classes $\cC_k(H)$ are controlled at some
finite radius.

\begin{lemma}\label{lem:metric-subgraph-colouring}
Let $J$ be a spanning subgraph of $G$.
Suppose $J$ has a tree decomposition whose bags are unions of at most
$k$ sets of weak diameter at most $D$ in $G$, where $k,D\ge1$. Then
\[
 \chi(J)\le(k+1)\chi^{D\cdot W_{\ref{lem:two-Distel}}(k)}(G).
\]
\end{lemma}
\begin{proof}
Complete each of the indicated bag parts to a clique, obtaining $F$.
The same decomposition is valid for $F$, and
\[
 uv\in E(F)\Longrightarrow \dist_G(u,v)\le D,
 \qquad \dist_G(x,y)\le D\cdot\dist_F(x,y)\quad(x,y\in V(G)).
\]
By \cref{lem:two-Distel}, $F$ has a partition $\mathcal P$ with
$\tw(F/\mathcal P)\le k$ and
$\diam_G(P)\le D\cdot W_{\ref{lem:two-Distel}}(k)$.
Color $F/\mathcal P$ with $k+1$ colors and each $G[P]$ with
$\chi^{D\cdot W_{\ref{lem:two-Distel}}(k)}(G)$ colors. Assigning to each vertex the pair of
its part color and its color within the part properly colors $J$.
The empty graph is immediate.
\end{proof}

\begin{theorem}\label{thm:finite-control-packing}
Let $H$ be connected with at least two vertices, let $r$ be its
\emph{radius}, i.e., the least integer $r$ such that some vertex of $H$ is
within distance $r$ of every vertex of $H$, and suppose that $f_H$
binds the $H$-free graphs. For $k\ge1$, let
$\rho=2(3r+1)W_{\ref{lem:two-Distel}}(k)$. Every graph $G$ with
$\sigma_H(G)\le k$ satisfies
\[
 \chi(G)\le(k+1)\chi^\rho(G)f_H(\omega(G)).
\]
In particular, $\cC_k(H)$ is $\rho$-controlled.
\end{theorem}
\begin{proof}
Graphs with at most one vertex are immediate, so assume $|V(G)|\ge2$,
and fix a node-width-$k$ decomposition $(T,\lambda)$. Let $E_H$ consist of the edges belonging to induced
copies of $H$, and let $J=(V(G),E_H)$. For each $e\in E_H$, choose such
a copy $Q_e$, and let $z_e$ be the image in $Q_e$ of a vertex of $H$
within distance $r$ of every vertex of $H$. Then
$\dist_G(z_e,v)\le r$ for $v\in V(Q_e)$. Let
\[
 S_u=\{\lambda(u)\}\cup\bigcup_{u\in e\in E_H}T[Q_e],
 \qquad \beta(x)=\{u\mid x\in V(S_u)\}.
\]
Each $S_u$ is connected, and $e=uw\in E_H$ implies
$T[Q_e]\subseteq S_u\cap S_w$. Hence $(T,\beta)$ is a tree
decomposition of $J$.

For an internal node $x$ with $\beta(x)\ne\varnothing$, choose a
maximal induced packing $Q_{e_1},\ldots,Q_{e_p}$ among the selected
copies whose spans contain $x$. Then $1\le p\le k$. Every
$u\in\beta(x)$ belongs to one such $Q_e$, and maximality gives
$\dist_G(Q_e,Q_{e_i})\le1$ for some $i\in[p]$. All vertices of $Q_e$
lie within distance $2r$ of $u$, through $z_e$, and all vertices of
$Q_{e_i}$ lie within distance $r$ of $z_{e_i}$. Hence
\[
 \dist_G(u,z_{e_i})\le2r+1+r=3r+1.
\]
At a leaf $\lambda(v)$, $\beta(\lambda(v))\subseteq B_G(v,2r)$.
Thus every bag is a union of at most $k$ sets of weak diameter at
most $2(3r+1)$. By \cref{lem:metric-subgraph-colouring},
$\chi(J)\le(k+1)\chi^\rho(G)$.
Each color class of $J$ induces an $H$-free subgraph of $G$, since
all edges of an induced copy of $H$ belong to $E_H$. Coloring these classes
with disjoint palettes proves the bound. Finally,
$\omega(G)\le\chi^\rho(G)$ and \cref{lem:node-span-width} give
heredity and the control function $(k+1)\tau f_H(\tau)$.
\end{proof}

Every part of \cref{thm:common-two-control} except the one for stars
starts from the following consequence.

\begin{lemma}\label{lem:control-start}
Let $H$ be a subdivision of a star with at least two vertices. Then
$\cC_0(H)$ is $R$-controlled for every $R\ge1$. For every $k\ge1$,
the class $\cC_k(H)$ is controlled at some finite radius.
\end{lemma}
\begin{proof}
The class $\cC_0(H)$ consists of the $H$-free graphs. By
\cref{lem:Scott}, it is $\chi$-bounded, and hence $R$-controlled for
every $R\ge1$, since $\omega(G)\le\chi^R(G)$ for nonempty $G$.
For $k\ge1$, \cref{lem:Scott} gives a function binding the $H$-free
graphs, so \cref{thm:finite-control-packing} applies.
\end{proof}

\subsection{Radius-reduction tools}
\label{subsec:reduction-tools}

\subsubsection{Multicovers}
\label{subsec:multicovers}

The first reduction is the nested-layer argument from the proof
of~\cite[Theorem~4.2]{CSSChandeliers}.

\begin{lemma}\label{lem:multicover-control-reduction}
Let $\rho\ge2$, and let $\cC$ be a $\rho$-controlled
class. Suppose that, for every integer $\tau\ge1$, there are integers
$N=N(\tau)\ge1$ and $c=c(\tau)\ge0$ such that no $G\in\cC$ with
$\chi^{\rho-1}(G)\le\tau$ has an ordered $\rho$-multicover of
length $N$ whose target has chromatic number greater than $c$. Then
$\cC$ is $(\rho-1)$-controlled.
\end{lemma}
\begin{proof}
Fix $\tau\ge1$, let $\psi$ witness $\rho$-control of $\cC$, and let
$N=N(\tau)$ and $c=c(\tau)$. Define $c_0=c$ and
$c_{i+1}=\psi(c_i+\tau)$ for $0\le i<N$. Suppose that $G\in\cC$
satisfies $\chi^{\rho-1}(G)\le\tau$ and $\chi(G)>c_N$.

Let $C_0=V(G)$. For $j=1,\ldots,N$, suppose that
$\chi(G[C_{j-1}])>c_{N-j+1}$. Since $G[C_{j-1}]\in\cC$ by heredity and
$c_{N-j+1}=\psi(c_{N-j}+\tau)$, control gives $z_j\in C_{j-1}$ with
\[
 \chi\bigl(G[B_{G[C_{j-1}]}(z_j,\rho)]\bigr)>c_{N-j}+\tau.
\]
Let $L_{j,a}$ be the distance-$a$ layer of $z_j$ in $G[C_{j-1}]$, and
let $C_j=L_{j,\rho}$. The first $\rho$ layers lie in
$B_G(z_j,\rho-1)$, so they induce a $\tau$-colorable subgraph. Hence
$\chi(G[C_j])>c_{N-j}$.

Let $\mathcal L_j=(L_{j,0},\ldots,L_{j,\rho-1})$. Its base covers
$C_j$ and its earlier levels are anticomplete to $C_j$, so
$\mathcal L_j$ is a $\rho$-cover of $C_j$, and hence of
$C_N\subseteq C_j$. For $i<j$, the cover $\mathcal L_j$ lies in
$C_{j-1}\subseteq C_i$, which is disjoint from $\mathcal L_i$ and
anticomplete to its core. Thus $(\mathcal L_j)_{j=1}^N$ is an ordered
$\rho$-multicover of $C_N$ of length $N$ with $\chi(G[C_N])>c_0=c$, a
contradiction. Hence $\chi(G)\le c_N$, where $c_N$ depends only on
$\tau$, and the nondecreasing envelope of these bounds witnesses
$(\rho-1)$-control.
\end{proof}

For $\rho\ge4$, the exact path lengths in \cref{lem:css-exact-paths} give a
reduction that needs fewer excluded subdivisions than
\cref{lem:css-two-control}.

\begin{lemma}\label{lem:exact-path-reduction}
Let $\rho\ge4$ and $m\ge1$. Every $\rho$-controlled class
$\cC$ that excludes $K_{m,m}^{(q)}$ for every $2\le q\le\rho$ is
$(\rho-1)$-controlled.
\end{lemma}
\begin{proof}
Fix $\tau\ge1$, and let $G\in\cC$ satisfy $\chi^{\rho-1}(G)\le\tau$.
Then $\omega(G)\le\tau$. By \cref{lem:css-impressions} with
$\nu=\tau$, there is $\mu$, depending only on $\rho$, $m$, and
$\tau$, such that $G$ has no impression of $K_{\mu,\mu}$ of order at most
$\rho$. Let $c_\mu=0$ and $c_{j-1}=c_j+(\mu(\rho-1)+1)\tau$ for $j\in[\mu]$.

Suppose that $\mathcal M_0=(\mathcal L_i^0)_{i=1}^\mu$ is an ordered
$\rho$-multicover of $C_0$ with $\chi(G[C_0])>c_0$, and let $x_i$ be
the apex of $\mathcal L_i^0$. For $j=1,\ldots,\mu$, apply
\cref{lem:css-exact-paths} to $\mathcal M_{j-1}$ and $C_{j-1}$, with $N=\mu$
and $c'=c_j$. This gives $C_j\subseteq C_{j-1}$ with
$\chi(G[C_j])>c_j$, an ordered $\rho$-multicover
$\mathcal M_j=(\mathcal L_i^j)_{i=1}^\mu$ of $C_j$ contained in
$\mathcal M_{j-1}$, a vertex $z_j\in C_{j-1}\setminus C_j$, and induced
$z_j$--$x_i$ paths $P_{i,j}$ with $\rho$ edges, since containment keeps
the apices. They satisfy
$V(P_{i,j})\subseteq V(\mathcal L_i^{j-1})\cup\{z_j\}$ and, for
$W_j=C_j\cup V(\mathcal M_j)$,
\[
 \{z_j\}\perp W_j,\qquad
 V(P_{i,j})\perp W_j\setminus V(\mathcal L_i^j).
\]

We claim that $x_i\mapsto x_i$, $z_j\mapsto z_j$, and
$x_iz_j\mapsto P_{i,j}$ define an impression of $K_{\mu,\mu}$ with sides
$\{x_1,\ldots,x_\mu\}$ and $\{z_1,\ldots,z_\mu\}$. The apices are distinct
and pairwise nonadjacent, since each lies in the core of its cover in
$\mathcal M_0$. For $j<j'$, we have $z_{j'}\in C_{j'-1}\subseteq C_j$,
so $z_j$ and $z_{j'}$ are distinct and nonadjacent. Every apex lies in
$W_j$, so no apex equals or is adjacent to $z_j$. Finally, let
$i\ne i'$ and $j<j'$. Containment gives
\[
 V(P_{i',j'})\subseteq V(\mathcal L_{i'}^{j'-1})\cup\{z_{j'}\}
 \subseteq V(\mathcal L_{i'}^j)\cup C_j
 \subseteq W_j\setminus V(\mathcal L_i^j),
\]
so $V(P_{i,j})\perp V(P_{i',j'})$. This proves the claim. The
impression has order $\rho$, a contradiction. Hence
\cref{lem:multicover-control-reduction} applies with $N=\mu$ and
$c=c_0$.
\end{proof}

The excluded subdivisions needed by
\cref{lem:css-two-control,lem:exact-path-reduction} come from
\cref{lem:short-biclique-obstruction}.

\begin{corollary}\label{cor:spider-biclique}
Let $r\ge1$, $s\ge3$, and $k\ge0$, and let $m=(sr+2)(k+1)$. Then
$\cC_k(K_{1,s}^{(r)})$ excludes $K_{m,m}^{(q)}$ for every
$2\le q\le2r$.
\end{corollary}
\begin{proof}
By \cref{lem:short-biclique-obstruction} with $k+1$ in place of $k$
and the first inequality of \cref{lem:node-span-width},
\[
 \sigma_{K_{1,s}^{(r)}}(K_{m,m}^{(q)})\ge
 \branchpi_{K_{1,s}^{(r)}}(K_{m,m}^{(q)})\ge k+1.
\]
Since $\sigma_{K_{1,s}^{(r)}}$ does not increase under taking induced
subgraphs (\cref{lem:node-span-width}), no graph in
$\cC_k(K_{1,s}^{(r)})$ contains an induced copy of $K_{m,m}^{(q)}$.
\end{proof}

We use \cref{lem:css-two-control} for paths and, at radius three, for
arms of length at least three, and \cref{lem:exact-path-reduction} at
radii at least four. Exact path lengths are needed because
\cref{cor:spider-biclique} provides excluded subdivisions only for
$q\le2r$. In particular, \cref{lem:css-two-control} applies at no
radius to $\cC_k(K_{1,3}^{(2)})$.

Next, we clean ordered multicovers. This gives a one-step radius
reduction, \cref{lem:multicover-descent}, whose two conditions are what
the radius reductions below have to verify. It uses the following sets.

\begin{observation}\label{obs:multicover-levels}
Let $\rho\ge3$ and $\theta,N\ge1$, and let $(\mathcal L_i)_{i=1}^N$ be
an ordered $\rho$-multicover of $C$ in $G$, where
$\mathcal L_i=(L_{i,0},\ldots,L_{i,\rho-1})$ and $L_{i,0}=\{x_i\}$.
Define
\[
\begin{aligned}
 C'&=\{c\in C\mid\dist(c,x_i)\ge\rho\text{ for all }i\in[N]\},&
 Z_i&=\bigcup_{a=0}^{\rho-2}L_{i,a},\\
 B_i&=L_{i,\rho-1}\cap N(C'),&
 E_i&=\{b\in B_i\mid|\{j>i\mid N(b)\cap Z_j\ne\varnothing\}|\ge\theta\}.
\end{aligned}
\]
Then $B_i\perp\bigcup_{a=0}^{\rho-3}L_{j,a}$ whenever $i<j$.
\end{observation}
\begin{proof}
If $b\in B_i$ sees $L_{j,a}$ with $i<j$ and $a\le\rho-3$, then a
neighbor $c\in C'$ of $b$ satisfies $\dist(x_j,c)\le a+2\le\rho-1$,
a contradiction.
\end{proof}

\begin{lemma}\label{lem:multicover-descent}
Let $\rho\ge3$, and let $\cC$ be a $\rho$-controlled class.
Suppose that, for every $\tau\ge1$, there are integers $m,\theta\ge1$
and $c_{\mathrm{free}},c_{\mathrm{rich}}\ge0$ with the following
property. Let $N=\theta m$. For every $G\in\cC$ with
$\chi^{\rho-1}(G)\le\tau$ and every ordered $\rho$-multicover of length
$N$ in $G$, with $C'$, $Z_i$, $B_i$, and $E_i$ as in
\cref{obs:multicover-levels}, the following conditions hold.
\begin{enumerate}[label=(C\arabic*)]
 \item\label{item:clean-free} For every $I\in\binom{[N]}m$, every
 $C_0\subseteq C'$, and every choice of sets $P_i\subseteq B_i$,
 $i\in I$, if each $P_i$ covers $C_0$ and $P_i\perp Z_j$ for all
 distinct $i,j\in I$, then $\chi(G[C_0])\le c_{\mathrm{free}}$.
 \item\label{item:clean-rich} $\chi(G[C'\cap N(E_i)])\le c_{\mathrm{rich}}$ for every $i\in[N]$.
\end{enumerate}
Then $\cC$ is $(\rho-1)$-controlled.
\end{lemma}
\begin{proof}
Fix $\tau\ge1$, let $G\in\cC$ satisfy $\chi^{\rho-1}(G)\le\tau$, and
consider an ordered $\rho$-multicover of $C$ of length $N$ in $G$. By
\cref{lem:multicover-control-reduction}, it suffices to show that
\[
 \chi(G[C])\le N\tau+Nc_{\mathrm{rich}}+\binom{N}{m}c_{\mathrm{free}}.
\]
The set $C\setminus C'$ lies in the balls $B(x_i,\rho-1)$, which cost
at most $N\tau$ colors.
Let $C''=C'\setminus\bigcup_iN(E_i)$. Then
$\chi(G[C'\setminus C''])\le Nc_{\mathrm{rich}}$
by~\ref{item:clean-rich}.
For $y\in C''$, choose $p_i(y)\in B_i\cap N(y)$ for each $i$.
Since $p_i(y)\notin E_i$, at most $\theta-1$ later cores meet
$N(p_i(y))$. Repeatedly select the least remaining index and
discard it and all its later contacts. Each step removes at most
$\theta$ indices, so, as $N=\theta m$, at least $m$ indices are selected;
let $I(y)\in\binom{[N]}m$ consist of the first $m$ of them. Then
\[
 i<j,\quad i,j\in I(y)\quad\Longrightarrow\quad
 \{p_i(y)\}\perp Z_j.
\]
For $I\in\binom{[N]}m$, let $C_I=\{y\in C''\mid I(y)=I\}$.
For $i\in I$, let $P_i=\{p_i(y)\mid y\in C_I\}$. Each $P_i$ covers
$C_I$. The displayed relation and the cover order give
$P_i\perp Z_j$ for distinct $i,j\in I$, so condition~\ref{item:clean-free} gives
$\chi(G[C_I])\le c_{\mathrm{free}}$. Summing over the at most
$\binom Nm$ sets $I$ proves the bound.
\end{proof}

\subsubsection{Supports and projection}
\label{subsec:projection}

An \emph{$H$-support} is a finite nonempty family of induced
copies of $H$ whose intersection graph is connected. For a branch
decomposition $(T,\lambda)$ of $G$ and a family $\mathcal R$ of
subgraphs of $G$, let
\[
 S_T(\mathcal R)=\bigcup_{Q\in\mathcal R}T[Q].
\]
If $\mathcal R$ is a support, then $S_T(\mathcal R)$ is a subtree of
$T$, since copies with a common vertex have intersecting spans.
A support $\mathcal R$ \emph{joins} $x,y$ if $x,y\in V(\mathcal R)$.
Two supports with a common
vertex form a support together. Hence, if $K[U]$ is connected for an
auxiliary graph $K$, every $x\in U$ has a support containing $x$, and
every edge of a spanning tree of $K[U]$ has a support joining its ends,
then the union of all these supports is a support.

We use repeatedly the elementary bound
\begin{equation}\label{eq:chi-balls}
 \chi^R(G)\le\tau
 \quad\Longrightarrow\quad
 \chi(G[B(X,R)])\le |X|\tau.
\end{equation}
The following maximal-packing estimate will be used at cuts.

\begin{lemma}\label{lem:local-crossing-packing}
Let $H$ have $h\ge1$ vertices, let $R,\tau\ge1$ and $k,c\ge0$,
and let $(A,B)$ be a cut of $G$ with $\pi^\times_H(G;A,B)\le k$.
Suppose $\chi^R(G)\le\tau$, and let $Y\subseteq V(G)$. Suppose that
every $Y'\subseteq Y$ with $\chi(G[Y'])>c$ admits a crossing induced
copy of $H$ contained in $B(Y',R-1)$.
Then $\chi(G[Y])\le c+hk\tau$.
\end{lemma}
\begin{proof}
Choose a maximal crossing induced $H$-packing $\mathcal Q$, and let
$Y'=Y\setminus B(V(\mathcal Q),R)$. If $\chi(G[Y'])>c$, the hypothesis
gives a crossing copy in $B(Y',R-1)$. It has distance at least two from
$V(\mathcal Q)$, contradicting maximality. Since
$|V(\mathcal Q)|\le hk$, we have
$\chi(G[Y\cap B(V(\mathcal Q),R)])\le hk\tau$ by \eqref{eq:chi-balls},
and the bound follows.
\end{proof}

For a connected graph $H$ with at least two vertices and a branch
decomposition $(T,\lambda)$ of $G$, a nonempty subtree $S\subseteq T$
is \emph{reserved against}
$X\subseteq V(G)$ if
\[
 \forall x\in V(S)\ \exists Q_x\cong H:
 \quad x\in V(T[Q_x]),\quad V(Q_x)\perp X,
\]
where every $Q_x$ is induced in $G$. Then
$\lambda(X)\cap V(S)=\varnothing$, since a leaf of $T$ in a span
$T[Q]$ corresponds to a vertex of $Q$.

\begin{observation}\label{obs:reserved}
If $\mathcal R$ is an $H$-support with $V(\mathcal R)\perp X$, then
$S_T(\mathcal R)$ is reserved against $X$.
\end{observation}
\begin{proof}
Each node of the subtree $S_T(\mathcal R)$ lies in the span of a member
of $\mathcal R$.
\end{proof}

Let $S\subseteq T$ be a subtree and let $T'$ be a component of
$T-V(S)$. The edge of $T$ joining $T'$ to $S$ displays the
\emph{attachment cut} $(A,B)$ of $T'$, where $B=\lambda^{-1}(V(T'))$.
For $Y\subseteq V(G)$, the set $Y\cap B$ is the \emph{attachment fiber}
of $Y$ at $T'$.

\begin{lemma}[Projection]\label{lem:layered-projection}
Let $H$ be connected with at least two vertices, let
$R,\tau,q,b,k\ge1$, and let $(T,\lambda)$ be a node-width-$k$
decomposition of $G$.
Suppose $\chi^R(G)\le\tau$, and let $\varphi$ witness $R$-control
of $\cC_{k-1}(H)$. Let $X\subseteq V(K)\subseteq V(G)$, where
$E(G[X])\subseteq E(K[X])$, and take distance layers in each
component of $K$.
Suppose that, for each layer $\Lambda$, either $\chi(G[X\cap\Lambda])\le\tau$,
or $X\cap\Lambda$ has a partition into at most $q$ sets $Y$, each with
a subtree $S\subseteq T$ reserved against $Y$ such that every
attachment fiber of $Y$ has chromatic number at most $b$.
Then
\[
 \chi(G[X])\le2qb\varphi(\tau).
\]
\end{lemma}
\begin{proof}
We first show that $\chi(G[Y])\le b\varphi(\tau)$ for every set $Y$ as
in the statement. Color each attachment fiber of $Y$ with $[b]$ and
reuse the labels, obtaining $Y=Y_1\mathbin{\dot\cup}\cdots
\mathbin{\dot\cup}Y_b$. Fix an induced $H$-packing $\mathcal Q$ in
$G[Y_i]$ with a common span node $x$ in $T$.
If $x\in V(S)$, add $Q_x$, giving $|\mathcal Q|+1\le k$.
Otherwise $x$ lies in a component $T'$ of $T-V(S)$ attached to
$x'\in V(S)$; let $(A,B)$ be its attachment cut.
Every member of $\mathcal Q$ meets $B$, since a span of leaves
outside $T'$ avoids $T'$. No member is contained in $B$, because
$Y_i\cap B$ is independent and $H$ has an edge. Thus every span
contains $x'$, and adding $Q_{x'}$ again gives
$|\mathcal Q|\le k-1$.
In the restriction of $(T,\lambda)$ to $Y_i$, a node lies in the span
of a copy if and only if it lies in its span in $T$
(see~\cref{sec:prelim}), so $\sigma_H(G[Y_i])\le k-1$. Since
$B_{G[Y_i]}(v,R)\subseteq B_G(v,R)$, the control hypothesis
gives $\chi(G[Y_i])\le\varphi(\tau)$; empty sets cause no difficulty.

Since $\varphi(\tau)\ge\tau$, every $X\cap\Lambda$ has chromatic number
at most $c=qb\varphi(\tau)$.
If $(\Lambda_i)$ are the layers in one component, the map
\[
 v\in X\cap\Lambda_i\longmapsto(i\bmod2,c_i(v)),
 \qquad c_i\colon X\cap\Lambda_i\to[c]\text{ proper},
\]
is proper on $G[X]$. Indeed, a $G[X]$-edge is a $K$-edge and
joins equal or consecutive layers. Distinct components of $K$ are
anticomplete in $G[X]$, so the same palettes apply to all components.
\end{proof}

For a cut $(A,B)$, let
\[
 \partial_A^1=\{x\in A\mid \dist(x,B)\le1\},
 \qquad \partial_A^2=\{x\in A\mid \dist(x,B)\le2\},
\]
and define $\partial_B^1$ and $\partial_B^2$ symmetrically.
The next lemma bounds the second boundary of a cut in terms of its
first boundaries, provided that every set of large chromatic number at
distance two from the other side yields a crossing copy of $H$ with
exactly one vertex on that side.

\begin{lemma}\label{lem:second-boundary}
Let $H$ be connected with $h\ge2$ vertices, let $a,\tau\ge1$ and
$k,L\ge0$, suppose that $\chi^2(G)\le\tau$, and let $(A,B)$ be a cut
of $G$ with $\pi^\times_H(G;A,B)\le k$ and
$\chi(G[\partial_A^1]),\chi(G[\partial_B^1])\le a$.
Consider disjoint $D,C\subseteq A$ and $E\subseteq B$, where
$E=\{e_1,\ldots,e_\ell\}$ is independent, $E$ covers $D$, $D$ covers
$C$, and $E\perp C$. Assign each $p\in D$ to its first neighbor in
$E$, obtaining classes $D_1,\ldots,D_\ell$. For $C'\subseteq C$, call
$i$ \emph{active} for $C'$ if some induced copy $Q$ of $H$ satisfies
\[
 e_i\in V(Q)\subseteq\{e_i\}\cup(D_i\cap N(C'))\cup C'.
\]
Suppose that, for all such $D,C,E$, every $j\ge0$, and every
$C'\subseteq C\cap N_G\bigl(\bigcup_{i>j}D_i\bigr)$ with
$\chi(G[C'])>L$, some index $i>j$ is active for $C'$.
Then
\[
 \chi(G[\partial_A^2])\le a+a\bigl((k+1)L+hk\tau\bigr).
\]
\end{lemma}
\begin{proof}
Let $D^*=\partial_A^1$, $E^*=\partial_B^1$, and
$C^*=\partial_A^2\setminus D^*$. Choose a coloring $\kappa\colon E^*\to[a]$
and neighbor maps $f\colon D^*\to E^*$, $g\colon C^*\to D^*$. For each color
$\gamma\in[a]$, the sets
\[
 E=\kappa^{-1}(\gamma),\qquad D=f^{-1}(E),\qquad C=g^{-1}(D)
\]
satisfy the conditions in the statement, since $C^*$ has no neighbor
in $B$. These sets $C$ partition $C^*$, and $\chi(G[D^*])\le a$. It
therefore suffices to show $\chi(G[C])\le(k+1)L+hk\tau$ for all
$D,C,E$ as in the statement.

Suppose that $\chi(G[C])>(k+1)L+hk\tau$, and let $C_0=C$. For
$q=0,\ldots,k$, suppose that
$\chi(G[C_q])>(k+1-q)L+(k-q)h\tau$; in particular,
$\chi(G[C_q])>L$. Since $D$ covers $C$, the hypothesis with $j=0$
gives an index active for $C_q$. Let $j_q$ be the largest one,
witnessed by a copy $Q_q$, and let
\[
 C_{q+1}=C_q\setminus\Bigl(N_G\Bigl(\bigcup_{i>j_q}D_i\Bigr)\cup
 B(V(Q_q),2)\Bigr).
\]
An index active for a subset of $C_q$ is active for $C_q$. Hence, by
the hypothesis with $j=j_q$, the set
$C_q\cap N_G(\bigcup_{i>j_q}D_i)$ has chromatic number at most $L$.
With \eqref{eq:chi-balls}, this gives
$\chi(G[C_{q+1}])>(k-q)L+(k-q-1)h\tau$. Every vertex of $C_{q+1}$
has a neighbor in $D$, but none in $\bigcup_{i\ge j_q}D_i$, because
$N(D_{j_q})\cap C\subseteq B(e_{j_q},2)$. Since $\{e_i\}\perp C$ and
$H$ is connected, every copy witnessing an active index $i$ meets
$D_i\cap N(C')$. Hence every index active for $C_{q+1}$ is smaller
than $j_q$.

The copies $Q_0,\ldots,Q_k$ cross the cut. Let $q<q'$. The vertices
of $Q_{q'}$ in $C$ lie in $C_{q'}$, at distance greater than two from
$Q_q$. Those in $D$ have a neighbor in $C_{q'}$, and hence distance
at least two from $Q_q$. The vertex $e_{j_{q'}}$ misses
$e_{j_q}$ by independence, the vertices of $Q_q$ in $D_{j_q}$ because
$j_{q'}<j_q$ and they are assigned to their first neighbor, and the
vertices of $Q_q$ in $C$ because $E\perp C$. Thus the $k+1$ copies are
pairwise anticomplete, a contradiction.
\end{proof}

\subsubsection{Free and rich configurations}
\label{subsec:core-paths}

Let $r\ge2$. A \emph{configuration of depth $r$} consists of pairwise
disjoint sets $C,Z_1,\ldots,Z_m$, its \emph{target} and its
\emph{cores}, and additional parent sets disjoint from these, whose
elements are its \emph{parents}.
We require $Z_i\perp C$ and $Z_i\perp Z_j$ for $i\ne j$.
A vertex is \emph{noncore} if it lies in no core.
Each core $Z_i$ is the disjoint union of its \emph{entry} set $D_i$
and its \emph{inward} set $Z_i^-$, and every $d\in D_i$ starts a fixed
induced path
\[
 z_1(d)-\cdots-z_r(d),\qquad z_1(d)=d,\quad
 z_2(d),\ldots,z_r(d)\in Z_i^-.
\]
Let $Z^-=\bigcup_iZ_i^-$.
For $x\in C$, the parents adjacent to $x$ are the parents \emph{of
$x$}; they form the set $P(x)$, and $P(U)=\bigcup_{x\in U}P(x)$ for
$U\subseteq C$.

A \emph{free configuration of depth $r$} has disjoint parent sets
$P_1,\ldots,P_m$, each covering $C$, such that
\[
 N(p)\cap D_i\ne\varnothing\ (p\in P_i),\qquad
 P_i\perp Z_i^-,\qquad P_i\perp Z_j\ (i\ne j).
\]
The parent sets need not be independent. Define $\iota(v)=i$ for
$v\in P_i\cup Z_i$. Let $K$ be the graph on $C$ in which two distinct
vertices are adjacent if they are adjacent in $G$ or have a common
parent.
A \emph{rich configuration of depth $r$} instead has one parent set
$E$ covering $C$, with $E\perp Z_i^-$ for all $i$. Let
$G_E=G[E\cup C]$ and
\[
 \acc(p)=\{i\mid N(p)\cap D_i\ne\varnothing\}
 \qquad(p\in E).
\]
We say that $p$ \emph{accesses} the cores $Z_i$ with $i\in\acc(p)$,
and call $|\acc(p)|$ its number of \emph{accesses}.

For a parent $p$ and $d\in D_i\cap N(p)$, the \emph{core path}
$p-z_1(d)-\cdots-z_r(d)$ is induced and has $r$ edges.
In a rich configuration, an induced copy of a graph whose only noncore
vertex is a parent $p$ is a \emph{pure copy} at $p$; for example,
$s$ core paths at $p$ in distinct accessed cores induce a pure copy of
$K_{1,s}^{(r)}$.
Cleaned multicovers yield configurations.

\begin{lemma}\label{lem:multicover-configurations}
Let $r\ge2$, $\rho\ge r+1$, and $\theta,m,\tau\ge1$, and let $G$
satisfy $\chi^{\rho-1}(G)\le\tau$. Consider an ordered
$\rho$-multicover of length $N=\theta m$ in $G$, with the notation of
\cref{obs:multicover-levels}. Give each core $Z_i$ the entry set
$L_{i,\rho-2}$ and the inward set $\bigcup_{a=0}^{\rho-3}L_{i,a}$.
\begin{enumerate}[label=(\alph*)]
 \item\label{item:mc-free} If $I$, $C_0$, and $(P_i)_{i\in I}$ are as
 in condition~\ref{item:clean-free} of \cref{lem:multicover-descent},
 then they form, with the cores $Z_i$ for $i\in I$, a free
 configuration of depth $r$ with $m$ parent sets and target $C_0$.
 \item\label{item:mc-rich} For every $i\in[N]$ and every independent
 set $E\subseteq E_i$, the cores $Z_j$ with $j>i$ form, with the parent
 set $E$, a rich configuration of depth $r$ with target
 $C'\cap N(E)$, in which every parent accesses at least $\theta$
 cores. Moreover, $\chi(G[E_i])\le\tau$.
\end{enumerate}
\end{lemma}
\begin{proof}
For every $d\in L_{i,\rho-2}$, choose vertices $z_a\in L_{i,\rho-1-a}$,
$a\in[r]$, with $z_1=d$ and $z_az_{a+1}\in E(G)$; they exist because
each level covers the next one, and $\rho-1-r\ge0$. The path
$z_1-\cdots-z_r$ is induced, since nonconsecutive levels are
anticomplete. Distinct cores are anticomplete by the cover order, and
all cores are anticomplete to $C$.

For~\ref{item:mc-free}, every $P_i$ covers $C_0$, every vertex of
$P_i\subseteq L_{i,\rho-1}$ has a neighbor in $L_{i,\rho-2}$, $P_i$ is
anticomplete to the inward set of $Z_i$ because nonconsecutive levels
are anticomplete, and $P_i\perp Z_j$ for $j\ne i$ by assumption.

For~\ref{item:mc-rich}, $E$ covers $C'\cap N(E)$. By
\cref{obs:multicover-levels}, $E$ is anticomplete to
the inward sets of all later cores. Hence a parent that sees a later core
$Z_j$ sees its entry set $L_{j,\rho-2}$, and the definition of $E_i$
gives at least $\theta$ accesses. Finally,
$E_i\subseteq L_{i,\rho-1}\subseteq B(x_i,\rho-1)$, so
$\chi(G[E_i])\le\chi^{\rho-1}(G)\le\tau$.
\end{proof}

The free and rich target estimates therefore give a common radius
reduction.

\begin{lemma}\label{lem:configuration-descent}
Let $\cC$ be controlled at some finite radius, and let
$R\ge r\ge2$. Suppose that, for each $\tau\ge1$, there are integers
$m,\theta\ge1$ and $c_{\mathrm{free}},c_{\mathrm{rich}}\ge0$ such that every $G\in\cC$ with
$\chi^R(G)\le\tau$ satisfies the following.
\begin{enumerate}[label=(\roman*)]
 \item\label{item:descent-free} Every free configuration of depth $r$
 with at least $m$ parent sets has target chromatic number at most
 $c_{\mathrm{free}}$.
 \item\label{item:descent-rich} Every rich configuration of depth $r$
 whose parent set is independent and whose parents each access at
 least $\theta$ cores has target chromatic number at most
 $c_{\mathrm{rich}}$.
\end{enumerate}
Then $\cC$ is $R$-controlled.
\end{lemma}
\begin{proof}
Suppose that $\cC$ is $\rho$-controlled, where $\rho>R$. Let
$G\in\cC$ satisfy $\chi^{\rho-1}(G)\le\tau$, so that
$\chi^R(G)\le\tau$, and consider an ordered $\rho$-multicover of
length $\theta m$ in $G$. Since $\rho\ge r+1$,
\cref{lem:multicover-configurations} applies. Part~\ref{item:mc-free}
of that lemma together with hypothesis~\ref{item:descent-free} gives
condition~\ref{item:clean-free} of \cref{lem:multicover-descent}. By part~\ref{item:mc-rich}, each $E_i$ splits
into at most $\tau$ independent sets $E$, and the targets $C'\cap N(E)$
of their rich configurations have union $C'\cap N(E_i)$;
by~\ref{item:descent-rich},
condition~\ref{item:clean-rich} holds with $\tau c_{\mathrm{rich}}$ in
place of $c_{\mathrm{rich}}$.
By \cref{lem:multicover-descent}, $\cC$ is $(\rho-1)$-controlled.
Descending from the given finite radius proves $R$-control; if that
radius is already at most $R$, use the monotonicity of control noted
in \cref{subsec:local-control}.
\end{proof}

The target bounds in
\cref{subsec:long-arms,subsec:arms-two,subsec:broom-three} all start
from reserved subtrees, which are given by supports of the following
kind. Let $H$ be connected with at least two vertices. A free
configuration with target $C$ is \emph{$H$-supported} if every $x\in C$
has an $H$-support containing $\{x\}\cup P(x)$ with all noncore vertices
in $\{x\}\cup P(x)$, and every edge $xy$ of $K$ has an $H$-support
joining $x,y$ with all noncore vertices in $\{x,y\}\cup P(x)\cup P(y)$.
A rich configuration with target $C$ and independent parent set $E$ is
\emph{$H$-supported} if every $p\in E$ lies in a pure copy of $H$, and
every induced path of length two or three with ends in $E$ and interior
in $C$ has an $H$-support joining its ends, each member of which
contains an end of the path and has all noncore vertices on the path.

\begin{lemma}\label{lem:layered-supports}
Let $H$ be connected with at least two vertices, and let $j\ge2$.
\begin{enumerate}[label=(\alph*)]
 \item\label{item:layered-free} Consider an $H$-supported free
 configuration. Take distance layers $(\Lambda_i)$ in a component of
 $K$, and let $U=\bigcup_{i\le j-2}\Lambda_i$. Then the union
 $\mathcal R$ of the supports of the vertices and edges of $K[U]$ is an
 $H$-support with $U\cup P(U)\subseteq V(\mathcal R)$, all noncore
 vertices in $U\cup P(U)$, and $V(\mathcal R)\perp\Lambda_j$.
 \item\label{item:layered-rich} Consider an $H$-supported rich
 configuration with target $C$ and parent set $E$. Take distance layers
 $(\Lambda_i)$ in a component of $G_E$ rooted at $z\in E$, and let
 $U=\bigcup_{i\le j-2}\Lambda_i$. Then the family $\mathcal R$ of all
 induced copies of $H$ that meet $E\cap U$ and have all noncore vertices
 in $U$ is an $H$-support with $E\cap U\subseteq V(\mathcal R)$ and
 $V(\mathcal R)\perp C\cap\Lambda_j$.
\end{enumerate}
\end{lemma}
\begin{proof}
For~\ref{item:layered-free}, the graph $K[U]$ is connected, so
$\mathcal R$ is a support, and it contains $U\cup P(U)$. Its
noncore vertices lie in $U\cup P(U)$. The set $U$ misses $\Lambda_j$
by the layer gap in $K$, and so does $P(U)$, since a parent of $u\in U$
adjacent to $w\in\Lambda_j$ would make $uw$ an edge of $K$. All cores
miss $C\supseteq\Lambda_j$.

For~\ref{item:layered-rich}, let $E'=E\cap U$. The family $\mathcal R$
contains, for each $p\in E'$, a pure copy at $p$. We first show that the
graph on $E'$ in which two vertices are adjacent when they are the ends of an
induced path of $G$ of length two or three with interior in $C\cap U$
is connected. Let $e\in E'$, and take a shortest $z$--$e$ path
$v_0-\cdots-v_a$ in $G_E$; it lies in $U$. Let $p_b=v_b$ if
$v_b\in E$, and otherwise let $p_b$ be a vertex of $E$ adjacent to
$v_b$. Since $v_b\in C$ implies $b<a$, the vertex $p_b$ lies in a
layer of index at most $b+1\le a$, and hence in $U$. If
$p_b\ne p_{b+1}$, the walk through $p_b,v_b,v_{b+1},p_{b+1}$, with
repeated vertices omitted, has length at most three and interior in
$\{v_b,v_{b+1}\}\cap C\subseteq C\cap U$. A shortest $p_b$--$p_{b+1}$
path in the subgraph induced by its vertices is an induced path of
$G$; since $E$ is independent, it has length two or three. Thus
$p_0=z$ and $p_a=e$ are connected in this graph.
For each of these paths, the members of its support belong to
$\mathcal R$ and join the pure copies at its ends. Every other member of $\mathcal R$ contains some
$p\in E'$ and thus meets the pure copy at $p$. Hence $\mathcal R$ is a
support. Its noncore vertices lie in $U$, which misses $\Lambda_j$ by
the layer gap in $G_E=G[E\cup C]$, and its cores miss $C$.
\end{proof}

In both cases of \cref{lem:layered-supports}, \cref{obs:reserved} shows
that $S_T(\mathcal R)$ is reserved against $\Lambda_j$, respectively
$C\cap\Lambda_j$, for every branch decomposition $(T,\lambda)$ of $G$.

\subsection{Paths and stars}
\label{subsec:paths-stars}

Part~(\ref{item:common-paths}) of \cref{thm:common-two-control}
follows from part~(\ref{item:common-long-arms}): for
$r=\max\{3,\lceil(t-1)/2\rceil\}$, every edge of $K_{1,3}^{(r)}$
belongs to an induced copy of $P_t$, so
\cref{lem:edge-extension,lem:node-span-width} give
$\cC_k(P_t)\subseteq\cC_{\lfloor3k/2\rfloor}(K_{1,3}^{(r)})$. We
include a direct proof, since it is much simpler.

\begin{proof}[Proof of \cref{thm:common-two-control}(\ref{item:common-paths})]
Fix $t\ge2$. By \cref{lem:control-start}, $\cC_0(P_t)$ is
$2$-controlled. Let $k\ge1$; then $\cC_k(P_t)$ is controlled at some
finite radius $\rho\ge2$. Let $r=\rho+t$ and $m=(2r+2)(k+1)$.
Since $t\le2r+1$, every edge of $K_{1,2}^{(r)}=P_{2r+1}$ belongs to
an induced copy of $P_t$. For $2\le q\le\rho+3$, we have $q\le2r$, so
\cref{lem:node-span-width,lem:edge-extension,lem:short-biclique-obstruction}
give
\[
 \sigma_{P_t}(K_{m,m}^{(q)})\ge\branchpi_{P_t}(K_{m,m}^{(q)})
 \ge\branchpi_{K_{1,2}^{(r)}}(K_{m,m}^{(q)})\ge k+1.
\]
Thus $\cC_k(P_t)$ excludes $K_{m,m}^{(q)}$ for $2\le q\le\rho+3$, and
\cref{lem:css-two-control} shows that it is $2$-controlled.
\end{proof}

For stars, both boundaries of every displayed cut are controlled
directly. We first partition a covered set into sets with private
neighbors.

\begin{lemma}\label{lem:private-neighbor-partition}
Let $s\ge3$, $k\ge0$, and $\Delta,\tau\ge1$, and suppose that
$\sigma_{K_{1,s}}(G)\le k$ and $\omega(G)\le\tau$.
Let $D,C\subseteq V(G)$ be disjoint, with $D$ covering $C$ and
$|N(d)\cap C|\le\Delta$ for every $d\in D$. Then $C$ has a
partition into at most $2\Delta-1$ sets $Z$ satisfying
\[
 \sigma_{K_{1,s-1}}(G[Z])\le \binom{\tau+k}{\tau}-1.
\]
\end{lemma}
\begin{proof}
Choose $p(x)\in D\cap N(x)$ for $x\in C$. Join distinct $x,y$
in an auxiliary graph when $yp(x)\in E(G)$ or $xp(y)\in E(G)$.
For every $Y\subseteq C$, the subgraph of the auxiliary graph induced
by $Y$ has at most $(\Delta-1)|Y|$ edges, since each $x\in Y$ has at
most $\Delta-1$ neighbors $y\ne x$ with $yp(x)\in E(G)$. Hence the
auxiliary graph is $(2\Delta-2)$-degenerate. In each of its
$2\Delta-1$ color classes $Z$,
\[
 x\ne y\Longrightarrow p(x)\ne p(y),\qquad N(p(x))\cap Z=\{x\}.
\]
Fix a branch decomposition $\mathcal D=(T,\lambda)$ of $G$ with
$\sigma_{K_{1,s}}(G,\mathcal D)\le k$. In its restriction to $Z$, a
node lies in the span of a copy if and only if it lies in its span in
$T$ (see~\cref{sec:prelim}), so it suffices to show that no
$\binom{\tau+k}{\tau}$ pairwise anticomplete induced copies of
$K_{1,s-1}$ in $G[Z]$ have spans in $T$ with a common node. Suppose
otherwise. Add to each copy the private parent $p(x)$ of its center
$x$; this gives an induced copy of $K_{1,s}$ and only enlarges the
span. By Ramsey's theorem, $k+1$ of the added parents are pairwise
nonadjacent, since $\omega(G)\le\tau$. The corresponding $k+1$
copies of $K_{1,s}$ are pairwise anticomplete and have spans with a
common node, a contradiction.
\end{proof}

\begin{lemma}\label{lem:star-boundaries}
Let $s\ge3$, and let $g$ be nondecreasing in both arguments. Suppose
that, for every graph $G'$ and every $k'\ge0$, $\tau'\ge1$,
\begin{equation}\label{eq:smaller-star-control}
 \sigma_{K_{1,s-1}}(G')\le k',\quad\chi^2(G')\le\tau'
 \quad\Longrightarrow\quad\chi(G')\le g(k',\tau').
\end{equation}
For $k\ge0$ and $\tau\ge1$, define
\[
\begin{gathered}
 h=s+1,\quad \Delta=(s-2)\tau,\quad \hat k=\binom{\tau+k}{\tau}-1,\quad
 M=(2\Delta-1)g(\hat k,\tau),\\
 a=(hk+2s-3)\tau+1,\qquad b_s(k,\tau)=a\bigl(1+(k+1)M+hk\tau\bigr).
\end{gathered}
\]
If $\sigma_{K_{1,s}}(G)\le k$, $\chi^2(G)\le\tau$, and
$(A,B)$ is a cut with $\pi_{K_{1,s}}^\times(G;A,B)\le k$, then
\[
 \chi(G[\partial_A^2]),\chi(G[\partial_B^2])\le b_s(k,\tau).
\]
\end{lemma}
\begin{proof}
For the first boundaries, apply \cref{lem:local-crossing-packing}
to $Y=\partial_A^1$ with $R=2$ and $c=(2s-3)\tau+1$. Let
$Y'\subseteq Y$ satisfy $\chi(G[Y'])>c$. Then some $v\in Y'$ has more
than $(2s-3)\tau$ neighbors in $Y'$. Choose $w\in B\cap N(v)$. Since
$\chi(G[N(v)])\le\tau$, either $(N(v)\cap Y')\setminus N(w)$
contains an independent $(s-1)$-set, giving a crossing star centered
at $v$ with a leaf $w$, or $N(v)\cap N(w)\cap Y'$ contains an
independent $s$-set, giving a crossing star centered at $w$. Both
stars lie in $B(Y',1)$. Hence $\chi(G[\partial_A^1])\le c+hk\tau=a$,
and the bound for $\partial_B^1$ follows by symmetry.

Second, let disjoint $P,Z\subseteq V(G)$ satisfy that $P$ covers $Z$
and $\alpha(G[N(p)\cap Z])\le s-2$ for every $p\in P$. Then
$|N(p)\cap Z|\le\Delta$. By \cref{lem:private-neighbor-partition}, $Z$
splits into at most $2\Delta-1$ sets $Z'$ with
$\sigma_{K_{1,s-1}}(G[Z'])\le\hat k$. Consequently $\chi(G[Z])\le M$.

Third, we verify the hypothesis of \cref{lem:second-boundary} with
$L=M$. Let $D,C,E$ be sets as there, let $j\ge0$, and let
$C'\subseteq C\cap N(\bigcup_{i>j}D_i)$ satisfy $\chi(G[C'])>M$.
The set $P=N(C')\cap\bigcup_{i>j}D_i$ covers $C'$, so, by the second
step, some $p\in P$ has an independent $(s-1)$-set
$Y\subseteq N(p)\cap C'$. If $p\in D_i$, then $i>j$ and
$G[\{e_i,p\}\cup Y]\cong K_{1,s}$, because $E\perp C$; hence $i$
is active for $C'$. Thus \cref{lem:second-boundary} gives
$\chi(G[\partial_A^2])\le a+a((k+1)M+hk\tau)=b_s(k,\tau)$, and the bound for
$\partial_B^2$ follows by symmetry.
\end{proof}

\begin{proof}[Proof of \cref{thm:common-two-control}(\ref{item:common-stars})]
We induct on $s\ge3$. Suppose that $\cC_k(K_{1,s-1})$ is
$2$-controlled for every $k$; for $s=3$ this is
part~(\ref{item:common-paths}), since $K_{1,2}=P_3$. Let $\varphi'_k$
witness this. Replacing $\varphi'_k$ by $\max_{i\le k}\varphi'_i$,
we may assume that $\varphi'_k$ is nondecreasing in $k$. Then the
function $g(k',\tau')=\varphi'_{k'}(\tau')$ satisfies
\eqref{eq:smaller-star-control} and is nondecreasing in both
arguments, so the function $b_s$ of \cref{lem:star-boundaries} is
nondecreasing in $k$ and $\tau$.

Induct on $k$, with $\varphi_k(0)=0$ and, for $\tau\ge1$,
\[
 \varphi_0(\tau)=(s-1)\tau+1,\qquad
 \varphi_k(\tau)=\max\{(2s-3)\tau+2,
                       2b_s(k,\tau)\varphi_{k-1}(\tau)\}.
\]
If $G\in\cC_0(K_{1,s})$ and $\chi^2(G)\le\tau$, then
$\alpha(G[N(v)])\le s-1$ and $\deg(v)\le(s-1)\tau$ for every
$v$, giving the base case.
For $k\ge1$, take $G\in\cC_k(K_{1,s})$ with $\chi^2(G)\le\tau$.
We may assume that $\chi(G)>(2s-3)\tau+2$ and pass to a vertex-critical
induced subgraph of the same chromatic number. It is connected and has
minimum degree at least $\chi(G)-1>(2s-3)\tau+1$.
Fix a node-width-$k$ decomposition $(T,\lambda)$ and distance layers
$\Lambda_i$ in $G$ from a vertex $z$. For $j\ge3$, let
\[
 U=\bigcup_{i\le j-3}\Lambda_i,\qquad U^+=\bigcup_{i\le j-2}\Lambda_i.
\]
For every $u\in V(G)$ and $v\in N(u)$, there is a $K_{1,s}$-support
contained in $N[u]$ and joining $u,v$. Indeed, either
$N(u)\setminus N[v]$ contains $s-1$ independent vertices, giving a
star centered at $u$ with leaf $v$, or this set has size at most
$(s-2)\tau$, since $\chi(G[N(u)])\le\tau$. In the latter case,
$\deg(u)>(2s-3)\tau+1$ gives $|N(u)\cap N(v)|>(s-1)\tau$, and an
independent $s$-set there gives two intersecting stars centered at
$u$ and $v$. For each $u\in U$ and $v\in N(u)$, take such a support.
Their union $\mathcal R$ is a support: the
families for one $u$ meet at $u$, and the edges of the connected graph
$G[U]$ join these families. Since
\[
 \lambda(U^+)\subseteq V(S_T(\mathcal R)),\qquad
 V(\mathcal R)\subseteq U^+\perp\Lambda_j,
\]
\cref{obs:reserved} shows that $S_T(\mathcal R)$ is reserved against
$\Lambda_j$. Each attachment fiber of $\Lambda_j$ is contained in the
second boundary of its attachment cut, because its vertices are at
distance at most two from $U^+$. Its chromatic number is at most
$b_s(k,\tau)$ by \cref{lem:star-boundaries,lem:node-span-width}. The first three layers
lie in $B(z,2)$. \Cref{lem:layered-projection}, with $X=V(G)$,
$K=G$, $q=1$, $R=2$, and $b=b_s(k,\tau)$, gives
\[
 \chi(G)\le2b_s(k,\tau)\varphi_{k-1}(\tau)\le\varphi_k(\tau).
\]
This completes both inductions.
\end{proof}

\subsection{Arms of length at least three}
\label{subsec:long-arms}

The proof of \cref{thm:common-two-control}(\ref{item:common-long-arms})
lowers the control radius in three steps: from the finite radius of
\cref{lem:control-start} to $\max\{r+1,6\}$ by
\cref{lem:configuration-descent} with \cref{lem:spider-target-bounds},
then to $3$ by \cref{lem:exact-path-reduction}, and finally to $2$ by
\cref{lem:css-two-control}.

We use the configurations from \cref{subsec:core-paths}. The next two
results specify both the vertices of the supports and the core indices
they use; these restrictions will ensure anticompleteness to previously
selected copies.
Copies are built by completing short paths through cores. For a
parent $p$, an entry $d\in D_i\cap N(p)$, and a path
$\Pi=v_0-\cdots-v_a$ with $v_a=p$, where $0\le a\le r$ and
$V(\Pi)\cap Z_i=\varnothing$, the \emph{completion} of $\Pi$ through
$d$ is
\[
 \Pi^+=v_0-\cdots-v_a-z_1(d)-\cdots-z_{r-a}(d);
\]
for $a=r$, no vertex is added. For initial paths
$\Pi_1,\ldots,\Pi_s$ with a common first vertex, write
\[
 J(\Pi_1,\ldots,\Pi_s)=G\left[\bigcup_{a=1}^sV(\Pi_a^+)\right],
\]
with the specified distinct completion cores. Each construction
specifies these cores and checks inducedness.

\begin{lemma}\label{lem:free-short-support}
Let $s\ge3$, $r\ge2$, and $m,\tau\ge1$. Suppose that
$\chi^1(G)\le\tau$ and that $G$ contains a free configuration of
depth $r$, with target $C$ and parent sets $P_1,\ldots,P_m$.
Let $I\subseteq[m]$ satisfy $m-|I|\ge(3s-5)\tau+1$, and call the
indices outside $I$ available. The following $K_{1,s}^{(r)}$-supports
exist.
\begin{enumerate}[label=(\arabic*)]
 \item\label{item:free-local-support}
 Each $x\in C$ has a support containing $x$, every available parent
 of $x$, and all entries adjacent to those parents. Its noncore
 vertices lie in $\{x\}$ or are available parents of $x$.
 \item\label{item:free-edge-support}
 If $r\ge3$ or $I=\varnothing$, every edge $xy\in E(K)$ has a
 support joining $x,y$. Its noncore vertices are $x,y$, available
 parents of $x$ or $y$, or one common parent of $x,y$.
\end{enumerate}
All cores used have available indices.
\end{lemma}
\begin{proof}
Let $P'(x)$ be the set of available parents of $x$.
Representatives from $q\tau$ parent sets in $N_G(x)$ contain $q$
independent vertices. Thus, after excluding one index, we can choose
$3s-5$ independent available parents with distinct indices.
For $Q=\{q_1,\ldots,q_{s-1}\}$, write
\[
 J_c(Q;\Pi)=J(c-q_1,\ldots,c-q_{s-1},\Pi),
\]
completing terminal parents through their own cores.

For (1), fix $p\in P'(x)$ and an entry $e$ of $p$.
Choose an independent $(2s-3)$-set $Q\subseteq P'(x)$ with distinct
indices different from $\iota(p)$. If $p$ is complete to an $(s-1)$-set
$Q'\subseteq Q$, choose $b\in Q\setminus Q'$ and use
\[
 \{J_p(Q';p),\ J_x(Q';x-b)\},
\]
with entry $e$ for the own-core arm at $p$.
Otherwise choose $Q'\in\binom{Q\setminus N_G(p)}{s-1}$ and use
$J_x(Q';x-p)$, again through $e$.
The first two copies intersect in $Q'$; the second construction is
one copy. Each family contains $x,p,e$. Their union over $(p,e)$
is the required support $\mathcal R_x$.

For (2), a common available parent joins $\mathcal R_x,\mathcal R_y$.
Otherwise
\[
 \{y\}\perp P'(x),\qquad \{x\}\perp P'(y).
\]
It suffices to construct one copy meeting both local supports.
If $xy\in E(G)$, choose $b\in P'(y)$ and an independent
$(2s-3)$-set $Q\subseteq P'(x)$ avoiding $\iota(b)$.
Since $|Q|=2s-3$, some $(s-1)$-set $Q'\subseteq Q$ lies in $N_G(b)$
or is anticomplete to $b$. Use, respectively,
\[
 J_b(Q';b)\qquad\text{or}\qquad J_x(Q';x-y-b).
\]
The latter initial path has length two.

If $xy\notin E(G)$, choose a common parent $p$; then
$\iota(p)\in I$, so $r\ge3$.
Suppose first that $p$ sees an available parent of one endpoint;
by symmetry take $a\in P'(x)\cap N_G(p)$.
Choose an independent $(3s-5)$-set $Q\subseteq P'(y)$ avoiding
$\iota(a)$. If $a$ is complete to an $(s-1)$-set $Q'\subseteq Q$, use
$J_a(Q';a)$. Otherwise $Q_0=Q\setminus N_G(a)$ has size at least
$2s-3$. Choose an $(s-1)$-set $Q'\subseteq Q_0$ wholly adjacent or
wholly nonadjacent to $p$, and use, respectively,
\[
 J_p(Q';p-a)\qquad\text{or}\qquad J_y(Q';y-p-a).
\]
Finally, if $p$ misses $P'(x)\cup P'(y)$, choose $b\in P'(y)$
and an independent $(2s-3)$-set $Q\subseteq P'(x)$ avoiding
$\iota(b)$. For an $(s-1)$-set $Q'\subseteq Q$ wholly adjacent or
wholly nonadjacent to $b$, use, respectively,
\[
 J_b(Q';b)\qquad\text{or}\qquad J_x(Q';x-p-y-b).
\]
All initial arms have length at most three. The selected independent
sets and the stated contact conditions exclude noncore chords.
Completion cores are distinct, so the anticompleteness conditions in
the definition of a free configuration exclude every remaining chord. In the last cases no core
has index $\iota(p)$. Every displayed copy meets both local supports,
so adjoining them proves (2).
\end{proof}

\begin{lemma}\label{lem:rich-short-support}
Let $s\ge3$ and $r\ge2$, and consider a rich configuration of depth $r$
with target $C$ and independent parent set $E$. Let $I$ be a set of core
indices such that every parent accesses at least $|I|+2s-1$ cores. Then
every induced path $\Pi$ of length $q\in\{2,3,4\}$ with ends $u,v\in E$
and interior in $C$ has a $K_{1,s}^{(r)}$-support joining $u,v$ whose
members are centered at $u$ or $v$, use at most $s$ cores, none with
index in $I$, and have all noncore vertices on $\Pi$. If $q>r$, this
support can consist of the two copies centered at $u$ and $v$ whose path
arms are the first $r$ edges of $\Pi$ from their centers. In particular,
the configuration is $K_{1,s}^{(r)}$-supported.
\end{lemma}
\begin{proof}
Every parent has at least $2s-1\ge s$ accesses, so $s$ core paths at it
in distinct accessed cores induce a pure copy of $K_{1,s}^{(r)}$; hence
the last assertion follows from the first.
Let $\Pi=p_0-\cdots-p_q$, where $p_0=u$ and $p_q=v$, and let
$A_x=\acc(x)\setminus I$; thus $|A_x|\ge2s-1$.
For $q>r$, use the initial $r$-edge segment at each endpoint and
$s-1$ core arms there. The path vertices other than $u,v$ belong to
$C$, so each copy is induced. Their path arms intersect because
$q\le4\le2r$.

Suppose $q\le r$, and let
$A_{uv}=\{i\notin I\mid D_i\cap N(u)\cap N(v)\ne\varnothing\}$.
If $|A_{uv}|\ge s$, pure copies at $u,v$ using $s$ common entries intersect.
Otherwise choose $s-1$ indices $I_u\subseteq A_u\setminus A_{uv}$ and
entries at $u$ in those cores. These entries miss $v$, and their
inward paths miss $E$. If $q=r$, use $\Pi$ as the remaining arm.
If $q<r$, choose
\[
 i_v\in A_v\setminus(A_{uv}\cup I_u),\qquad
 |A_v\setminus(A_{uv}\cup I_u)|\ge(2s-1)-2(s-1)=1,
\]
and complete $\Pi$ through an entry at $v$ in core $i_v$.
That entry misses $u$. In either case the $s$ arms induce $K_{1,s}^{(r)}$:
different cores are anticomplete and miss $C$, and the endpoint and
inward-set conditions exclude all remaining chords.
\end{proof}

In rich configurations, crossing copies are also obtained by extending
a parent into a target set of large chromatic number.

\begin{lemma}\label{lem:rooted-path-extension-local}
Let $r,\tau\ge1$, let $F$ be connected and induced in $G$, and let $v\notin V(F)$
have a neighbor in $F$. If $\chi^1(G)\le\tau$ and
$\chi(F)>(r-1)\tau$, then an induced $r$-edge path starts at $v$
and has all other vertices in $F$.
\end{lemma}
\begin{proof}
Induct on $r$, the case $r=1$ being immediate.
For $r\ge2$, some component $F'$ of $F-N(v)$ has
$\chi(F')>(r-2)\tau$. Connectedness gives
$x\in N(v)\cap V(F)$ with a neighbor in $F'$.
Prepend $v$ to the path starting at $x$ given by induction in
$G[V(F')\cup\{x\}]$. Since $v$ misses $F'$, the result is induced.
\end{proof}

\begin{lemma}\label{lem:spider-target-bounds}
Let $s,r\ge3$, $k,\tau\ge1$, $h=sr+1$, and $R=\max\{r+1,6\}$.
Suppose that $\cC_{k-1}(K_{1,s}^{(r)})$ is $R$-controlled by
$\varphi$. Let $G\in\cC_k(K_{1,s}^{(r)})$ satisfy $\chi^R(G)\le\tau$.
For configurations of depth $r$ in $G$, the following hold.
\begin{enumerate}[label=(\alph*)]
 \item\label{item:spider-free-target}
 A free configuration with $m\ge hk+(3s-5)\tau+1$ parent sets
 satisfies $\chi(G[C])\le2hk\tau\varphi(\tau)$.
 \item\label{item:spider-rich-target}
 A rich configuration with an independent parent set and at least
 $sk+2s-1$ accesses per parent satisfies
 $\chi(G[C])\le2(hk+r-1)\tau\varphi(\tau)$.
\end{enumerate}
\end{lemma}
\begin{proof}
Fix a node-width-$k$ decomposition $(T,\lambda)$.

For~\ref{item:spider-free-target}, take distance layers $(\Lambda_i)$
in a component of $K$ rooted at $x\in C$.
For $j\ge2$, let $W=\Lambda_j$ and $U=\bigcup_{i\le j-2}\Lambda_i$.
By \cref{lem:free-short-support} with $I=\varnothing$ and
\cref{lem:layered-supports}\ref{item:layered-free}, there is a support
$\mathcal R$ such that $S=S_T(\mathcal R)$ contains $\lambda(U)$ and is
reserved against $W$.

For an attachment cut $(A,B)$, let $Y=W\cap B$. Choose a maximal
crossing packing $\mathcal Q$ of copies contained in
$C\cup\bigcup_i(P_i\cup Z_i)$, and let
\[
 V_0=V(\mathcal Q),\qquad
 I=\{\iota(v)\mid v\in V_0\cap\bigcup_i(P_i\cup Z_i)\}.
\]
Then $|V_0|\le hk$, $|I|\le hk$, and at least
$(3s-5)\tau+1$ indices remain available.
If $y\in Y\setminus B_G(V_0,R)$, take the last three vertices
$u,z',y$ of a shortest $x$--$y$ path in $K$. Then $u\in U$,
$\dist_G(y,u)\le4$, and $\dist_G(y,z')\le2$.
Apply \cref{lem:free-short-support}\ref{item:free-edge-support} to
$uz',z'y$, avoiding $I$.
The combined support joins opposite sides of the attachment cut,
so one member $Q$ crosses it. Every noncore vertex of $Q$ belongs
to $B_G(y,5)$ and hence misses $V_0$, since $R\ge6$.
Every core vertex has an index outside $I$ and misses old targets,
cores, and parents by the definition of a free configuration. Thus
$V(Q)\perp V_0$, contradicting
maximality. We obtain
\[
 Y\subseteq B_G(V_0,R),\qquad \chi(G[Y])\le hk\tau.
\]
The first two layers lie in a radius-two ball. Apply
\cref{lem:layered-projection} with $X=C$, the graph $K$ of the
configuration, $q=1$, and $b=hk\tau$ to obtain~\ref{item:spider-free-target}.

For~\ref{item:spider-rich-target}, take distance layers $(\Lambda_i)$
in a component of $G_E$ rooted at $z\in E$. For $j\ge3$,
let
\[
 W=C\cap\Lambda_j,\qquad U=\bigcup_{i\le j-2}\Lambda_i,\qquad E'=E\cap U.
\]
By \cref{lem:rich-short-support} with $I=\varnothing$, the
configuration is $K_{1,s}^{(r)}$-supported, so, by
\cref{lem:layered-supports}\ref{item:layered-rich}, there is a support
$\mathcal R$ such that $S=S_T(\mathcal R)$ contains $\lambda(E')$ and is
reserved against $W$.

For an attachment cut $(A,B)$, let $Y=W\cap B$. Choose a maximal
crossing packing $\mathcal Q$ of copies contained in
$E\cup C\cup\bigcup_iZ_i$ that use at most $s$ cores. Let
\[
 V_0=V(\mathcal Q),\quad I=\{i\mid Z_i\cap V_0\ne\varnothing\},
 \quad Y'=Y\setminus B_G(V_0,R).
\]
Then $|V_0|\le hk$, $|I|\le sk$, and every parent retains at least
$2s-1$ accesses. If $\chi(G[Y])>(hk+r-1)\tau$, some component
$F$ of $G[Y']$ has $\chi(F)>(r-1)\tau$.
Inward vertices in cores outside $I$ already miss $V_0$: they miss
$E$ and $C$ since $E\perp Z_i^-$ and $Z_i\perp C$, and all old cores
by their different indices.

If $e\in E\setminus B$ sees $F$, \cref{lem:rooted-path-extension-local}
gives an induced path $e-x_1-\cdots-x_r$ with $x_i\in V(F)$.
Add $s-1$ core arms at $e$, using distinct cores outside $I$.
The resulting $K_{1,s}^{(r)}$ crosses the cut. All its vertices
outside $Z^-$ are in $B_G(x_1,r)$, so they miss $V_0$ because
$\dist_G(x_1,V_0)>R\ge r+1$.

Otherwise every parent seeing $F$ lies in $B$.
Choose $y\in V(F)$, a parent $v$ of $y$, and the vertex
$a\in\Lambda_{j-3}$ on a shortest $z$--$y$ path in $G_E$.
Let $u=a$ if $a\in E$, and otherwise choose a parent $u$ of $a$.
Then $u\in E'\subseteq A$, and the resulting
$u$--$v$ walk has length at most five, with vertices in $B_G(y,4)$.
Shorten it and split at internal parents. If a segment has length
five, it is the entire induced path $u-c_1-c_2-c_3-y-v$.
Choose a parent $w$ of $c_2$, and shorten
$u-c_1-c_2-w$ and $w-c_2-c_3-y-v$ separately.
Inducedness gives $w\notin\{u,v\}$, and $\dist_G(w,y)\le3$.
We now have segments of lengths two, three, or four, with ends in
$E$, interior in $C$, and all vertices in $B_G(y,4)$.
Their supports from \cref{lem:rich-short-support}, avoiding $I$,
join $u,v$. A member $Q$ crosses the cut, uses at most $s$ cores,
and has all vertices outside $Z^-$ in $B_G(y,5)$. Since
$R\ge6$,
these vertices miss $V_0$; its inward vertices were checked above.

In both cases, the new copy is anticomplete to $V_0$ and could be
added to $\mathcal Q$, a contradiction. Thus every fiber has chromatic number at most
$(hk+r-1)\tau$. The first three layers lie in $B_G(z,2)$.
Apply \cref{lem:layered-projection} with $X=C$, $K=G_E$, $q=1$,
and $b=(hk+r-1)\tau$.
\end{proof}

\begin{proof}[Proof of \cref{thm:common-two-control}(\ref{item:common-long-arms})]
Fix $s,r\ge3$, let $h=sr+1$ and $R=\max\{r+1,6\}$, and write
$\cC_k=\cC_k(K_{1,s}^{(r)})$. By \cref{lem:control-start}, $\cC_0$ is
$R$-controlled, and every $\cC_k$ is controlled at some finite radius. We first prove $R$-control of $\cC_k$ by induction on
$k$. For $k\ge1$, let $\varphi$ witness $R$-control of $\cC_{k-1}$.
Apply \cref{lem:configuration-descent} with depth $r$ and the bounds of
\cref{lem:spider-target-bounds}:
\[
\begin{aligned}
 m&=hk+(3s-5)\tau+1,&\theta&=sk+2s-1,\\
 c_{\mathrm{free}}&=2hk\tau\varphi(\tau),&c_{\mathrm{rich}}&=2(hk+r-1)\tau\varphi(\tau).
\end{aligned}
\]
This proves $R$-control.
For fixed $k$, \cref{cor:spider-biclique} shows that $\cC_k$ excludes
$K_{m_0,m_0}^{(q)}$ for $2\le q\le2r$, where $m_0=(sr+2)(k+1)$. Since
$R\le2r$, \cref{lem:exact-path-reduction} at radii $R,R-1,\ldots,4$
gives $3$-control. Since $6\le2r$, \cref{lem:css-two-control} at
radius three gives $2$-control.
\end{proof}

\subsection{Three arms of length two}
\label{subsec:arms-two}

The proof of \cref{thm:common-two-control}(\ref{item:common-short-claw})
lowers the control radius in four steps: from the finite radius of
\cref{lem:control-start} to $5$ by \cref{lem:configuration-descent}
with \cref{lem:s222-free-support,lem:s222-rich-support}, then to $4$
by \cref{lem:multicover-descent} with
\cref{lem:multicover-configurations}\ref{item:mc-free},
\cref{lem:s222-free-support}, and \cref{lem:five-cover-rich-impression},
then to $3$ by \cref{lem:exact-path-reduction}, and finally to $2$ by
\cref{lem:s222-three-multicover-bound,lem:multicover-control-reduction}.

Let $h=7$ be the number of vertices of $K_{1,3}^{(2)}$. All
configurations in this subsection have depth two, and all supports
consist of induced copies of $K_{1,3}^{(2)}$. By \cref{lem:long-holes},
there is a binding function $g_6$ for graphs with no induced cycle of
length at least six. Let
\[
 M_k(\tau)=\tau g_6(\tau)+hk\tau.
\]

\begin{lemma}\label{lem:s222-protected-tails}
Let $k\ge0$ and $\tau\ge1$, and let $(A,B)$ be a cut of $G$ with
$\pi^{\times}_{K_{1,3}^{(2)}}(G;A,B)\leq k$, and suppose that $\chi^4(G)\leq\tau$.
Let $Y\subseteq B$, and let $Q\subseteq V(G)\setminus Y$ cover $Y$.
Suppose that every $q\in Q$ has one of the following certificates:
\begin{enumerate}
\item an induced path $q-a-b$, where $\{a,b\}$ is anticomplete to $Y$
and at least one of $q,a$ belongs to $A$; or
\item two induced paths $q-a-b$ and $q-d-d'$, where the sets
$\{a,b\}$ and $\{d,d'\}$ are anticomplete to one another and to $Y$,
and $b\in A$.
\end{enumerate}
Then $\chi(G[Y])\leq M_k(\tau)$.
\end{lemma}

\begin{proof}
We verify the hypothesis of \cref{lem:local-crossing-packing} with
$R=4$ and $c=\tau g_6(\tau)$.
Let $Y'\subseteq Y$ have chromatic number greater than $c$.
Order $Q=(q_1,\ldots,q_N)$ and assign each vertex of $Y'$ to its first
neighbor, obtaining $Y_i'\subseteq N(q_i)$.
Color each $Y_i'$ with $[\tau]$. Some union $X$ of equal color
classes has $\chi(G[X])>g_6(\tau)$, and therefore contains an
induced cycle $z_0-\cdots-z_{\ell-1}-z_0$ with $\ell\ge6$.
Choose the first parent $q=q_i$ seeing this cycle. Then
$N(q)\cap\{z_0,\ldots,z_{\ell-1}\}\subseteq Y_i'\cap X$ is independent.

If $q$ has at least two neighbors on the cycle, orient and index them as $z_0,z_j$
with $2\le j\le\lfloor\ell/2\rfloor$. Let
$F=\{z_0,z_{\ell-1},z_2,z_3\}$ if $j=2$, and
$F=\{z_0,z_1,z_j,z_{j+1}\}$ if $j\ge3$.
Then $G[F]\cong2P_2$, and $q$ sees precisely one endpoint of each
edge. Together with the path $q-a-b$, whose vertices $a,b$ miss $Y$,
these give a
crossing $K_{1,3}^{(2)}$ centered at $q$.
If $q$ has exactly one neighbor $z_0$ on the cycle, distinguish the two
certificates.
\begin{enumerate}
 \item For the first certificate, take the copy centered at $z_0$
 with arms
 \[
  z_0-z_{\ell-1}-z_{\ell-2},\qquad z_0-z_1-z_2,\qquad z_0-q-a.
 \]
 \item For the second certificate, take the copy centered at $q$
 with arms
 \[
  q-z_0-z_1,\qquad q-a-b,\qquad q-d-d'.
 \]
\end{enumerate}
The independent-neighbor condition, the induced cycle, and the
certificate exclude all extra edges. Each copy meets $A$, contains
$z_0\in B$, and lies in $B(z_0,3)$.
Thus \cref{lem:local-crossing-packing} gives
$\chi(G[Y])\le\tau g_6(\tau)+hk\tau=M_k(\tau)$.
\end{proof}

\begin{lemma}\label{lem:s222-free-support}
Let $k,\tau\ge1$ and $\rho\ge4$, and let $\varphi$ witness
$\rho$-control of $\cC_{k-1}(K_{1,3}^{(2)})$. Suppose that
$G\in\cC_k(K_{1,3}^{(2)})$ satisfies $\chi^\rho(G)\le\tau$ and
contains a free configuration of depth two with $m\ge4\tau+1$
parent sets and target $C$. Then
$\chi(G[C])\leq4M_k(\tau)\varphi(\tau)$.
\end{lemma}

\begin{proof}
Assume $C\ne\varnothing$, and fix a node-width-$k$ decomposition
$(T,\lambda)$.
Take distance layers $(\Lambda_i)$ in a component of $K$ rooted at
$x\in C$. For $j\ge2$, let
\[
\begin{gathered}
 W=\Lambda_j,\quad U=\bigcup_{i\le j-2}\Lambda_i,\quad
 D(U)=\bigcup_{p\in P(U)}(D_{\iota(p)}\cap N_G(p)).
\end{gathered}
\]
The layer gap gives $(U\cup P(U))\perp W$.
Applying \cref{lem:layered-supports}\ref{item:layered-free} with the
supports of \cref{lem:free-short-support} for $s=3$, $r=2$, and
$I=\varnothing$ gives a support $\mathcal R$ such that
$S=S_T(\mathcal R)$ is reserved against $W$. The local supports of
\cref{lem:free-short-support}\ref{item:free-local-support} also contain
the entries adjacent to the parents of $U$, so
$\lambda(U\cup P(U)\cup D(U))\subseteq V(S)$.

For an attachment cut $(A,B)$, let $Y=W\cap B$. Then
$U\cup P(U)\cup D(U)\subseteq A$.
For $z\in\Lambda_{j-1}$, choose a predecessor $u_z\in U\cap N_K(z)$.
If $u_zz\in E(G)$, choose any parent $p_z$ of $u_z$; otherwise
choose a common parent of $u_z,z$. If $p_zz\in E(G)$, choose an
entry $d_z$ of $p_z$ and use $z-p_z-d_z$; otherwise use $z-u_z-p_z$.
These are first-type certificates for $q=z$ whose vertices other
than $z$ miss $W$.
Hence \cref{lem:s222-protected-tails} gives
$\chi(G[Y_1])\le M_k(\tau)$ for $Y_1=Y\cap N_G(\Lambda_{j-1})$.

Let $Y_0=Y\setminus Y_1$. For each $y\in Y_0$, choose a
$K$-predecessor $z\in\Lambda_{j-1}$ and a common parent $q$ of $y,z$.
Let $u=u_z$ and $p=p_z$, and choose an entry $d$ of $q$. The core
vertices $d,z_2(d)$ miss $Y_0\subseteq C$.
If $q\in A$ or some entry of $q$ lies in $A$, it is a first-type
certificate, choosing that entry when needed.
Otherwise $q$ and all its entries lie in $B$.
The layer gap gives $qu\notin E(G)$, so $q\ne p$.
The remaining possibilities give the following certificates.
\begin{enumerate}
 \item If $uz\in E(G)$, use the second-type certificate consisting
 of $q-z-u$ and $q-d-z_2(d)$.
 \item If $uz\notin E(G)$ and $pq\in E(G)$, use the first-type
 certificate $q-p-u$.
 \item If $uz,pq\notin E(G)$, use the second-type certificate
 consisting of $q-z-p$ and $q-d-z_2(d)$.
\end{enumerate}
Here $u,p\in A$ and both miss $Y_0$, while $z\perp Y_0$.
For the last pair, $pz_2(d)\notin E(G)$, since $P_i\perp Z_i^-$ and
$P_i\perp Z_j$ for $i\ne j$; also $pd\notin E(G)$: for different
indices since $P_i\perp Z_j$, and equal indices would imply
$d\in D(U)\subseteq A$, a contradiction.
All other required nonedges follow from the configuration.
Thus \cref{lem:s222-protected-tails} also gives
$\chi(G[Y_0])\le M_k(\tau)$.

Every fiber has chromatic number at most $2M_k(\tau)$.
The first two $K$-layers lie in $B_G(x,2)$.
Apply \cref{lem:layered-projection} with $X=C$, the graph $K$ of the
configuration, $q=1$, and $b=2M_k(\tau)$.
\end{proof}

\begin{lemma}\label{lem:s222-rich-support}
Let $k,\tau\ge1$, and let $\varphi$ witness $5$-control of
$\cC_{k-1}(K_{1,3}^{(2)})$. Suppose that
$G\in\cC_k(K_{1,3}^{(2)})$ satisfies $\chi^5(G)\le\tau$ and
contains a rich configuration of depth two with target $C$ and
independent parent set $E$, where every parent accesses at least
$3k+5$ cores. Then $\chi(G[C])\leq2M_k(\tau)\varphi(\tau)$.
\end{lemma}

\begin{proof}
Fix a node-width-$k$ decomposition $(T,\lambda)$.
Use \cref{lem:rich-short-support} with $s=3$, $r=2$, avoiding any
at most $3k$ core indices. For paths of length three or four, always
use the two endpoint copies specified there.

Take distance layers $(\Lambda_i)$ in a component of $G_E$
rooted at $z\in E$. For $j\ge3$, let
\[
 W=C\cap\Lambda_j,\qquad U=\bigcup_{i\le j-2}\Lambda_i,\qquad E'=E\cap U.
\]
By \cref{lem:rich-short-support}, the configuration is
$K_{1,3}^{(2)}$-supported, so, by
\cref{lem:layered-supports}\ref{item:layered-rich}, the family
$\mathcal R$ of all induced copies of $K_{1,3}^{(2)}$ that meet $E'$ and
have all noncore vertices in $U$ is a support. Then $S=S_T(\mathcal R)$
contains $\lambda(E')$ and is reserved against $W$.

For an attachment cut $(A,B)$, let $Y=W\cap B$,
$Y_1=Y\cap N_G(E\cap A)$, and $Y_0=Y\setminus Y_1$.
Choose a maximal crossing packing $\mathcal Q$ among copies centered
in $E$ and contained in $E\cup C\cup\bigcup_iZ_i$. Let
\[
 V_0=V(\mathcal Q),\quad I=\{i\mid Z_i\cap V_0\ne\varnothing\},
 \quad Y_i'=Y_i\setminus B_G(V_0,5)\quad(i=0,1).
\]
Then $|V_0|\le hk$. Each arm of a copy centered in $E$ meets at most
one core, since distinct cores are anticomplete; hence $|I|\le3k$, and
every parent retains five accesses.
Inward vertices in cores outside $I$ miss $V_0$, since they miss
$E\cup C$ and all old cores.

If $\chi(G[Y_1'])>\tau$, a component $F$ of $G[Y_1']$ has
$\chi(F)>\tau$. Choose $e\in E\cap A$ seeing $F$.
By \cref{lem:rooted-path-extension-local}, there is an induced
path $e-y-y'$ with $y,y'\in V(F)$. Add two core arms at $e$
in distinct cores outside $I$. This crossing copy has all vertices
outside $Z^-$ within distance two of $y$, and is therefore anticomplete
to $V_0$, a contradiction. Hence $\chi(G[Y_1'])\le\tau$.

Suppose $y\in Y_0'$. Let $a,b,c,y$ be the last four vertices of
a shortest $z$--$y$ path in $G_E$. All parents of $y$ lie in $B$.
We construct a support joining $E'\subseteq A$ to a parent
$v$ of $y$, avoiding $I$, such that every member either belongs to
$\mathcal R$ or has all vertices outside $Z^-$ in $B_G(y,4)$.

If $\{a,b\}\cap E\ne\varnothing$, take its last vertex as $u$,
follow the path to $y$, and append $v$. Delete repetitions, split at
internal parents, and shorten each segment. The resulting paths
have lengths two, three, or four; their centers are within distance
three of $y$, and their entries within distance four.

If $a,b\in C$ and $c\in E$, take $v=c$ and a parent $u$ of $a$.
Then $u\in U$ and $\dist_G(u,y)\le4$. When this distance is at most
three, the shortened path $u-a-b-v$ has the required support.
Otherwise it is induced: $ub$ is excluded by distance, $av$ by the
layer gap, and $uv$ by independence. Its endpoint copy at $u$ has
noncore vertices $u,a,b\in U$, so belongs to $\mathcal R$; its
companion has all vertices outside $Z^-$ within distance three of $y$.

Finally, let $a,b,c\in C$. Choose parents $u,w,v$ of $a,b,y$,
respectively, and use the walks $u-a-b-w$ and $w-b-c-y-v$,
omitting coincident-end walks. Here $u\in U$, $\dist_G(w,y)\le3$,
and $\dist_G(v,y)=1$. If $\dist_G(u,y)\le3$, shortening gives the required
supports. Otherwise $\dist_G(u,y)=4$, so $ub\notin E(G)$ and $u\ne w$.
If $wa\in E(G)$, then $w\in U$ and the support of $u-a-w$ lies
in $\mathcal R$. Otherwise $u-a-b-w$ is induced: its endpoint
copy at $u$ lies in $\mathcal R$, and its companion has all
vertices outside $Z^-$ in $B_G(y,4)$. The second walk has the same
distance property after shortening.

The support joins opposite sides, so some member $Q$ crosses.
No member of $\mathcal R$ crosses, since all its leaves lie in $S$.
Thus $V(Q)\setminus Z^-\subseteq B_G(y,4)$, and these vertices
miss $V_0$ because $\dist_G(y,V_0)>5$; its inward vertices were checked
above. This contradicts maximality, proving $Y_0'=\varnothing$.
Consequently
\[
 \chi(G[Y])\le(hk+1)\tau\le M_k(\tau).
\]
The first three layers lie in $B_G(z,2)$.
Apply \cref{lem:layered-projection} with $X=C$, $K=G_E$, $q=1$,
and $b=M_k(\tau)$.
\end{proof}

For the reduction from radius five to radius four, the next estimate
verifies condition~\ref{item:clean-rich} of \cref{lem:multicover-descent}
directly from the cover levels.

\begin{lemma}\label{lem:five-cover-rich-impression}
Let $\mu,m,\tau\ge1$, and let $N=\mu m$. Let $G$ satisfy
$\chi^4(G)\le\tau$ and have no impression of $K_{\mu,\mu}$ of order at
most four. For an ordered $5$-multicover of length $N$ in $G$, define
$C'$, $Z_i$, $B_i$, and $E_i$ as in \cref{obs:multicover-levels},
with $\rho=5$ and $\theta=\mu$. Then
$\chi(G[C'\cap N(E_i)])\le2^N\mu\tau$ for every $i\in[N]$.
\end{lemma}

\begin{proof}
By \cref{obs:multicover-levels}, $B_i\perp\bigcup_{a=0}^2L_{j,a}$
whenever $i<j$.
Fix $i$. For each $e\in E_i$, choose a $\mu$-set
$I(e)\subseteq\{j>i\mid N(e)\cap Z_j\ne\varnothing\}$; let $E_I=\{e\mid I(e)=I\}$ and $C_I=C'\cap N(E_I)$.
There are at most $2^N$ choices of $I$.
If $\chi(G[C_I])>\mu\tau$, choose $y_\ell\in C_I$ and parents
$e_\ell\in E_I\cap N(y_\ell)$ for $\ell\in[\mu]$ successively, with
\[
 y_\ell\notin\bigcup_{a<\ell}B(e_a,3).
\]
The deleted set costs at most $(\ell-1)\tau$ colors, so the choices
exist, and distinct $e_\ell,e_{\ell'}$ have distance at least three.
For each $j\in I$, choose $d_{\ell,j}\in L_{j,3}\cap N(e_\ell)$
and successive predecessors $z_{\ell,j,2}\in L_{j,2}$,
$z_{\ell,j,1}\in L_{j,1}$. Then
\[
 e_\ell-d_{\ell,j}-z_{\ell,j,2}-z_{\ell,j,1}-x_j
\]
is induced, since $E_I$ is anticomplete to levels zero, one, and two
of later covers. Map the two parts of $K_{\mu,\mu}$ to $(e_\ell)$ and
$(x_j)$. For nonincident edges, their core vertices lie in distinct
anticomplete cores. The parents miss the first three levels of the
other path, and an edge $e_\ell d_{\ell',j}$ with $\ell'\ne\ell$ would give
distance at most two between $e_\ell$ and $e_{\ell'}$. Hence these paths form a forbidden impression
of order four. Thus $\chi(G[C_I])\le\mu\tau$, and summing over $I$
proves the bound.
\end{proof}

\begin{lemma}\label{lem:s222-three-multicover-bound}
Let $k\ge0$ and $\tau\ge1$, and let $G$ satisfy
$\sigma_{K_{1,3}^{(2)}}(G)\le k$ and $\chi^2(G)\le\tau$.
Let $N=3\tau+1$ and $N_0=(k+1)N$. If $C$ admits an
ordered $3$-multicover of length $3N_0$, then
$\chi(G[C])\leq2((N_0+hk)\tau+2N\tau^3)$.
\end{lemma}

\begin{proof}
We first prove two claims, which use only $\chi^2(G)\le\tau$.

\begin{claim}\label{clm:s222-rooted-three-cover}
Let $(\{x\},L_1,L_2)$ be a $3$-cover of a set $X$ with
$\chi(G[X])>2\tau^3$. Then $G[\{x\}\cup L_1\cup L_2]$ contains an induced
copy of $K_{1,3}^{(2)}$ centered at $x$.
\end{claim}
\begin{claimproof}
Color $L_1$ and $L_2$ with $\tau$ colors, obtaining independent parts
$A_i\subseteq L_1$ and $B_j\subseteq L_2$. For each $i,j$, let
$D_{ij}\subseteq A_i$ be an inclusion-minimal set covering
$B_j\cap N(A_i)$.
Each $a\in D_{ij}$ has a private neighbor $y_a$ in that set.
If $|D_{ij}|\ge3$, the paths $x-a-y_a$ for three such vertices
induce $K_{1,3}^{(2)}$: both selected layers are independent, the private
neighbors exclude cross-edges, and $x\perp L_2$.
Otherwise $D=\bigcup_{i,j}D_{ij}$ covers $L_2$ and has size at most
$2\tau^2$. Hence $X\subseteq B(D,2)$ and
$\chi(G[X])\le2\tau^3$, a contradiction.
\end{claimproof}

\begin{claim}\label{clm:s222-three-cover-cut}
Let $(A,B)$ be a cut of $G$, and let $X\subseteq B$ admit an
ordered $3$-multicover $(\{x_i\},L_{i,1},L_{i,2})_{i=1}^N$ with
$x_i\in A$ for every $i\in[N]$. Suppose that every vertex of
$\bigcup_iL_{i,2}$ with a
neighbor in $X$ misses $\{x_1,\ldots,x_N\}$. If $\chi(G[X])>2N\tau^3$,
then there is an induced copy of $K_{1,3}^{(2)}$ meeting both $A$ and
$B$, all of whose vertices belong to the covers or to $X$, and each of
whose vertices in a base has a neighbor in $X$.
\end{claim}
\begin{claimproof}
Restrict $L_{i,2}$ to $L_{i,2}\cap N(X)$, and let
$X_i=X\cap N(L_{i,2}\cap B)$.
If $\chi(G[X_i])>2\tau^3$, apply \cref{clm:s222-rooted-three-cover}
to $(\{x_i\},L_{i,1},L_{i,2}\cap B)$, obtaining a copy centered in $A$
with arm endpoints in $B$.
Otherwise choose $c\in X\setminus\bigcup_iX_i$ and parents
$v_i\in L_{i,2}\cap N(c)$. All $v_i\in A$. Since
$N=3\tau+1$ and $\chi(G[N(c)])\le\tau$, four are independent.
In increasing cover order denote them by $p_1,\ldots,p_4$, their
apices by $y_1,\ldots,y_4$, and choose neighbors $a_j$ of $p_j$ in
level one of the same cover. The apices miss all selected parents; different cores are
anticomplete and miss $c$; and $p_j$ misses $a_i$ for $i<j$.

If $p_i a_j\in E(G)$ for some $1\le i<j\le3$, use
\[
 p_i-a_i-y_i,\qquad p_i-a_j-y_j,\qquad p_i-c-p_4.
\]
The parents are independent, and $p_4$ misses $a_i$ and $a_j$,
so these paths induce $K_{1,3}^{(2)}$. Otherwise the paths
$c-p_j-a_j$, $j\in[3]$, induce $K_{1,3}^{(2)}$.
Both copies cross through $c\in B$ and their parents in $A$.
All selected parents see $c$; in the first case of the proof, they
belong to the restricted bases.
\end{claimproof}

Now let $f=(N_0+hk)\tau+2N\tau^3$, and suppose $\chi(G[C])>2f$.
In a node-width-$k$ decomposition, \cref{lem:balanced-core} gives a
displayed cut with at least $N_0$ of the $3N_0$ apices on each side.
Name its sides $A,B$ so that $C^*=C\cap B$ has chromatic number
greater than $f$. Retain $N_0$ covers with apices in $A$, and let
\[
 C^0=C^*\setminus\bigcup_{i=1}^{N_0}B(x_i,2),\qquad
 \chi(G[C^0])>hk\tau+2N\tau^3.
\]
Every parent seeing $C^0$ misses all retained apices, since otherwise
its target neighbor would be in one of the deleted balls.
Partition the retained covers into $k+1$ consecutive blocks of size
$N$ and process them in reverse order.
After $q\le k$ choices with union $V_q$, let
$C^q=C^0\setminus B(V_q,2)$. Then
\[
 \chi(G[C^q])>h(k-q)\tau+2N\tau^3\ge2N\tau^3.
\]
The covers of the current block form an ordered $3$-multicover of
$C^q$. By \cref{clm:s222-three-cover-cut}, they give a crossing copy,
each vertex of which in a base has a neighbor in $C^q$.
New target vertices have distance greater than two from $V_q$, and
new parents therefore have distance greater than one. The remaining
vertices belong to earlier cores, which miss all old cover vertices
by the cover order and all targets by the core condition.
Thus the new copy is anticomplete to the previous copies.
The resulting $k+1$ crossing copies contradict the width bound.
\end{proof}

\begin{proof}[Proof of \cref{thm:common-two-control}(\ref{item:common-short-claw})]
Write $\cC_k=\cC_k(K_{1,3}^{(2)})$. By \cref{lem:control-start},
$\cC_0$ is $R$-controlled for every $R\ge1$, and every $\cC_k$ is
controlled at some finite radius. We prove $4$-control
by induction on $k$; $5$-control is an intermediate conclusion in each
induction step. Let $k\ge1$, and let $\varphi$ witness $4$-control of
$\cC_{k-1}$. The same function witnesses $5$-control, because
$\chi^4\le\chi^5$.

Apply \cref{lem:configuration-descent} with depth two, $R=5$, and
\[
 m=4\tau+1,\quad \theta=3k+5,\quad
 c_{\mathrm{free}}=4M_k(\tau)\varphi(\tau),\quad
 c_{\mathrm{rich}}=2M_k(\tau)\varphi(\tau),
\]
using \cref{lem:s222-free-support,lem:s222-rich-support}.
Thus $\cC_k$ is $5$-controlled.

By \cref{cor:spider-biclique} with $s=3$ and $r=2$, $\cC_k$ excludes
$K_{m_0,m_0}^{(q)}$ for $q=2,3,4$, where $m_0=8(k+1)$. Now fix
$\tau\ge1$.
For $G\in\cC_k$ with
$\chi^4(G)\le\tau$, \cref{lem:css-impressions} gives
$\mu=\mu(k,\tau)$ such that $G$ has no impression of
$K_{\mu,\mu}$ of order at most four. We apply
\cref{lem:multicover-descent} with $\rho=5$, $\theta=\mu$, and
$m=4\tau+1$; let $N=\mu m$. By
\cref{lem:multicover-configurations}\ref{item:mc-free} with $r=2$ and
\cref{lem:s222-free-support} at radius four,
condition~\ref{item:clean-free} of \cref{lem:multicover-descent} holds
with $c_{\mathrm{free}}=4M_k(\tau)\varphi(\tau)$. By
\cref{lem:five-cover-rich-impression}, condition~\ref{item:clean-rich}
holds with $c_{\mathrm{rich}}=2^N\mu\tau$. Since $\cC_k$ is $5$-controlled, it is
$4$-controlled. This completes the induction on $k$.

For each fixed $k$, the same three excluded subdivisions and
\cref{lem:exact-path-reduction} give $3$-control. Finally, for
$\tau\ge1$, \cref{lem:s222-three-multicover-bound} gives $N_0$ and a
bound, both depending only on $k$ and $\tau$, such that no $G\in\cC_k$
with $\chi^2(G)\le\tau$ has an ordered $3$-multicover of length $3N_0$
whose target has chromatic number greater than this bound.
\Cref{lem:multicover-control-reduction} gives $2$-control.
\end{proof}

\subsection{A handle on three vertices}
\label{subsec:broom-three}

Fix $d\ge2$, and let $H=B_{3,d}$ and $h=d+3$. All configurations in
this subsection have depth two, and all supports consist of induced
copies of $H$. Recall that, for a parent $p$ adjacent to an entry $e$, the
core path $p-e-z_2(e)$ is induced. Let
\[
 L_{3,d}(\tau)=\tau\bigl((d-1)\tau+1\bigr)\qquad(\tau\ge1).
\]

\begin{lemma}\label{lem:broom-critical-root}
Let $\tau\ge1$, and suppose that $\chi^2(G)\le\tau$.
Let $G[C]$ be vertex-critical and
$a\notin C$, with $\dist(a,C)\le2$ and
$\chi(G[C])>L_{3,d}(\tau)$. Then $G$ contains an induced copy of $B_{3,d}$
rooted at $a$.
\end{lemma}
\begin{proof}
Let $q=(d-1)\tau+1$, $X=C\cap B(a,2)$, and $\ell=\chi(G[C])$.
Then
\[
 X\ne\varnothing,\qquad \chi(G[X])\le\tau,\qquad
 \chi(G[C\setminus X])\le\ell-1,\qquad \ell-1\ge\tau q.
\]
Suppose $|N(x)\cap(C\setminus X)|<q$ for every $x\in X$.
Choose a proper map $f\colon C\setminus X\to[\ell-1]$, an independent
partition $X=X_1\mathbin{\dot\cup}\cdots\mathbin{\dot\cup}X_\tau$,
and disjoint $q$-sets $A_1,\ldots,A_\tau\subseteq[\ell-1]$.
For $x\in X_i$, the set
\[
 A_i\setminus f\bigl(N(x)\cap(C\setminus X)\bigr)
\]
is nonempty. Choosing $f(x)$ from it extends $f$ to an
$(\ell-1)$-coloring of $G[C]$, a contradiction.

Hence some $x\in X$ has at least $q$ neighbors in $C\setminus X$.
Their $\tau$-colorability gives an independent $d$-set $Y$ among them.
Since $Y\cap B(a,2)=\varnothing$, we have $\dist(a,x)=2$.
For a shortest path $a-z-x$, any edge from $\{a,z\}$ to $Y$
would place its endpoint in $B(a,2)$. Thus
\[
 \{a,z\}\perp Y,\qquad G[\{a,z,x\}\cup Y]\cong B_{3,d},
\]
with root $a$.
\end{proof}

For $k\ge0$ and $\tau\ge1$, let
\[
 b_{3,d}(k,\tau)=\bigl(L_{3,d}(\tau)+hk\tau\bigr)
 \bigl(1+(k+1)L_{3,d}(\tau)+hk\tau\bigr).
\]
\begin{lemma}\label{lem:broom-second-boundary}
Let $k\ge0$ and $\tau\ge1$.
Let $G$ satisfy $\chi^2(G)\le\tau$, and let $(A,B)$ be a cut of
$G$ with $\pi_H^\times(G;A,B)\le k$. Then
$\chi(G[\partial_A^2]),\chi(G[\partial_B^2])\le b_{3,d}(k,\tau)$.
\end{lemma}

\begin{proof}
Let $L=L_{3,d}(\tau)$ and $a=L+hk\tau$. For the first boundaries, apply
\cref{lem:local-crossing-packing} to $Y=\partial_A^1$ with $R=2$ and
$c=L$. If $Y'\subseteq Y$ satisfies $\chi(G[Y'])>L$, choose a
vertex-critical induced subgraph $G[X]$ of $G[Y']$ with
$\chi(G[X])=\chi(G[Y'])$, and an edge $ve$ with $v\in X$ and $e\in B$.
\Cref{lem:broom-critical-root}, applied in $G[X\cup\{e\}]$, gives an
induced copy of $H$ rooted at $e$. It crosses the cut and lies in
$B(Y',1)$. Hence $\chi(G[\partial_A^1])\le L+hk\tau=a$, and
symmetry gives the other first-boundary bound.

For the second boundaries, we verify the hypothesis of
\cref{lem:second-boundary}, in its notation. Let $D,C,E$ be sets as
there, let $j\ge0$, and let
$C'\subseteq C\cap N(\bigcup_{i>j}D_i)$ satisfy $\chi(G[C'])>L$.
Choose a vertex-critical induced subgraph $G[X]$ of $G[C']$ with
$\chi(G[X])=\chi(G[C'])$, and let $i$ be the largest index with
$D_i\cap N(X)\ne\varnothing$; then $i>j$. In
$\Gamma=G[\{e_i\}\cup(D_i\cap N(X))\cup X]$, the vertex $e_i$ has
distance two from $X$, because $D_i\subseteq N(e_i)$ and $E\perp C$.
By \cref{lem:broom-critical-root}, $\Gamma$ contains an induced copy of $H$
rooted at $e_i$, so $i$ is active for $C'$. Thus
\cref{lem:second-boundary} gives
$\chi(G[\partial_A^2])\le a\bigl(1+(k+1)L+hk\tau\bigr)=b_{3,d}(k,\tau)$,
and symmetry gives the bound on $\partial_B^2$.
\end{proof}

\begin{lemma}\label{lem:broom-three-free-support}
Let $\tau\ge1$. Suppose that $\chi^1(G)\le\tau$ and
that $G$ contains a free configuration of depth two with target $C$
and $m\ge(3d-2)\tau+1$ parent sets. Then the following hold.
\begin{enumerate}[label=(\arabic*)]
\item For each $x\in C$, a support contains $x$ and all of $P(x)$;
its noncore vertices lie in $\{x\}\cup P(x)$.
\item For each oriented edge $xz$ of $K$, a support joins $x,z$
and has all its noncore vertices in $\{x,z\}\cup P(x)$.
\end{enumerate}
\end{lemma}

\begin{proof}
After excluding one index, representatives in $N_G(x)$ give
$3d-2$ independent parents with distinct indices; note that
$3d-2\ge2d$, since $d\ge2$.
For a parent $p$, choose an entry $e_p$ in its own core.
Write $J(c,u,v;Y)=G[\{c,u,v\}\cup Y]$ for the broom with center
$c$, handle $c-u-v$, and independent leaf set $Y$ of size $d$.

For (1), fix $p\in P(x)$.
Choose an independent $2d$-set $Q\subseteq P(x)$ with distinct
indices, avoiding $\iota(p)$. If $|Q\setminus N_G(p)|\ge d$, choose a
$d$-set $Y$ there and use $J(x,p,e_p;Y)$.
Otherwise $|Q\cap N_G(p)|\ge d+1$; choose $q_0\in Q\cap N_G(p)$ and a
$d$-set $Y\subseteq(Q\cap N_G(p))\setminus\{q_0\}$, and use
\[
 \{J(p,e_p,z_2(e_p);Y),\ J(x,q_0,e_{q_0};Y)\}.
\]
These copies intersect in $Y$. Their union over $p\in P(x)$ is a
support $\mathcal R_x$ containing $x$ and $P(x)$.

For (2), it suffices to find one copy containing $z$ and meeting
$\mathcal R_x$. If $xz\notin E(G)$, fix a common parent $p$.
Choose an independent $(3d-2)$-set $Q\subseteq P(x)$ with distinct
indices and, if $xz\notin E(G)$, avoiding $\iota(p)$, and let
$Q_1=Q\cap N_G(z)$, $Q_0=Q\setminus N_G(z)$.
If $|Q_1|\ge d+1$, choose $q_0\in Q_1$ and a $d$-set
$Y\subseteq Q_1\setminus\{q_0\}$, and use $J(z,q_0,e_{q_0};Y)$.
Otherwise $|Q_0|\ge2d-2\ge d$, and one of the following cases applies.
\begin{enumerate}[label=(\alph*)]
 \item If $xz\in E(G)$, choose $q_0\in Q_0$ and
 $Y\in\binom{Q_0\setminus\{q_0\}}{d-1}$. Use
 \[
  J(x,q_0,e_{q_0};\{z\}\cup Y).
 \]
 \item If $xz\notin E(G)$ and $|Q_0\setminus N_G(p)|\ge d$,
 choose $Y\in\binom{Q_0\setminus N_G(p)}d$ and use $J(x,p,z;Y)$.
 \item If $xz\notin E(G)$ and $|Q_0\setminus N_G(p)|<d$, then
 $|Q_0\cap N_G(p)|\ge d-1$. Choose
 $Y\in\binom{Q_0\cap N_G(p)}{d-1}$ and use
 \[
  J(p,e_p,z_2(e_p);\{z\}\cup Y).
 \]
\end{enumerate}
The displayed choices exclude all noncore chords; distinct indices
and the anticompleteness conditions in the definition of a free
configuration exclude all core chords.
Each copy meets $\mathcal R_x$ at $x$, $p$, or a selected parent,
and all its noncore vertices lie in $\{x,z\}\cup P(x)$.
Adjoin $\mathcal R_x$ to obtain the required support.
\end{proof}

In a rich configuration, entries $e_0,\ldots,e_d$ at $p$ in $d+1$
distinct accessed cores give a pure copy of $H$ at $p$, with handle
$p-e_0-z_2(e_0)$ and leaf set $\{e_1,\ldots,e_d\}$.

\begin{lemma}\label{lem:broom-three-rich-support}
Consider a rich configuration of depth two with target $C$ and
independent parent set $E$, each vertex of which accesses at least $2d$
cores. Then every induced path of length two or three with endpoints in
$E$ and interior in $C$ has a $B_{3,d}$-support joining its endpoints,
each member of which is centered at an endpoint and has all noncore
vertices on the path. Moreover, every edge $px$ with $p\in E$ and
$x\in C$ belongs to an induced copy of $B_{3,d}$ centered at $p$ whose
only noncore vertices are $p,x$. In particular, the configuration is
$B_{3,d}$-supported.
\end{lemma}

\begin{proof}
Since $2d\ge d+1$, every $p\in E$ lies in a pure copy, as noted before
the lemma; hence the last assertion follows from the first.
For an edge $px$, choose entries $e_0,\ldots,e_{d-1}$ at $p$ in
distinct cores. The handle $p-e_0-z_2(e_0)$, with leaves
$x,e_1,\ldots,e_{d-1}$, gives the required copy.
For a path $p-x-y-q$, take $d$ entry leaves at each endpoint in
distinct cores. The two brooms with handles $p-x-y$ and $q-y-x$
intersect in $\{x,y\}$.

For $p-x-q$, choose a set $D$ of $2d$ entries at $p$ in distinct cores.
If $D\setminus N(q)$ contains a $d$-set $Y$, use the handle
$p-x-q$ with leaves $Y$.
Otherwise $|D\cap N(q)|\ge d+1$; choose $e_0\in D\cap N(q)$ and a
$d$-set $Y\subseteq(D\cap N(q))\setminus\{e_0\}$.
The copies with handles $p-e_0-z_2(e_0)$ and $q-e_0-z_2(e_0)$ and leaves
$Y$ intersect. All displayed copies are induced because distinct
cores are anticomplete, cores miss $C$, and inward vertices miss $E$.
\end{proof}

\begin{lemma}\label{lem:broom-three-targets}
Let $k,\tau\ge1$. Suppose that $\cC_{k-1}(B_{3,d})$
is $2$-controlled by $\varphi$.
Let $G\in\cC_k(B_{3,d})$ satisfy $\chi^2(G)\le\tau$, and let
$b=b_{3,d}(k,\tau)$. For configurations of depth two in $G$, the
following hold.
\begin{enumerate}[label=(\alph*)]
\item\label{item:broom-three-free} A free configuration with
$m\ge(3d-2)\tau+1$ parent sets has $\chi(G[C])\le4b\varphi(\tau)$.
\item\label{item:broom-three-rich} A rich configuration with an
independent parent set and at least $2d$ accesses per parent has
$\chi(G[C])\le2b\varphi(\tau)$.
\end{enumerate}
\end{lemma}

\begin{proof}
Fix a node-width-$k$ decomposition $(T,\lambda)$. By
\cref{lem:node-span-width,lem:broom-second-boundary}, both second boundaries of every
displayed cut have chromatic number at most $b$.

For~\ref{item:broom-three-free}, take distance layers $(\Lambda_i)$
in a component of $K$ rooted at $x\in C$. For $j\ge2$, let
\[
 W=\Lambda_j,\quad U=\bigcup_{i\le j-2}\Lambda_i.
\]
By \cref{lem:broom-three-free-support} and
\cref{lem:layered-supports}\ref{item:layered-free}, there is a support
$\mathcal R$ with $U\cup P(U)\subseteq V(\mathcal R)$, all noncore
vertices in $U\cup P(U)$, and $V(\mathcal R)\perp W$.
Each $z\in\Lambda_{j-1}$ has a predecessor $u_z\in U$, and either the edge
$u_zz$ or a common parent gives $\dist(z,U\cup P(U))\le1$.
Let $W_1=W\cap N_G(\Lambda_{j-1})$ and $W_0=W\setminus W_1$. Then
\[
 W_1\subseteq B_G(V(\mathcal R),2).
\]
Add the oriented supports for $u_zz$, for all $z\in\Lambda_{j-1}$,
obtaining $\mathcal R'$. Then $\Lambda_{j-1}\subseteq V(\mathcal R')$,
and all its noncore vertices lie in $U\cup P(U)\cup\Lambda_{j-1}$,
anticomplete to $W_0$; with the core conditions, this gives
$V(\mathcal R')\perp W_0$. Every vertex of $W_0$ shares a parent with
its $K$-predecessor in $\Lambda_{j-1}$, and hence
\[
 W_0\subseteq B_G(V(\mathcal R'),2).
\]
By \cref{obs:reserved}, $S_T(\mathcal R)$ and $S_T(\mathcal R')$
are reserved against $W_1$ and $W_0$, respectively. Each attachment fiber
lies in the second boundary of its cut and has chromatic number at
most $b$. The first two layers lie in a radius-two ball. Apply
\cref{lem:layered-projection} with $X=C$, the graph $K$ of the
configuration, and $q=2$ to prove~\ref{item:broom-three-free}.

For~\ref{item:broom-three-rich}, take distance layers $(\Lambda_i)$
in a component of $G_E$ rooted at $z\in E$.
For $j\ge2$, let $W=C\cap\Lambda_j$, $U=\bigcup_{i\le j-2}\Lambda_i$,
and $E'=E\cap U$. By \cref{lem:broom-three-rich-support}, the
configuration is $B_{3,d}$-supported, so, by
\cref{lem:layered-supports}\ref{item:layered-rich}, the family
$\mathcal R$ of all induced copies of $H$ that meet $E'$ and have all
noncore vertices in $U$ is a support with $E'\subseteq V(\mathcal R)$
and $V(\mathcal R)\perp W$.

We also have $C\cap U\subseteq V(\mathcal R)$. For $x\in C\cap U$,
let $v$ precede
$x$ on a shortest $z$--$x$ path. Take $p=v$ if $v\in E$, and
otherwise a parent $p$ of $v$; then $p\in E'$.
If $px\in E(G)$, use the edge assertion of
\cref{lem:broom-three-rich-support}. Otherwise the handle $p-v-x$
and $d$ entry leaves at $p$ give a broom in $\mathcal R$ containing
$x$.
Therefore
\[
 U\subseteq V(\mathcal R),\qquad V(\mathcal R)\perp W,
 \qquad W\subseteq B_{G_E}(U,2)\subseteq B_G(V(\mathcal R),2).
\]
By \cref{obs:reserved}, $S_T(\mathcal R)$ is reserved
against $W$. Each attachment fiber lies in the second boundary of its
cut and has chromatic number at most $b$. The first two layers have
chromatic number at most $\tau$.
Apply \cref{lem:layered-projection} with $X=C$, $K=G_E$, and $q=1$ to
obtain $\chi(G[C])\le2b\varphi(\tau)$.
\end{proof}

\begin{proof}[Proof of \cref{thm:common-two-control}(\ref{item:common-short-broom})]
Write $\cC_k=\cC_k(B_{3,d})$. By
\cref{lem:control-start}, $\cC_0$ is $2$-controlled, and every
$\cC_k$ is controlled at some finite radius. Induct on
$k$. For $k\ge1$, let $\varphi$ witness $2$-control of $\cC_{k-1}$.
Apply \cref{lem:configuration-descent} with depth two, $R=2$, and
\[
 m=(3d-2)\tau+1,\qquad \theta=2d,\qquad
 c_{\mathrm{free}}=4b_{3,d}(k,\tau)\varphi(\tau),\qquad
 c_{\mathrm{rich}}=2b_{3,d}(k,\tau)\varphi(\tau),
\]
using \cref{lem:broom-three-targets}.
\end{proof}

This completes the proof of \cref{thm:common-two-control}.

\section{Polynomial degree-boundedness}\label{sec:degreeboundedness}

In this section we prove~\cref{thm:forest-induced-biclique-chi},
which immediately implies~\cref{cor:forest-induced-biclique-polychi}.

\polydegree*

\subsection{Graphs with no \texorpdfstring{$C_4$}{C4} subgraph}
We first bound the degeneracy of graphs with $\branchpi_H\le k$ and
no $C_4$ as a subgraph.
Large minimum degree forces large treewidth, so some displayed cut is
crossed by every edge of a large matching. Since there is no $C_4$
subgraph, this matching contains a large induced matching, and the
minimum degree lets us extend its edges to an induced $H$-packing
whose members all cross the same cut.

\begin{lemma}\label{lem:biclique-c4-tree}
Let $H$ be a tree on $h\ge2$ vertices, let $k\in\Nzero$, and let
$
 D=\max\{2h(k+1),24k^2\}.
$
If $\branchpi_H(G)\le k$ and $G$ contains no $C_4$ as a subgraph,
then $\dgn(G)\le D-1$.
\end{lemma}

\begin{proof}
Suppose that a nonempty induced subgraph $F$ of $G$ has $\delta(F)\ge D$.
Then $\branchpi_H(F)\le k$ by \cref{lem:basic}\ref{item:heredity}.
Moreover, $F$ contains no $C_4$ as a
subgraph, so
\begin{equation}\label{eq:biclique-codegree}
 |N_F(x)\cap N_F(y)|\le1
 \qquad\text{for all distinct }x,y\in V(F).
\end{equation}
Indeed, two distinct common neighbors would form a cycle of length
four with $x$ and $y$.

We claim the following:

\begin{claim}\label{clm:growingmatching}
Any induced matching of size $k+1$ in $F$
extends to $k+1$ pairwise anticomplete induced copies of $H$, each
containing a different edge of the matching.
\end{claim}

\begin{claimproof}
 Fix an edge $uv$ of $H$ and an ordering $u,v,w_3,\dots,w_h$ of
$V(H)$ in which each $w_i$ has exactly one neighbor among its
predecessors; such an ordering exists because $H$ is a tree.
Identify a copy of $uv$ with each prescribed matching edge, and grow
the $k+1$ copies one vertex at a time according to this ordering.
If $S$ is the set of vertices already chosen and $x\in S$ is the
required neighbor of the next vertex, the forbidden choices in
$N_F(x)$ are contained in
\[
 S\cup\bigcup_{z\in S\setminus\{x\}}N_F(z).
\]
By \eqref{eq:biclique-codegree}, the number of forbidden choices is
at most
\[
 |S|+\sum_{z\in S\setminus\{x\}}
 |N_F(x)\cap N_F(z)|
 \le2|S|-1<2h(k+1)\le\delta(F).
\]
Here $|S|<h(k+1)$, since an extension is still required.
Thus there is an unused neighbor of $x$ with no other neighbor in
$S$. This preserves both inducedness of every partial copy and
anticompleteness between different copies, proving the claim.
\end{claimproof}

Since $\delta(F)\ge D\ge4$, the graph $F$ has an edge. Fix a branch
decomposition of $F$ whose displayed cuts $(A,B)$ all satisfy
$\pi^\times_H(F;A,B)\le k$. By
\cref{lem:mmw-tw} and $\delta(F)\le\dgn(F)\le\tw(F)$, we have
\[
 3\mmw(F)\ge\tw(F)+1\ge D+1>24k^2.
\]
Hence some displayed cut $(A,B)$ is crossed by every edge of
a matching $M$ with $m\ge8k^2+1$ edges.

Let $J=F[V(M)]$, and let $L$ be the graph with vertex set $M$ in
which two matching edges are adjacent when an edge of $F$ joins their
endpoint sets. Every edge of $L$ is witnessed by an edge of $J$
outside $M$, and each edge of $J$ outside $M$ witnesses at most one
edge of $L$; hence $|E(L)|\le|E(J)|-m$. As $J$ has $2m$ vertices and
no $C_4$ as a subgraph, Reiman's bound~\cite{Reiman} gives
$|E(J)|\le m\bigl(1+\sqrt{8m-3}\bigr)/2$, so
$\overline d(L)=2|E(L)|/m\le\sqrt{8m-3}-1$. By Tur\'an's
theorem~\cite{Turan}, applied to the complement of $L$,
\[
 \alpha(L)\ge\frac{m}{\overline d(L)+1}\ge\frac{m}{\sqrt{8m-3}}>k,
\]
where the last inequality holds because
$m^2-k^2(8m-3)=m(m-8k^2)+3k^2>0$.
Therefore $L$ has an independent set of size $k+1$, which is an
induced matching of $F$ whose edges all cross $(A,B)$. By
\cref{clm:growingmatching}, its edges extend to $k+1$ pairwise
anticomplete induced copies of $H$, each of which crosses $(A,B)$
because it contains a matching edge. Hence
$\pi^\times_H(F;A,B)\ge k+1$, a contradiction.
\end{proof}

To deal with the disconnected case, we use the following simple
observation.

\begin{lemma}\label{lem:biclique-disconnected-packing}
Let $H$ be a disconnected graph and let $k\in\Nzero$. If
$\branchpi_H(G)\le k$, then every induced $H$-packing in $G$ has at
most $3k$ members.
\end{lemma}

\begin{proof}
Let $Q_1,\ldots,Q_q$ be an induced $H$-packing in $G$. The
assertion is immediate if $q=0$. If $q=1$, then $G$ is not $H$-free,
so $k\ge1$ by \cref{lem:basic}\ref{item:zero}.
Assume $q\ge2$, fix a branch decomposition of $G$ whose displayed cuts
$(A,B)$ all satisfy $\pi^\times_H(G;A,B)\le k$, and choose one
representative from each $Q_i$. By \cref{lem:balanced-core}, some
displayed cut $(A,B)$ has at least $\lceil q/3\rceil$ representatives
on each side.
Let $m$ be the number of copies that cross the cut, and let $a$ and
$b$ be the numbers of copies contained in $A$ and in $B$,
respectively. Since a representative on $A$
belongs to one of the $a+m$ copies meeting $A$, and similarly for
$B$, we have
\[
 a+m\ge\lceil q/3\rceil,
 \qquad b+m\ge\lceil q/3\rceil.
\]

Write $H=H_1+H_2$, where $H_1,H_2$ are nonempty
unions of components. Pair $\min\{a,b\}$ copies contained in $A$
with distinct copies contained in $B$. From each pair take an
induced copy of $H_1$ in its member contained in $A$ and an induced
copy of $H_2$ in its member contained in $B$.
Their union induces $H$ and crosses the cut. Together with the
original $m$ crossing copies, these form an induced $H$-packing:
different new members use disjoint collections of the original
pairwise anticomplete copies. Its size is
\[
 m+\min\{a,b\}=\min\{a+m,b+m\}\ge\lceil q/3\rceil.
\]
Thus $k\ge\pi^\times_H(G;A,B)\ge\lceil q/3\rceil$, giving $q\le3k$.
\end{proof}

\begin{corollary}\label{cor:biclique-c4-forest}
Let $H$ be a forest on $h\ge2$ vertices, let $k\in\Nzero$, and let
$D_{h,k}=\max\bigl\{2\bigl((3k+1)h+1\bigr),24k^2\bigr\}$.
If $\branchpi_H(G)\le k$ and $G$ contains no $C_4$ as a subgraph,
then $\dgn(G)\le D_{h,k}-1$.
\end{corollary}

\begin{proof}
If $H$ is connected, apply \cref{lem:biclique-c4-tree}, noting that
$2h(k+1)\le2((3k+1)h+1)$.
Otherwise $H$ is disconnected (possibly edgeless), and by
\cref{lem:biclique-disconnected-packing}, $G$ has no
induced $H$-packing with $3k+1$ members. In other words, $G$
contains no induced copy of the disjoint union $H'$ of $3k+1$ copies
of $H$.
Form a tree $T$ from $H'$ by adding one new vertex adjacent to
exactly one chosen vertex in each component of $H'$.
The graph $T$ is connected and acyclic, has $(3k+1)h+1$ vertices,
and contains $H'$ as an induced subgraph. Hence $G$ is $T$-free,
so $\branchpi_T(G)=0$ by \cref{lem:basic}\ref{item:zero}.
Applying \cref{lem:biclique-c4-tree} with $T$ in place of $H$ and
$0$ in place of $k$, we obtain
\[
 \dgn(G)\le2\bigl((3k+1)h+1\bigr)-1
 \le D_{h,k}-1.
\]
\end{proof}

\subsection*{Excluding a biclique as a subgraph}

The following induced-subgraph extraction theorem is the only
additional structural result needed.

\begin{theorem}[Hunter, Milojevi\'c, Sudakov, and Tomon
{\cite[Theorem~1.6]{hunter2024k}}]
\label{thm:biclique-induced-extraction}
For every integer $d\ge1$, there is an integer $b(d)\ge2$ such
that, for every integer $s\ge b(d)$, a graph of average degree at
least $s^{1000d^2}$ with no $K_{s,s}$ as a subgraph has an induced
subgraph of average degree at least $d$ containing no $C_4$ as a
subgraph.
\end{theorem}

It is essential here that the extracted subgraph is induced:
\cref{lem:basic}\ref{item:heredity} then preserves the bound on
$\branchpi_H$.

\begin{proposition}\label{prop:biclique-polynomial-degree}
For all $h\ge2$ and $k\in\Nzero$, there is a positive integer $b$
such that, with $\ell=4000D_{h,k}^2$ for $D_{h,k}$ as in
\cref{cor:biclique-c4-forest}, the following holds for every forest
$H$ on $h$ vertices and every integer $s\ge2$: every graph $G$ with
$\branchpi_H(G)\le k$ and no $K_{s,s}$ as a subgraph satisfies
$\dgn(G)<(s+b)^\ell$.
\end{proposition}

\begin{proof}
Let $D=D_{h,k}$ and $d=2D$, and take $b=b(d)$ from
\cref{thm:biclique-induced-extraction}; thus $\ell=1000d^2$.
We may assume that $G$ is nonempty. Let $F$ be a nonempty induced
subgraph of $G$ with $\delta(F)=\dgn(G)$, and suppose that
$\dgn(G)\ge(s+b)^\ell$. Then $F$ has average degree at least
$(s+b)^{1000d^2}$. Since $F$ contains no $K_{s,s}$ as a subgraph, it
also contains no $K_{s+b,s+b}$ as a subgraph. Apply
\cref{thm:biclique-induced-extraction} with $s+b$ in place of $s$
(allowed, since $s+b\ge b(d)$) to obtain
an induced subgraph $F'$ of $F$ with no $C_4$ as a subgraph and with
average degree at least $2D$.
By \cref{lem:basic}\ref{item:heredity},
$\branchpi_H(F')\le k$.
But \cref{cor:biclique-c4-forest} gives $\dgn(F')\le D-1$, so
$\overline d(F')\le2\dgn(F')\le2(D-1)<2D$, a contradiction.
\end{proof}

Since $\branchpi_H(G)=0$ exactly when $G$ is $H$-free
(\cref{lem:basic}\ref{item:zero}), the case $k=0$ of
\cref{prop:biclique-polynomial-degree} recovers the theorem of Scott,
Seymour, and Spirkl~\cite{SSS23} that, for every forest
$H$, the degeneracy of $H$-free graphs with no $K_{s,s}$ subgraph is
bounded by a polynomial in~$s$.

\subsection{Induced bicliques and the chromatic bound}

We now convert exclusion of an induced biclique into exclusion of
a larger biclique as a subgraph, with a polynomial loss in the
clique number.

\begin{proof}[Proof of \cref{thm:forest-induced-biclique-chi}]
Take $b$ and $\ell=4000D_{h,k}^2$ from
\cref{prop:biclique-polynomial-degree}, and set $c=(b+1)^\ell$.
The empty graph satisfies the assertion, so assume $G$ is nonempty,
write $\omega=\omega(G)\ge1$, and put
\[
 t=\binom{\omega+s-1}{s-1}.
\]
By the bound of Erd\H{o}s and Szekeres~\cite{ES} on Ramsey numbers,
every graph with at least $t$ vertices contains a clique of size
$\omega+1$ or an independent set of size $s$.

If $G$ contained a $K_{t,t}$ as a subgraph with parts $A,B$,
then $G[A]$ and $G[B]$ would each contain an independent set
of size $s$, since neither contains a clique of size $\omega+1$.
These two independent sets would induce a $K_{s,s}$ in $G$,
a contradiction. Since $\omega\ge1$, we also have
$t\ge\binom{s}{s-1}=s\ge2$. Thus
\cref{prop:biclique-polynomial-degree}, applied with $t$ in place of
$s$, gives
\[
 \chi(G)\le\dgn(G)+1\le(t+b)^\ell.
\]
Finally,
$t=\binom{\omega+s-1}{s-1}\le(\omega+1)^{s-1}$:
the binomial coefficient counts the nondecreasing sequences of length
$s-1$ over a set of size $\omega+1$, while there are
$(\omega+1)^{s-1}$ sequences of length $s-1$ over such a set in
total. Moreover, $t+b\le(b+1)t$ because $t\ge1$. Therefore
\[
 (t+b)^\ell
 \le(b+1)^\ell t^\ell
 \le c(\omega+1)^{(s-1)\ell}.
\]
This proves the theorem.
\end{proof}

\begin{remark}
The assumption in \cref{thm:forest-induced-biclique-chi} that $H$ is
a forest cannot be dropped, even for $k=0$. Indeed, let $H$ be a
graph containing a cycle. By a theorem of Erd\H{o}s~\cite{ErdosGirth},
there are graphs of girth greater than $\max\{4,|V(H)|\}$ with
arbitrarily large chromatic number. Every such graph is
triangle-free and contains no $C_4$ as a subgraph. It therefore has
clique number at most two and no induced copy of $K_{s,s}$ for any $s\ge2$.
It is also $H$-free, since an induced copy of $H$ would contain a cycle of
length at most $|V(H)|$. Hence its induced $H$-packing branchwidth
is zero by \cref{lem:basic}\ref{item:zero}. Consequently, for every
$k\in\Nzero$ and $s\ge2$, the class of graphs with
$\branchpi_H\le k$ and no induced copy of $K_{s,s}$ is not $\chi$-bounded.
\end{remark}

\section{Dominated balanced separators}\label{sec:separators}

In this section, we prove the following theorem.
Its first assertion is exactly~\cref{thm:dominated}, so~\cref{thm:dominated}
follows immediately.
The second assertion is used in~\cref{sec:linegraphs}.

\begin{theorem}\label{thm:dominated-full}
Let $t\ge2$, let $k\ge1$, and let $G$ satisfy $\branchpi_{P_t}(G)\le k$.
For every vertex-weighting $\wei$ of $G$, the following hold.
\begin{enumerate}[label=(\arabic*)]
\item There is $D\subseteq V(G)$ with $|D|\le t\lfloor3k/2\rfloor+t-2$
such that $N[D]$ is a $1/2$-balanced separator of $(G,\wei)$.\label{item:half-sep}
\item There is $D\subseteq V(G)$ with $|D|\le tk+t-2$ such that
$N[D]$ is a $2/3$-balanced separator of $(G,\wei)$ and $G[D]$ has at most $k+1$
components, and at most $k$ components if $t=2$.\label{item:twothirds-sep}
\end{enumerate}
\end{theorem}

The next lemma finds a dominated balanced separator in a graph in which
no induced copy of $P_t$ meets two parts of a partition into light parts. Its
proof is a variant of the path argument of
Gy\'arf\'as~\cite{gyarfas1987problems}.

\begin{lemma}\label{lem:clean-partition-balanced}
Let $t\ge2$, $1/2\le\delta<1$, and $W\ge0$.
Let $F$ be a graph with a vertex-weighting $\wei$ satisfying
$\wei(V(F))\le W$, and let $A_1,\ldots,A_q$ be a partition of $V(F)$,
with some parts possibly empty, such that $\wei(A_i)\le\delta W$ for
every $i$.
Suppose that no induced copy of $P_t$ in $F$ meets two parts.
Then either every component of $F$ has weight at most $\delta W$,
or $t\ge3$ and $F$ contains an induced path $Q$ on at most $t-2$
vertices such that every component of $F-N[V(Q)]$ has weight
at most $\delta W$.
\end{lemma}

\begin{proof}
If $t=2$, there are no edges between distinct parts, so every component
of $F$ is contained in a part and has weight at most $\delta W$.
We may therefore assume that $t\ge3$ and that $F$ has a component
$C_0$ of weight greater than $\delta W$.
Such a component is unique, since $\delta\ge1/2$ and $\wei(V(F))\le W$.

Suppose for a contradiction that no induced path on at most $t-2$
vertices satisfies the conclusion.
Choose $x_1\in V(C_0)$.
We inductively construct an induced path $x_1,\ldots,x_{t-1}$.
For $1\le j\le t-2$, once $x_1,\ldots,x_j$ have been chosen, let
$C_j$ be the unique component of
$F-N[\{x_1,\ldots,x_j\}]$ of weight greater than $\delta W$.
Its existence follows from our assumption.
Deleting vertices cannot increase component weights, so $C_j$ is
a component of $C_{j-1}-N[x_j]$.

For $j=1$, we have $x_1\in V(C_0)$; for $j\ge2$, the preceding
choice ensures that $x_j$ has a neighbor in $C_{j-1}$.
In either case, connectedness of $C_{j-1}$ gives a vertex
$x_{j+1}\in N(x_j)\cap V(C_{j-1})$ with a neighbor in $C_j$.
Since $V(C_{j-1})$ is anticomplete to
$\{x_1,\ldots,x_{j-1}\}$, each extension preserves inducedness.

The resulting path $x_1,\ldots,x_{t-1}$ has its last vertex adjacent
to $C_{t-2}$, while $V(C_{t-2})$ is anticomplete to
$\{x_1,\ldots,x_{t-2}\}$.
Since $\wei(V(C_{t-2}))>\delta W$, this component meets at least two
parts. Choose $y\in V(C_{t-2})$ in a different part from $x_1$,
and extend $x_1,\ldots,x_{t-1}$ by a shortest $x_{t-1}$--$y$ path
in $F[V(C_{t-2})\cup\{x_{t-1}\}]$.
The extended path is induced, has at least $t$ vertices, and meets
two parts. A subpath on $t$ consecutive vertices containing an edge
between distinct parts is therefore a forbidden induced copy of $P_t$.
This contradiction proves the lemma.
\end{proof}

We can now prove the main result of this section.

\begin{proof}[Proof of~\cref{thm:dominated-full}]
Let $W=\wei(V(G))$.
If $W=0$, take $D=\varnothing$ for both assertions.
If some vertex $v$ has weight at least $W/2$, then $D=\{v\}$
also satisfies both assertions, as $k\ge1$.
We may therefore assume that $W>0$ and $\wei(v)<W/2$ for every
$v\in V(G)$; in particular, $G$ has at least three vertices.

Fix a branch decomposition $(T,\lambda)$ of $G$ with
$\branchpi_{P_t}(G,(T,\lambda))\le k$.
Put the weights of $G$ on the corresponding leaves of $T$ and
assign zero to internal nodes.
A centroid $z$ is internal, since deleting a leaf leaves a component
of weight greater than $W/2$.
The components of $T-z$ give a partition $A_1,\ldots,A_q$ of
$V(G)$, where $q\le3$ and $\wei(A_i)\le W/2$ for every $i$.

For~\ref{item:half-sep}, use this partition.
Every induced copy of $P_t$ meeting two parts crosses at least two of the
cuts displayed by the edges incident with $z$, and
$\pi^\times_{P_t}(G;A_e,B_e)\le k$ for each such edge $e$.
Counting these incidences shows that every induced $P_t$-packing
whose members meet two parts has at most $\lfloor3k/2\rfloor$
members.

For~\ref{item:twothirds-sep}, observe that the weights of the parts
sum to $W$, so some part, say $A_1$, has weight at least $W/3$.
Use the partition $\{A_1,V(G)\setminus A_1\}$. The corresponding edge
incident with $z$ displays the cut $(A_1,V(G)\setminus A_1)$, both
sides have weight at most $2W/3$, and
$\pi^\times_{P_t}(G;A_1,V(G)\setminus A_1)\le k$.

To complete either assertion, let $\delta=1/2$ for~\ref{item:half-sep}
and $\delta=2/3$ for~\ref{item:twothirds-sep}, and choose an
inclusion-maximal induced $P_t$-packing $\mathcal P$ whose members
meet two parts of the corresponding partition. Put
$F=G-N_G[V(\mathcal P)]$. The intersections of the parts with
$V(F)$ form a partition of $V(F)$, with some parts possibly empty,
into sets of weight at most $\delta W$, and
no induced copy of $P_t$ in $F$ meets two of them, since such a path could be
added to $\mathcal P$.

If every component of $F$ has weight at most $\delta W$, take
$D=V(\mathcal P)$.
Otherwise, \cref{lem:clean-partition-balanced}, applied with the
same $W$, gives $t\ge3$ and an induced path $Q$ on at most $t-2$
vertices such that every component of $F-N_F[V(Q)]$ has weight
at most $\delta W$.
Taking $D=V(\mathcal P)\cup V(Q)$ gives
$G-N_G[D]=F-N_F[V(Q)]$, so $N_G[D]$ is a $\delta$-balanced separator of $(G,\wei)$.

In both cases, $|D|\le t|\mathcal P|+t-2$.
Moreover, $D$ is the union of the vertex sets of the members of
$\mathcal P$, together with $V(Q)$ in the second case (which requires
$t\ge3$).
Each of these sets induces a connected subgraph, so $G[D]$ has at
most $|\mathcal P|+1$ components, and at most $|\mathcal P|$ if $t=2$.
The two bounds on $|\mathcal P|$ therefore give both assertions.
This finishes the proof.
\end{proof}
\section{Line graphs}\label{sec:linegraphs}

The goal of this section is to prove~\cref{thm:linegraph}.

\linegraphtheorem*

We first turn the second assertion of~\cref{thm:dominated-full} into a
bound on tree-independence number under star exclusion. We then
transfer a tree decomposition of $L(G)$ to a branch decomposition
of $G$.

The following lemma bounds the independence number of the closed
neighborhood of each connected set used to dominate our separators.

\begin{lemma}\label{lem:connected-neighborhood}
	Let $d\ge2$. If $Q$ is $K_{1,d}$-free and $S\subseteq V(Q)$ is nonempty
	and connected, then $\alpha(Q[N[S]])\le(d-2)|S|+1$.
\end{lemma}

\begin{proof}
	Order the vertices of $S$ as $v_1,\ldots,v_{|S|}$ so that every $v_i$, $i\ge2$, has an
	earlier neighbor $p_i$.
	Let $I\subseteq N[S]$ be independent, and assign each vertex of
	$I$ to the first closed neighborhood $N[v_i]$ containing it.
	Write $I_i$ for the vertices assigned to $v_i$.
	
	Since $Q$ is $K_{1,d}$-free, we have $|I_1|\le d-1$.
	For $i\ge2$, the set $I_i$ avoids $N[p_i]$, so $v_i\notin I_i$
	and $I_i\cup\{p_i\}$ is an independent subset of $N(v_i)$.
	Thus $|I_i|\le d-2$.
	Summing gives
	$|I|\le(d-1)+(|S|-1)(d-2)=(d-2)|S|+1$.
\end{proof}

The next lemma combines the preceding estimate with our dominated
balanced separators to bound tree-independence number under star
exclusion.

\begin{lemma}\label{lem:star-quantitative}
	Let $t,d\ge2$ and $k\ge1$. If $Q$ is $K_{1,d}$-free and
	$\branchpi_{P_t}(Q)\le k$, then
	$\treealpha(Q)\le7\bigl(((d-2)t+1)k+(d-2)(t-2)+1\bigr)$.
	For $t=2$, the stronger bound
	$\treealpha(Q)\le7(2d-3)k$ holds.
\end{lemma}

\begin{proof}
	Fix a normalized vertex-weighting $\wei$ of $Q$, and choose $D$ as in
	the second assertion of~\cref{thm:dominated-full}.
	Let $p$ be the number of components of $Q[D]$.
	Applying~\cref{lem:connected-neighborhood} to their vertex sets
	and using subadditivity gives
	$\alpha(Q[N[D]])\le(d-2)|D|+p$.

	Since $|D|\le tk+t-2$ and $p\le k+1$, we obtain
	$\alpha(Q[N[D]])
	\le((d-2)t+1)k+(d-2)(t-2)+1$.
	For $t=2$, we have $|D|\le2k$ and $p\le k$, so
	$\alpha(Q[N[D]])\le(2d-3)k$.
	
	In both cases, $N[D]$ is a $2/3$-balanced separator of $(Q,\wei)$.
	Since the weighting was arbitrary,~\cref{lem:sep-vs-ta} gives the corresponding bound on
	$\treealpha(Q)$.
	This proves the lemma.
\end{proof}

The next lemma allows us to bound the treewidth of $G$ using a tree
decomposition of $L(G)$, via the maximum matching width defined
in~\cref{sec:prelim}.

\begin{lemma}\label{lem:linegraph-transfer}
	For every graph $G$ with at least one edge,
	$\mmw(G)\le\treealpha(L(G))$. Consequently,
	$\tw(G)+1\le3\treealpha(L(G))$.
\end{lemma}

\begin{proof}
	We may discard isolated vertices of $G$, since this changes none
	of the parameters in the statement.
	Let $a=\treealpha(L(G))$, and take a tree decomposition $(T,\beta)$
	of $L(G)$ of independence number $a$.
	
	For each $v\in V(G)$, the edges incident with $v$ form a nonempty
	clique in $L(G)$ and hence are contained together in a bag
	$\beta(\eta(v))$ for some node $\eta(v)$.
	Since $G$ has at least two vertices, \cref{lem:tree-to-branch}
	gives a branch decomposition of $G$ in which every displayed cut
	$(A,B)$ has a node $x$ of $T$ lying on the path between $\eta(u)$
	and $\eta(v)$ for all $u\in A$ and $v\in B$.
	If $uv\in E(G)$ crosses $(A,B)$, then the vertex $uv$ of $L(G)$
	belongs to $\beta(\eta(u))$ and $\beta(\eta(v))$.
	The nodes whose bags contain $uv$ therefore induce a subtree that
	contains the path between $\eta(u)$ and $\eta(v)$, and in particular
	contains $x$.
	Thus all edges of $G$ crossing $(A,B)$ belong to the single bag
	$\beta(x)$.
	A matching among them is an independent set in
	$L(G)[\beta(x)]$, so it has size at most $a$.
	
	Since the cut was arbitrary, $\mmw(G)\le a=\treealpha(L(G))$.
	The stated treewidth bound follows from \cref{lem:mmw-tw}.
\end{proof}

We can now put everything together.

\begin{proof}[Proof of~\cref{thm:linegraph}]
	By \cref{lem:basic}\ref{item:comparison},
	$k\le\treepi_{P_t}(L(G))$. Moreover, $\treepi_{P_t}\le\treealpha$:
	for every bag $\beta(x)$ of a tree decomposition, choosing one
	vertex of $\beta(x)$ from each member of a $\beta(x)$-rooted
	induced $P_t$-packing gives an independent set, since the members
	are pairwise anticomplete. Dallard et al.~\cite[Theorem~5.3]{DKKMMW}
	proved that $\treealpha(L(G))\le\tw(G)+1$. Hence
	$k\le\tw(G)+1$.
	
	We may assume that $G$ is nonempty.
	Suppose first that $k=0$.
	Then $L(G)$ is $P_t$-free by~\cref{lem:basic}\ref{item:zero}, so
	$G$ contains no path of
	$t$ edges: the edges of such a path would induce a $P_t$ in $L(G)$.
	Take a depth-first-search forest of $G$.
	Every root-to-leaf path has at most $t$ vertices, and every edge
	of $G$ has ancestor--descendant endpoints.
	At each vertex, take the bag consisting of that vertex and all
	its ancestors.
	Every edge is covered in the bag of its deeper endpoint, and a
	vertex occurs precisely in the bags indexed by itself and its
	descendants.
	Thus these bags, after joining the decomposition trees for distinct
	components, form a tree decomposition of $G$ with bags of size
	at most $t$.
	Hence $\tw(G)+1\le t$, which also gives the general upper bound
	since $t\le21(t-1)$.
	
	Now assume that $k\ge1$.
	The graph $L(G)$ is $K_{1,3}$-free, since the neighbors of a vertex
	corresponding to an edge $uv$ are covered by the two cliques of
	edges incident with $u$ and with $v$.
	Applying~\cref{lem:star-quantitative} with $d=3$ and then
	\cref{lem:linegraph-transfer} gives
	$\tw(G)+1\le3\treealpha(L(G))
	\le21((t+1)k+t-1)$.
	For $t=2$, the sharper bound in~\cref{lem:star-quantitative}
	instead gives
	$\tw(G)+1\le3\treealpha(L(G))\le63k$.
	This proves all assertions.
\end{proof}

The next lemma shows that the linear dependence on $t$
in~\cref{thm:linegraph} is necessary up to an absolute factor.

\begin{lemma}\label{lem:sharpness}
	Let $t\ge2$ and $c>0$. If
	$\tw(G)+1\le c\,(1+\branchpi_{P_t}(L(G)))$ for every graph $G$,
	then $c\ge t+1$.
\end{lemma}

\begin{proof}
	If edges $e_1,\ldots,e_t$ of a graph induce a $P_t$ in its line graph,
	then consecutive edges meet and nonconsecutive edges are disjoint.
	Moreover, each internal edge meets its two neighbors at distinct
	endpoints, since otherwise those neighbors would be adjacent in
	the line graph.
	Thus $e_1,\ldots,e_t$ form a path on $t+1$ vertices.

	Different members of an induced packing in a line graph thus yield
	vertex-disjoint paths, because edges belonging to different
	members share no endpoint.
	Consequently,
	\[
	 \branchpi_{P_t}(L(K_n))\le\lfloor n/(t+1)\rfloor.
	\]
	Since $\tw(K_n)+1=n$, taking $n=q(t+1)$ gives
	$q(t+1)\le c\,(q+1)$.
	Letting $q$ tend to infinity proves $c\ge t+1$.
\end{proof}

By~\cref{thm:linegraph}, one may take $c=21(t+1)$ in
\cref{lem:sharpness}. Thus the dependence on $t$ is linear up to an
absolute factor. The same comparison holds for induced packing
treewidth: for nonempty $G$, the inequalities
$\branchpi_{P_t}(L(G))\le\treepi_{P_t}(L(G))
\le\treealpha(L(G))\le\tw(G)+1$
show that $1+\treepi_{P_t}(L(G))=\Theta_t(1+\tw(G))$.
\section{Metric properties}\label{sec:metric}

In this section, we prove~\cref{thm:qisometri}.

\quasiso*

\subsection{Coarse treewidth and quasi-isometries}

For integers $a\ge1$ and $R\ge0$, a set $X\subseteq V(G)$ is
\emph{$(a,R)$-coverable} if $X\subseteq B(Z,R)$ for some
$Z\subseteq V(G)$ with $|Z|\le a$. The centers in $Z$ need not
belong to $X$. A class has \emph{bounded coarse treewidth} if there
exist $a,R$ such that every graph in the class has a tree decomposition
with $(a,R)$-coverable bags. Equivalently, there are fixed $a,D\ge1$ such
that every bag can be partitioned into at most $a$ sets of weak diameter
at most $D$: a cover by radius-$R$ balls gives $D=\max\{2R,1\}$, and a partition
into sets of weak diameter at most $D$ gives a cover by radius-$D$ balls.

For $L\ge1$ and $C\ge0$, an \emph{$(L,C)$-quasi-isometry} from a
graph $G$ to a graph $Q$ is a map $\varphi:V(G)\to V(Q)$ such that
\[
\begin{aligned}
	\dist_Q(\varphi(u),\varphi(v))&\le L\cdot\dist_G(u,v)+C,\\
	\dist_G(u,v)&\le L\cdot\dist_Q(\varphi(u),\varphi(v))+C
\end{aligned}
\qquad(u,v\in V(G)),
\]
and $V(Q)=B_Q(\varphi(V(G)),C)$. The inequalities are interpreted
with extended distances. A class is \emph{uniformly quasi-isometric
	to graphs of bounded treewidth} if there exist $L,C,b$ such that every
member admits an $(L,C)$-quasi-isometry to a graph of treewidth at most
$b$.

We use the following characterization, which identifies bounded coarse
treewidth with uniform quasi-isometry to graphs of bounded treewidth.
Its first part is due to Nguyen, Scott, and
Seymour~\cite[Theorem~1.2]{NSSCoarse}. Hickingbotham obtained it
independently, with a weaker bound on the treewidth, together with the
converse~\cite[Theorem~1 and Lemma~6]{Hick};
see also~\cite[Theorem~2.1]{NSSCoarse}.

\begin{lemma}[Nguyen, Scott, and Seymour~\cite{NSSCoarse};
Hickingbotham~\cite{Hick}]\label{lem:NSS-coarse-tw}
For all integers $a,D\ge1$, there exist $L,C\ge1$ such that if $G$ has a tree
decomposition in which every bag is the union of at most $a$ sets of
weak diameter at most $D$, then $G$ admits an $(L,C)$-quasi-isometry
to a graph of treewidth at most $a$. Conversely, for all $L\ge1$,
$C\ge0$, and $b\ge0$, if $G$ admits an
$(L,C)$-quasi-isometry to a graph of treewidth at most $b$, then $G$
has a tree decomposition in which every bag is the union of at most
$b+1$ sets of weak diameter at most $3c^2$, where
$c=\lceil\max\{L,C\}\rceil$.
\end{lemma}

In the first part of~\cref{lem:NSS-coarse-tw}, one may take
\begin{equation}\label{eq:metric-centred-conversion}
	L=D(2a+2),\qquad C=D(2a-1).
\end{equation}
To see this, we use a notion of Nguyen, Scott, and Seymour~\cite{NSSCoarse}.
A \emph{pseudo-tree-decomposition} of a graph $Q$ is defined like a tree
decomposition $(T,\beta)$, except that an edge $uv$ of $Q$ may instead be
covered by an edge $xy$ of $T$ with $\beta(x)\setminus\beta(y)=\{u\}$ and
$\beta(y)\setminus\beta(x)=\{v\}$. Its width is $\max_{x\in V(T)}|\beta(x)|$,
and the \emph{pseudo-tree-width} of $Q$ is the minimum width of such a
decomposition.

First let $G$ be connected. Turning each of the sets that form the bags
into a clique gives a graph $G'$ with the same tree decomposition, in which every bag
has independence number at most $a$, and the identity map is a
$(D,0)$-quasi-isometry from $G$ to $G'$. By~\cite[Theorem~2.3]{NSSCoarse},
$G'$ admits a $(2a+2,2a-1)$-quasi-isometry to a graph $Q$ of
pseudo-tree-width at most $a$; composing the two maps
gives~\eqref{eq:metric-centred-conversion}, as in the proof
of~\cite[Theorem~1.5]{NSSCoarse}. For disconnected $G$, take the disjoint
union of the graphs $Q$ obtained for the components. Finally, $\tw(Q)\le a$:
in a pseudo-tree-decomposition $(T_Q,\beta_Q)$ of $Q$ of width at most $a$,
subdivide every edge $xy$ of $T_Q$ that covers an edge of $Q$ in the pseudo
sense, and give the new node the bag $\beta_Q(x)\cup\beta_Q(y)$, which has
$|\beta_Q(x)|+1\le a+1$ vertices.

\subsection{Distanced packings}

Let $H$ have at least two vertices, and let $r\ge1$ be an integer.
For a cut $(A,B)$, let
$\pi^\times_{H,r}(G;A,B)$ be the maximum number of induced
copies of $H$ crossing the cut whose vertex sets have pairwise distance
greater than $r$. Define
\[
\branchpi_{H,r}(G)
=\min_{\mathcal D=(T,\lambda)}\;
\max_{e\in E(T)}
\pi^\times_{H,r}(G;A_e,B_e).
\]
The value is zero when $|V(G)|\le1$. Since vertex sets at distance
greater than $1$ are disjoint and anticomplete,
$\branchpi_{H,1}=\branchpi_H$. For $r>1$, this parameter need not be
monotone under taking induced subgraphs, since distances can grow.
For integers $t\ge2$, $r\ge1$, and $k\ge0$, let
\[
a_k=\max\left\{1,\left\lfloor\frac{3k}{2}\right\rfloor\right\},
\qquad D_{t,r}=2r+3t-3.
\]

We are now ready to prove the main result of this section.

\begin{proof}[Proof of~\cref{thm:qisometri}]
	The equivalence of the first two assertions follows
	from~\cref{lem:NSS-coarse-tw}. We prove the equivalence of the
	first and third assertions, retaining quantitative bounds.
	We first show that
	\begin{equation}\label{eq:metric-centred-decomposition}
		\branchpi_{P_t,r}(G)\le k
		\quad\Longrightarrow\quad
		\begin{gathered}
			G\text{ has a tree decomposition whose bags are unions of}\\
			\text{at most }a_k\text{ sets of weak diameter at most }D_{t,r}.
		\end{gathered}
	\end{equation}

	Suppose that $G$ is connected and $|V(G)|\ge2$, and fix a branch
	decomposition $(T,\lambda)$ witnessing the bound. For each $v\in V(G)$,
	let $T_v$ be the minimal subtree of $T$ containing $\lambda(N[v])$,
	and let
	\[
	\beta(x)=\{v\in V(G)\mid x\in V(T_v)\}\qquad(x\in V(T)).
	\]
	The nodes whose bags contain $v$ are exactly the nodes of $T_v$, and,
	for $uv\in E(G)$, $\lambda(u)\in V(T_u)\cap V(T_v)$. Thus $(T,\beta)$
	is a tree decomposition. Moreover,
	\[
	\beta(\lambda(u))=N[u],
	\]
	so every leaf bag has weak diameter at most two.

	Fix an internal node $x$ of $T$, and let $V_1,V_2,V_3$ be the vertex
	classes corresponding to the three components of $T-x$. Choose an
	inclusion-wise maximal family $\mathcal P$ of induced copies of $P_t$ meeting
	at least two classes and having pairwise distance greater than $r$.
	Every member crosses at least two of the three cuts displayed by the
	edges incident with $x$, so
	\begin{equation}\label{eq:metric-node-count}
		2|\mathcal P|
		\le\sum_{i=1}^{3}
		\pi^\times_{P_t,r}(G;V_i,V(G)\setminus V_i)
		\le3k.
	\end{equation}
	Let $J=G-B(V(\mathcal P),r)$.
	By maximality, no induced copy of $P_t$ in $J$ meets two classes.

	Suppose first that $\mathcal P\ne\varnothing$. We claim that
	$\beta(x)\subseteq B(V(\mathcal P),r+t-1)$. Let $v\in\beta(x)\cap V_i$.
	Since $x\in V(T_v)$, the set $N_G[v]$ meets at least two classes.
	Thus $v$ has a neighbor $w\notin V_i$. If $v$ or $w$ belongs to
	$B_G(V(\mathcal P),r)$, then $\dist_G(v,V(\mathcal P))\le r+1$.
	Otherwise $v,w$ belong to one component $K$ of $J$. The cut
	\[
	\bigl(V(K)\cap V_i,\ V(K)\setminus V_i\bigr)
	\]
	is nontrivial and $P_t$-clean. This is impossible for $t=2$; for
	$t\ge3$, \cref{lem:clean-cut} gives
	$\dist_K(v,z)\le t-2$ for every $z\in V(K)$.
	Since $G$ is connected, some $z\in V(K)$ satisfies
	$\dist_G(z,V(\mathcal P))=r+1$. Hence
	\[
	\dist_G(v,V(\mathcal P))
	\le\dist_K(v,z)+\dist_G(z,V(\mathcal P))\le r+t-1,
	\]
	which proves the claim.

	Assign each vertex of $\beta(x)$ to a member $P\in\mathcal P$
	at distance at most $r+t-1$. Vertices assigned to the same path
	have distance at most
	\[
	2(r+t-1)+(t-1)=D_{t,r}.
	\]
	By~\eqref{eq:metric-node-count}, this partitions $\beta(x)$ into
	at most $a_k$ sets of the required weak diameter.

	If $\mathcal P=\varnothing$, no induced copy of $P_t$ crosses any of the three
	cuts displayed by the edges incident with $x$. Connectedness gives an
	edge between two classes. Thus $t\ge3$, and \cref{lem:clean-cut} gives
	$\diam(G)\le2t-4\le D_{t,r}$.
	In this case one part suffices.

	For disconnected $G$, take the restriction of the branch decomposition
	to each component with at least two vertices. Distances within a
	component are unchanged, and every cut displayed by a restriction is
	the restriction of a cut displayed by the original decomposition.
	Apply the preceding
	construction componentwise, use a singleton bag for each isolated
	vertex, and join the decomposition trees. Graphs with at most one
	vertex are immediate.
	This proves~\eqref{eq:metric-centred-decomposition}, and hence the third
	assertion implies the first.

	Conversely, graphs with at most one vertex are immediate. Suppose
	$(T,\beta)$ is a tree decomposition of $G$ with $(a,R)$-coverable
	bags, increasing $R$ to at least one if necessary.
	For every vertex $v$, let $\eta(v)$ be a node whose bag contains $v$.
	By \cref{lem:tree-to-branch}, $G$ has a branch decomposition in which,
	for every displayed cut, some bag of $(T,\beta)$ meets every connected
	subgraph of $G$ crossing the cut.

	Consider a family of induced copies of $P_t$ that cross a displayed cut and
	whose vertex sets have pairwise distance greater than $2R$. Fix a bag
	that meets every connected subgraph crossing this cut, and choose
	one vertex from its intersection with each path. Two chosen vertices
	cannot belong to the same radius-$R$ covering ball. Thus the family
	has at most $a$ members, and
	\[
	\branchpi_{P_t,2R}(G)\le a.
	\]
	This proves the remaining implication, and hence finishes the proof.
\end{proof}

\begin{corollary}\label{thm:metric-quasi-isometry}
	Let $t\ge2$, $r\ge1$, and $k\ge0$ be integers.
	If $\branchpi_{P_t,r}(G)\le k$, then $G$ admits an
	$(L_{t,r,k},C_{t,r,k})$-quasi-isometry to a graph of treewidth at most
	$a_k$, where
	\[
	L_{t,r,k}=(2a_k+2)D_{t,r},
	\qquad C_{t,r,k}=(2a_k-1)D_{t,r}.
	\]
	In particular, if $\branchpi_{P_t}(G)\le k$, one may take
	\[
	L_{t,1,k}=(2a_k+2)(3t-1),
	\qquad C_{t,1,k}=(2a_k-1)(3t-1).
	\]
\end{corollary}

\begin{proof}
	Apply~\eqref{eq:metric-centred-conversion} to the decomposition
	in~\eqref{eq:metric-centred-decomposition}.
\end{proof}

For $t=2$, the identity $\branchpi_{P_2}=\simw$ shows that
\cref{thm:metric-quasi-isometry} extends
\cite[Theorem~14]{Hick}, with target treewidth at most $a_k$ in place
of $12k-1$ for $k\ge1$. The hypothesis is strictly more general for
every $t\ge3$. Indeed, Munaro and Yang~\cite[Proposition~3 and
Theorem~9]{munaro2023algorithmic} imply that $(K_5,P_3+P_2+P_1)$-free graphs have unbounded
sim-width. An induced $P_3$-packing of size three contains an induced
copy of $P_3+P_2+P_1$, so every graph in this class satisfies
$\branchpi_{P_3}\le2$, and hence $\branchpi_{P_t}\le2$ for every
$t\ge3$ by \cref{lem:basic}\ref{item:zero}.

\bibliographystyle{abbrvurl}
\bibliography{ref}

\end{document}